\documentclass[11pt]{article}
\usepackage[textwidth=14.5cm,textheight=22cm]{geometry}

\usepackage{amsmath,amssymb}
\usepackage{color}
\usepackage{graphicx}
\usepackage[sectionbib]{chapterbib}
\usepackage{comment}
\usepackage{tikz}
\usetikzlibrary{shapes}
\usetikzlibrary{decorations.markings}
\usetikzlibrary{arrows.meta}
\usetikzlibrary{calc}
\usepackage{float}

\usepackage{etoolbox}
\usepackage{eurosym}
\usepackage{xurl}

\usepackage{hyperref} 
\hypersetup{colorlinks=true, linkcolor=blue,  anchorcolor=blue,  
citecolor=blue, filecolor=blue, menucolor=blue,   
urlcolor=blue}

\title{A Primer on Ethics in Mathematics}
\author{Maurice Chiodo\footnote{Email: mcc56@cam.ac.uk} \\University of Cambridge}

\begin{document}

\maketitle

\begin{abstract}
This is a course about \textit{Ethics in Mathematics}. That is, about the ethical issues that arise in mathematics itself, in mathematical work, and in the actions of mathematicians, broadly interpreted to cover all mathematically-powered technologies.

This is a \textit{novel course}. Traditionally, ethics has seldom been considered by, let alone taught within, the mathematics discipline. This course seeks to change that, by bringing to light the specific ethical issues relevant to mathematics, in a way that is clear and understandable to mathematicians.

This course is written \textit{by} mathematicians, \textit{for} mathematicians, and \textit{about} mathematicians. It presents new and original material, in the form of case studies covering harmful and unethical behaviour, analysis of the causal chain of events leading to such harm, and examination of the surrounding processes that drive such actions. All this is fully cited, and supported with interactive activities to aid learning. It develops much of the foundational \textit{theory} of ethics in mathematics: what can go wrong, why it goes wrong, and how to address it. Students and practitioners alike take away an understanding of the context of their mathematical work, and how to change their actions to ensure more appropriate and useful outputs.

In contrast to more traditional philosophy teaching, this course does not delve into the depths of ethical theory. Rather, it presents a practical, hands-on, pragmatic approach to ethics in mathematics, treating it as a branch of applied ethics. It uncovers the very particular ethical issues that arise from mathematical technologies, as well as show how general principles from ethics and responsible development apply even to  mathematical work. As such, this makes it highly relevant to adjacent areas such as computer science, physics, engineering, and any other numerate discipline.

\end{abstract}

\newpage
\tableofcontents
\newpage
\setcounter{section}{-1} 
\section{Introduction}

\fbox{
  \parbox{0.98\textwidth}{
This document -- a set of course notes -- lays out an ongoing course on Ethics in Mathematics at the University of Cambridge, delivered annually since 2016, and in its current form since 2018. Well-received by mathematics students over many years, it is one of the few examples of enduring ethics teaching to mathematicians.

This introduction describes who was involved in the course and how it came to be, how it was structured and delivered, what motivated the structure it took, and how to make use of these notes to create other similar courses.
}}

\subsection{How this course came about}

Everyone always asks me how this all started. How I came to give a course on ethics in mathematics, out of the blue. So I'll put it here, once and for all. This is the only part of the notes where I'll use the pronoun ``I'', because it is a very personal tale.

I had, for many years, become more and more concerned about how mathematicians and their work were impacting society for the worse. This sentiment started when I was an undergraduate, in particular when I was taught the governing equation for nuclear diffusion and calculated the minimum radius for a bomb from fissile material; at the time I felt it strange that this was presented with no further discussion about the consequences of such knowledge. As I progressed in my career, I saw the continued harmful impact of mathematical work around me. I saw the financial crisis of 2008 unfold. I saw the revelations of state-run mass surveillance systems in 2013. I saw the rise of `Big Tech' and the consequences it was beginning to have. And all the while, I knew that mathematicians, and their work, were powering the very sources of these concerns.

It was October 2015, and I had just started as a postdoc in the mathematics faculty at Cambridge. I was delivering a new third-year course, ``Automata and Formal Languages'', covering the fundamentals of the theory of computation. It was a few weeks into the course, and I found myself lecturing on the day of the Careers Event being run for students in the faculty, whereby various employers set up stalls and recruited students. So, at the end of the lecture, having made clear the `official' part of the lecture was over, I told the students I had a few words to say about the careers event, but that they were of course free to leave. None left. So I gave a 9-minute\footnote{I know this because the lecture ended at 12:55, and I distinctly recall being in the faculty cafeteria queue at 13:05 to ensure I got \textbf{something} for lunch before the food ran out.} explanation about how, as mathematicians, we create tools of great power, and how, when they go to the careers event, they should think carefully about in whose hands they might be placing such tools. I also spoke about ways in which employers can try to pressure and manipulate employees to act against their own better judgement of what is right and wrong, including how to deal with the argument ``If you don't do it, someone else will''. 9 minutes. That's it.

Then some unexpected things started to happen. I received an applause at the end of the 9 minutes. This was a mid-course applause; the first, and to date only, time this has happened in my career. I had students come up to thank me directly after the course. I had students email me to thank me and say ``that was the first time I've heard ethics seriously discussed by any member of the maths faculty'', even though I never explicitly used the word `ethics' in my 9 minutes. And I had other faculty members come up to me the following week asking me about the ``25 minute talk'' I gave on ethics that their students had told them about; the content was so new that the students thought it went on much longer than it actually did.

And so I saw a gap. Not just in the students' understanding, but also in the teaching they were receiving. I realised that they didn't know about the ethics of their mathematics, that no-one was telling them about it, and that they had no way of finding out about it. And this made me start to realise \textit{why} the mathematical disasters I had been observing had been occurring, and moreover \textit{would keep occurring} unless someone did something about it. I tried to raise this concern with colleagues, suggesting that perhaps the faculty could somehow create some sort of seminar series on this, but (apart from those listed below) received practically zero interest. And that's when I realised that \textit{if I don't do it, no-one else will}. So I found someone who \textit{was} interested in the idea, we discussed it extensively, and on 22 April 2016 I gave my first lecture on ``Ethics in Mathematics''. The rest, as they say, is history\ldots

\subsection{The formation of this course and who was involved}\label{s:who-was}

This course on ethics in mathematics was put together and first delivered by the author, Dr Maurice Chiodo\footnote{\url{https://www.mauricechiodo.com}}, at the University of Cambridge in 2016. Maurice is currently an Assistant Research Professor at the Centre for the Study of Existential Risk (CSER),  University of Cambridge\footnote{\url{https://www.cser.ac.uk/team/maurice-chiodo}}. Prior to that, he was a postdoctoral researcher in the Faculty of Mathematics at Cambridge (2015--2019), and then Full Fellow and College Teaching Officer in mathematics at King's College Cambridge (2019--2023). A research mathematician by training, and with 2 PhDs in mathematics, he initiated this course while he was a postdoctoral researcher in the Cambridge mathematics faculty, with the aim of raising ethical awareness among mathematicians; something he saw as somewhat lacking at the time.

This was not a lone endeavour. Initially, Maurice worked with Dr Piers Bursill-Hall from 2016--2019, discussing, exploring, and refining, many of the ideas that made up the content of the course. Piers was the history of mathematics lecturer at Cambridge for over 40 years. After Maurice delivered some initial ethics in mathematics lectures in 2016, the Cambridge University Ethics in Mathematics Society\footnote{\url{https://cueims.soc.srcf.net}} (CUEiMS) was formed by students. It a student society, run by students to host events relating to Ethics in Mathematics which might be of interest or value to students, and operates  completely independent of the University of Cambridge.

In 2017, Dennis Müller began contributing to the ideas within the lectures. Dennis was a mature-age undergraduate mathematics student at Cambridge from 2017 to 2020. He served various functions in CUEiMS, including acting and vice-president in the academic year 2017--2018, and president in 2018--2019. In the summer of 2018, Maurice, Piers, and Dennis hosted 9 summer students, to help look up additional case studies related to ethics in mathematics for the lectures. It was at this time that the course was refined from the original 2016 version, to the revised 2018 version still in use and presented in these notes. It was also at this time that the Ethics in Mathematics Project\footnote{\url{https://www.ethics-in-mathematics.com}} (EiMP) was formed. It is an academic project, focusing on producing academic output related to Ethics in Mathematics, such as research papers, notices and articles in various magazines, tools for those wishing to operationalise ethics in their mathematical work, and resources for those wishing to teach Ethics in Mathematics in some form.

\subsection{The structure and delivery of the course}

The Faculty of Mathematics at the University of Cambridge has a longstanding tradition of optional and non-examinable courses at the undergraduate and graduate levels. In the past, these have included several introductory history and philosophy of mathematics courses, as well as more specialised courses on physics and its philosophy. These courses have been typically taught by researchers from within the faculty. The Ethics in Mathematics course presented in these notes has been one of these optional non-examinable courses since 2016. Unlike some of the other courses, this course has been primarily hosted by CUEiMS, with additional help from EiMP. Throughout the years, the 8 lectures described here have been complemented by several guest speaker series organised jointly by CUEiMS and EiMP (see also \cite{müller2022situating} for further discussion of Ethics in Mathematics at Cambridge). Being hosted under the auspices of a student society, the faculty and its teaching board had no influence on the content and structure of the ethics in mathematics course, and in return, the course is not an official part of any degree programme in the faculty. Maurice Chiodo designed and delivered all lectures in the course, wrote this set of notes, and was fully responsible for all content covered.

Some of the choices made when teaching this course (e.g., why there are only 8 lectures) are owed to this fact. Cambridge terms are only 8 weeks, and the course had to fit well with the other non-examinable courses and students' already very busy schedules. These constraints also explain the limited reading and homework that was given to students: every student attending the seminar did so voluntarily, at the expense of their own time. When the course was taught, each lecture lasted between 60--90 minutes, with some informal discussion afterwards so the the whole session ran for about 2 hours.

These notes focus on the course as it was taught in the first terms of 2018, 2019, 2020 (online-covid), 2021 (online-covid), 2022, 2023, 2024, 2025. The course was mostly unchanged in those years. An initial 8 lecture course was given in 2016 and 2017. It was very different in structure (though with some content overlap), to the current incarnation of the course which commenced in 2018.

The course was titled ``Ethics for the Working Mathematician'', and was advertised to all students enrolled in the mathematics undergraduate and masters degrees at Cambridge, done via an email on the first day of the teaching term. The course was also listed on the CUEiMS website, and posted on ``talks.cam'' -- the university's repository website for all open talks, seminars, etc.~happening across the university. Attendance was open to all, and as such the audience was not exclusively mathematics undergraduates or masters students, though they did make up the vast majority of attendees. There was no requirement to join CUEiMS, or to pay any fees, to attend.  No part of the lectures were assessable in any way, nor could any students claim coursework credit for attending. 

Attendees were provided with a syllabus for the course, which contained a (non-compulsory) reading list -- it was not assumed that students had read any material from this list before attending. This was made available on the CUEiMS website, the most recent version being in 2025\footnote{\url{https://cueims.soc.srcf.net/2025/intro}}. There were no printed or digital notes provided for students, nor any homework, exercises, or other tasks set for students between lectures. However, a recording of the core content (but not of the post-lecture discussion) was made of the 2018 lectures, and can be found on YouTube\footnote{Playlist at \url{www.youtube.com/playlist?list=PLhjUVXp9e11Qu0mw6HlsIu_RJpT6N6N-J}}. These were not `traditional' mathematics lectures, and there was no mathematics taught in them. However, there was an assumption that attendees had at least started an undergraduate degree in mathematics, and mathematical concepts were occasionally discussed in a technical way. As such, attendees from outside mathematics would have come across various aspects of the lectures for which they may not have understood all the technical details.

\subsection{Motivation for this course design}

The primary audience for this course was undergraduate and masters mathematics students at Cambridge. As such, it was taken as a given that these students had an extremely firm grasp of abstract mathematical concepts. Wile the course referred to mathematics at various places, no mathematical theory was `taught'. Instead, the aim of the course was to raise ethical awareness among mathematicians. Indeed, this required a long a thorough treatment in the first 4 lectures that there even exists ethical questions associated with mathematics; this was not a concept that many mathematics students had encountered before.

And so, this was not an ethics course, at least not in the way that a philosopher might expect. It did not cover any theories of ethics. Deep problems of, say, \textit{how} to choose an optimisation framework and \textit{whose} interests to take into account, or \textit{which} mathematical problems should a mathematician refuse to work on, were not `solved' in this course. Rather, the course was designed to show that such problems and choices exist and need to be made, and that the solution to these is not always obvious and can sometimes be very hard. This was an introductory applied ethics course, of the most fundamental type: what \textit{are} the ethical questions that need to be addressed?

But such a course is still immensely valuable, even if if does not cover the frontier of ethics research. Bringing even an elementary understanding of ethics to mathematicians is just as valuable as, say,  bringing an elementary understanding of calculus and the trapezoidal rule to medical clinicians\footnote{Which was done in 1994 with a paper titled ``A Mathematical Model for the Determination of Total Area Under Glucose Tolerance and Other Metabolic Curves'' \cite{mary1994area}, and which currently has over 400 citations \cite{web-science}.}. And there was a concern that `cutting edge' ethics training might simply be off-putting for students who had never come across such concepts before; one would not try to teach multi-dimensional integral calculus to medical clinicians who had never heard of the trapezoidal rule.

With this, questions of the form ``What are the axioms for ethics?'', or indeed ``Whose ethics?'' where then easily deferred back to the students; now that they knew such ethical questions existed in mathematics, they could go and investigate different ethical theories and different frameworks for considering such questions, and combine those with their own world-view, to start deciding how \textit{they} wanted to do ethics. They understood mathematics extremely well, and the course showed them how their mathematics fits into (and potentially harms) the world and why this happens, thus they could then go away and consider how they wanted to change their own behaviour. The course was set up to connect their \textit{mathematical understanding} with their own  \textit{humanity}; two parts of them that, until now, they had perhaps not considered linked.

\subsection{How to use these notes}

These notes are divided into the 8 lectures as given in the course, and a read of the table of contents shows the breakdown of each lecture into various parts. Most lectures, and especially lectures 2--4, involve many `examples', signposted as such. These are short write-ups of an event or circumstance, written as one would present in a research paper, and heavily referenced. Analysis of how they fit into the broader picture is kept to a minimum there, and they are primarily written in a way to describe `what happened'. Alongside these, there are `activities' and `interactive problems', also signposted as such. These describe short activities run during lectures that had student involvement in some way, usually done to demonstrate a point about how mathematicians consider their work and the ethical questions that may arise from it.

The remainder of the notes is the overall analysis of how the examples, activities, and interactive problems, fit together to illustrate where ethical questions lie in mathematics, and what the sources are of such considerations. This analysis takes `what happened', and gives reasons for `why it happened', as well as guidance on `how to improve it or stop it happening'. The proportion of each lecture devoted to analysis increases in the latter lectures; the earlier ones are mostly devoted to presenting examples. 

These notes are highly modularised, with each example, activity, interactive problem, and piece of analysis presented in its own self-contained, labelled unit. The reason for presenting it like this is to enable someone creating their own course on ethics in mathematics to simply pick out and rearranges the pieces from this course, to produce a course relevant to their own circumstances and surrounds. This course was designed to be delivered to an audience of Cambridge undergraduate mathematicians; a different cohort, at a different university, in a different country, will no doubt require a different course. However, given the various common features of mathematical practice around the world, there will (hopefully) be some, if not many, parts of this course that can be used within ethics in mathematics courses elsewhere. Ethics is mathematics is best-addressed when considered at a local scale \cite{muller2025ethical}, and so any such course should have local concerns in mind \cite[p.~41]{muller2022hippocratic}.

And finally, given the large number of people involved in hosting the course, throughout the notes the pronoun ``we'' is used when describing what was done, even though the lecture content was primarily produced and delivered by Maurice Chiodo.

\subsection{Acknowledgements}

These notes would not have been possible without the thoughtful input of Dennis Müller. Throughout the writing of these notes, Dennis has provided invaluable commentary and suggestions for improvements. This is alongside his very valuable discussions during, and about, the lectures in 2017-2020. In addition, Dennis is a core contributor to the EiMP, having helped guide and carry out its outputs from 2017 to the present day.

Thanks also go to Piers Bursill-Hall, who helped form and shape all the early ideas that went into creating this course. Without him, this course would have never been formed.

Additional thanks go to Zoya Yasmine, for suggesting changes and additional references to Lecture 5 on law and regulation. 

And finally, thanks go to the 9 summer students of 2018, who helped look up additional case studies when the course was refined from the original 2016 version, to the revised 2018 version still in use and presented in these notes.

\begingroup
\let\section\subsection
\renewcommand{\bibname}{Bibliography (Introduction)}
\renewcommand{\refname}{Bibliography (Introduction)}
\addcontentsline{toc}{subsection}{Bibliography (Introduction)}
\renewcommand{\addcontentsline}[3]{}
\bibliographystyle{IEEEtran}
\bibliography{bibintro}
\endgroup

\newpage
\section{Lecture 1: An introduction to ethics in mathematics and why it is important}\label{L:1}

\fbox{
  \parbox{0.98\textwidth}{
Mathematicians sit at the heart of technological advancement and industrial progress. Mathematics is a universal tool. It is not neutral, and can be used help, and also to harm. But the world is complicated, and simple notions of `good' and `bad' mask the complexity that happens every time someone does mathematics.

This lecture looks at how mathematicians currently perceive what `ethics in mathematics' is, and how through some short  interactive problems this perception can be significantly broadened and refined, demonstrating that ethics in mathematics can be taught, and learned. It shows that mathematicians don't need to be scared of ethics and can certainly do it themselves, and moreover that deferring `the ethics' to someone else might not be very effective.
}}

\subsection{Introduction to the course}

\subsubsection{What is this course?}

The course was titled ``Ethics for the Working Mathematician'', and covered the ethical issues pertaining to the interplay between mathematics, mathematical work, and mathematicians, as researched by the Ethics in Mathematics Project\footnote{\url{https://www.ethics-in-mathematics.com}} at the University of Cambridge (see also \cite{chiodo2019teaching, muller2022situating, muller2025ethical}). We made it clear that this is primarily a course on raising ethical awareness and that the course is about teaching students the necessary skills to explore basic issues of ethics in mathematics. Some students arrived with the worrying assumption that ethics courses try to dictate to them a right or wrong way of looking at the world. We found it important to counter this early, and to tell them that while mathematical artefacts have politics \cite{muller2023mathematical}, this was not a course on (political) ideology. It was a course aimed at \textit{identifying} where and why ethical issues arise in mathematics, rather than \textit{solving} such issues.

\subsubsection{Course format}

As students might never have attended a course on, or otherwise read about, material such as this before, we found it helpful to give the rough breakdown of how the material was to be presented, so that they would know what to expect throughout the course and the logical flow that motivates it. We presented the course format as 8 lectures, broken down into the following parts:
\begin{itemize}
    \item Part 1: Showing mathematicians can learn ethics in mathematics -- L1.
    \item Part 2: Giving case studies showing there exists ethics in mathematics -- L2--4.
    \item Part 3: Explaining the origins and causes of such issues, to show that for all mathematics, there are ethical issues -- L5--7.
    \item Part 4: Outlining what mathematicians can do to address such issues -- L8.
\end{itemize}
This introductory lecture to Ethics in Mathematics -- part 1 -- sought to actively demonstrate that, through teaching and learning, mathematicians can increase their awareness and understanding of the ethical issues in mathematics; that is, they \textit{can learn ethics in mathematics}. The next three lectures -- part 2 -- aimed to show that \textit{there exists ethics in mathematics}, done through a long series of examples and case studies. The next three lectures -- part 3 -- aimed to show that \textit{for all mathematics, there are ethical issues}, done by exploring the origins and causes of the ethical issues raised earlier. The final lecture -- part 4 -- was an introduction to the next steps and what mathematicians can begin to \textit{do} to actively understand, and avoid or mitigate, the harms that might come from their work and output. We found that this structure was comparatively easy to digest for students, as some of them might struggle had we done the more theoretical \textit{for all}-arguments (covering origins and causes -- part 2) first. We also felt that students needed additional time to emotionally cope with the fact that their favourite area of science can lead to harm, and walking through some areas slowly ensured that not all students were hit with it at the same time.

\subsubsection{Why study ethics in mathematics?}\label{S:3reasons}

We first discussed with the students what might motivate each of us as individuals to study ethics in mathematics, and gave the following as potential reasons to do so, (also covered in \cite{chiodo2024field}): 

\begin{enumerate}
    \item Doing good and/or avoiding harm.
    \item Preventing negative personal consequences.
    \item Potential career opportunities.
\end{enumerate}
One reason given to study ethics in mathematics is the motivation to  \textit{do good} and \textit{prevent, or not cause, harm} with mathematical work. In this context, students were also presented with the idea that harm from doing mathematics can be accidental, and that it can even happen when the mathematics was done with the best intentions.  This was covered more explicitly in Part 2 of the course (case studies).

Another motivation given to explore ethics in mathematics was \textit{negative personal consequences}, ranging from reputational and financial harm to more stringent legal consequences such as going to jail (e.g., the case studies in Section \ref{s:misbehaving} of mathematicians such as Tom Hayes \cite{Enrich.2018} going to jail, and the further legal aspects of mathematical work covered in Lecture \ref{s:law}).

In addition, we thought some students may be motivated to study ethics in mathematics by \textit{new job opportunities}, such as helping and advising others mathematicians on ethical and responsible mathematical development. Indeed, we explained that there are many opportunities arising as a result of greater public awareness and concern for the harm that mathematically-powered technologies can bring. The recent EU AI Act was such an example of a factor which motivates companies to conduct their mathematical work in a more responsible way (see \cite{Chiodo.2024.lawpaper} for a deeper discussion on this and related issues).

When time permitted, these were put into a larger philosophical context by asking students how they relate to the question of ``what it means to live a good life''. Over the years, we have observed that students who find and understand their own motivation for studying ethics in mathematics are more likely to stay with the course, and have an easier time dealing with some of its challenges. Learning that they can do harm through mathematics can be hard for students and working mathematicians alike (cf. \cite{muller2022situating}), and thus we thought it necessary to incorporate such personal discussions early. We informed students that those in previous years who treated this as any other regular undergraduate course, by attending all lectures, taking notes, doing pre/post-lecture reading (such as \cite{chiodo2023manifesto}), attempting exercises (such as \cite{chiodo2023teaching}) etc., have reported back to us that they gained a lot from the course, with several explicitly saying ``this was one of the most valuable courses of my undergraduate degree''.

\subsection{Learning through case studies}

\subsubsection{Activity: What does ethics in mathematics mean to you?}
In this lecture, students were initially asked ``What does `Ethics in Mathematics' mean to you?''. The answers were solicited from the class and collected on the board. Typical responses included ``Nothing'', ``Mathematics can't hurt anyone'', ``Mathematics can give false objectivity'',  ``\textit{Someone else} might misuse/misappropriate our mathematics'', ``It's not the mathematics, it's what people do with it''. This was a very open solicitation of responses, and it was made clear that there is no wrong answer as we are specifically asking the students what it (currently) means to them. After this initial discussion, students were presented with two interactive mathematical problems where ethical issues might arise:

\begin{enumerate}
    \item Computing the `optimal' path for an oil pipe from a rig to a refinery.
    \item Computing the population survival rate of disease for a given mortality rate.
\end{enumerate}

\subsubsection{Interactive problem: The oil pipe problem}\label{s:oil-pipe}

\textit{This problem was initially published by Bonnie Shulman} \cite{Shulman2002}, \textit{and its extended formulation with a complete discussion can be found in} \cite[p.~56--57]{chiodo2023teaching}: 

Students were presented with a scenario whereby there is an oil rig in the ocean and a refinery on a straight shoreline. The students were tasked with describing how they would calculate the optimal path of an oil pipe connecting the two. Initially students were presented with incomplete `mathematical' information (e.g., a missing distance, missing unit costs of the pipe on land and under water, etc.) and were asked what further information they would need. Students typically mentioned the (missing) distances and the unit costs of building the pipe under water and on land as the first pieces of additional information that they would need, yet seldom raised anything else. They were then guided through other factors such as missing environmental costs and constraints (e.g., is there a marine habitat or coral reef in the ocean?). Students saw that before they can even `solve the mathematics', they need to ask further questions, make decisions, and provide proper reasoning. We have observed that many students like to rush into setting up and solving the `mathematical part', and only properly consider potential ethical aspects after explicitly being prompted to do so. We have found that, when presented with such an optimisation problem, students may tend to prefer to focus on those aspects and details which are easy to quantify, and choose to optimise by finding the minium-cost pipe, ignoring other factors. When this does happen, we find it useful to explore where this tendency is coming from and how they can overcome it. One pitfall we've seen in this context is students having very little (if any) background in philosophy, and thinking that a utilitarian perspective is (or even should be) the default stance for everyone. At this point, we briefly made  mention of other forms of ethics (consequentialism, deontology, virtue ethics, contract theory, \dots), and encouraged those who might want to explore it further to read more about it (e.g., \cite{rachels2012elements})\footnote{There is also the famous Harvard Justice course which is a useful resource for students: \url{https://pll.harvard.edu/course/justice}.}.

\subsubsection{Interactive problem: Disease modelling and mortality rates}\label{s:disease}

After having discussed the oil pipe problem, we gave another case study from disease modelling. One can model the evolution of an infectious disease via a differential equation (see \cite[pp.~39--41]{chiodo2023teaching} for full details and an in-depth discussion). This model suggests that for a disease with a rate of infection $3/4$, and a mortality rate $3/5$, the long-term surviving population will be $4/5$ of the original population. However, with a mortality rate $2/5$, the surviving population \textit{reduces} to $3/5$ of the original population (according to the model). Students are then presented with the following questions in an `open floor' discussion: If a drug is discovered that reduces the mortality rate from $3/5$ to $2/5$, what should be done? Is it a binary choice? Is it a purely utilitarian choice? We asked students to consider if their friend, parent, or child had the disease, whether they would advocate giving them the drug? Or if \textit{they} had the disease, whether they would take the drug? We also asked them what they would do \textit{in their capacity as a mathematician}, and what actual power and options they really have as a mathematician. And we asked them what they would say as the head of the National Health Service, or even as Prime Minister, and to consider what a drug company might say, and prefer to do? 

It was then explained that this is a version of the \textit{trolley problem} in philosophy \cite{bruers2014review} and law \cite{thomson1984trolley}, and that the questions given were carefully designed to guide students through some of the various options and their consequences. Students were able to see that, even with a highly mathematical problem of disease modelling, deep ethical questions can quickly arise.

\subsubsection{Activity: What does ethics in mathematics mean to you \textit{now}?}

The first problem given to students (oil pipe) was much more technical as students needed to select parameters and constraints, optimisation objectives, etc. Although it contained new lines of thinking, primarily around `other non-mathematical factors', it lay within the comfort realm of most students. The second problem (disease modelling) however, already took them beyond  technical concerns, and demonstrated the societal and political aspects that can arise from doing mathematics. In the lecture, this transition was very smooth, and they saw the connections back to the mathematical work very clearly.

At this point, students were again asked ``What does `Ethics in Mathematics' mean to you \textit{now}?'', with new answers recorded on the board alongside the previous ones. Students almost always demonstrated a substantially deeper understanding, even after only having seen just two case studies. They gave responses such as ``Interpretation and political opinions'', ``Mathematics is a tool, not a self-contained solution'', ``Quantification is hard'', ``Our responsibility as mathematicians'', ``Showing restraint'', ``What is the \textit{right} question to be answering?'', ``Ask others - impacted parties, and domain experts''.  We found it remarkable to see the dramatic improvement in their understanding in roughly 20 minutes. We observed that, even without learning solid (philosophical) theory, the mere act of considering mathematical problems with a prominent ethical component, and being asked some carefully thought-out questions, was enough for many of them to self-learn fundamental aspects of ethics in mathematics. In short: not only did the students \textit{learn} some ethics in mathematics, but they \textit{saw that they could learn} ethics in mathematics.

\subsection{Mathematicians and management}

\subsubsection{Management manipulating mathematicians}

In the final part of the lecture, we introduced responsibility gaps and the role of mathematicians in the workplace. Our concern was that even conscientious students may be tempted to assume that they can rely on management to guide them as mathematicians and give them the `correct' distillations of problems. We outlined that this might not always be the case, and indeed that some managers may try to be manipulative to achieve unethical goals (see, for example, the case study of Theranos \cite{Tourish.2023}). Foreshadowing the discussion to come in Section \ref{s:ifyou}, students were given an example of a technique  management may try and use to pressure mathematicians who object to an (unethical) request, by arguing ``If you don't do it, someone else will'', and we gave ways to combat this. Firstly, we explained the relevance of the skilled worker shortage in mathematics to this argument. We highlighted that there aren’t that many mathematicians (especially in a particular niche area), so a mathematician's strength of objection matters; if they don’t do it, there probably aren’t many who could come and take their place. Secondly, as mathematicians in training, we encouraged students to discuss the contrapositive, to see that the original statement is logically equivalent to “If no-one else does it, then you will.”, which is obviously absurd (explained in full in  \cite[p.~56--57]{chiodo2023teaching}). We presented a more accurate interpretation of this threat - ``If you don’t do it, someone else \textit{might}'' (or might not) - and discussed why it isn't formulated like this in the first place.

Ultimately, this part led students to see that their objection \textit{does} matter, and that if they say no to a task, it may well prevent the task from happening. Here, some students raise the (valid) argument of ``What if I could do the task in a \textit{less harmful} way than a colleague; doesn't that mean I should do it?'' We presented two lines of reasoning to help the students  work through this:

\begin{enumerate}
    \item This is indeed a possibility, and they should think carefully on a case-by-case basis about how likely it is that a colleague will actually take up the task, and how much  harm they could realistically prevent by carrying it out themselves.
    \item By openly, vocally objecting to a task, they may motivate their colleague to also object to doing it. There is a well-known saying which is relevant here: ``Courage is contagious''.
\end{enumerate}

\subsubsection{Management lacking skills, understanding, and time}

Foreshadowing the discussion to come in Section \ref{s:defer}, we concluded this section by explaining that, even when a manager \textit{is} benevolent, they may still lack enough understanding of the technical details and/or societal context to take on all the responsibility and be able to guide the mathematicians in their team through the ethics of the problem. We have observed that some students are more than happy to defer responsibility when it means that they can focus their attention on doing mathematics, and so we reasoned that an all-knowing manager us unlikely to exist . 

However, as a thought experiment, we proposed the  `magic manager'; one who possesses sufficient (technical, and societal) knowledge to guide mathematicians through all ethical problems, thus releasing them from all responsibility. However, we then made students aware of the complexities of the workplace and the myriad of other tasks that managers have to do on a constant basis, and that even with a `magic manager', such a person is likely to be under so other many other constraints  that they simply lack the time to check everything the mathematician does. 

We found it useful to end the lecture with this discussion, as it gently introduced some additional reasons for why it is necessary for mathematicians to take at least some responsibility for considering ethics, as it may be the case that nobody else can, or will, do it for them. We explained that working mathematicians often need to understand more than mathematics, and indeed, ``to be an effective member of the team, the mathematician must also understand the basic principles of the various disciplines which he is expected to discuss. He [sic] should be, in other words, the sort of man who a century ago was known as a natural philosopher -- a man who had a keen analytical mind, adequate mathematical training, and a broad and sympathetic interest in a wide range of natural phenomena'' \cite[p.~938]{Fry.1964}, and the many applications of modern mathematics have made it necessary that this includes ethical considerations and taking appropriate responsibility.

\begingroup
\let\section\subsection
\renewcommand{\bibname}{Bibliography (Lecture 1)}
\renewcommand{\refname}{Bibliography (Lecture 1)}
\addcontentsline{toc}{subsection}{Bibliography (Lecture 1)}
\renewcommand{\addcontentsline}[3]{}
\bibliographystyle{IEEEtran}
\bibliography{bibLec1}
\endgroup

\newpage
\section{Lecture 2: Financial mathematics and modelling}\label{s:financialmaths}

\fbox{
  \parbox{0.98\textwidth}{
Mathematical modelling is ubiquitous in understanding the way the world works, from finance to physics to disease transmission. Understanding how to develop and use a model, as well as its limitations and the way it interacts with the world, is indispensable in ensuring it is a tool for good, rather than a mechanism for harm.

This lecture looks at how mathematicians can avoid causing harm when doing mathematical modelling; harm to others, and harm to themselves. It shows how they can prevent such work tanking financial markets through misapplication, or propagating pathogens or injustices in society though poorly chosen assumptions and use cases, or even landing them in jail through unauthorised actions.
}}

\subsection{Mathematicians misbehaving in finance}\label{s:misbehaving}

We started the lecture as bluntly as possible, by presenting a list of mathematicians in finance who were caught up in significant legal proceeding, several of whom spent time in jail. While these were elaborated on more in Lecture \ref{s:law}, we found it useful to present them here as they were all from finance, and moreover it gave students a very sudden, sharp `shock' to help them realise that there were significant \textit{personal} consequences for irresponsible behaviour as a mathematician (as introduced in Section \ref{S:3reasons}), and that mathematicians were not invincible. We made it clear that it was the intentional choices and actions of these mathematicians that led to such unpleasant consequences, and were by no means `accidental'. We also highlighted that, by attending our course, students were helping to reduce the chance that they themselves would make similar decisions in the future, having seen some of the consequences.  Those we spoke about were:

\subsubsection{Example: Traditional finance -- LIBOR rigging}\label{s:libor}

Each day, 16 of the largest banks in London submit how much they would expect to pay to borrow money (i.e., a rate of interest) for a 3 month loan. The average of these is the \emph{London Interbank Offered Rate} (or \emph{LIBOR}) for that day \cite{.07072012}. These interest rates are used as variables in pricing trillions of dollars of financial products (approximately USD\$150trillion in 2008 \cite{Twaronite.17092008}). The banks were thus able to influence a number (LIBOR) that in turn determined the value of their portfolios. Several people working in finance were sufficiently tempted to rig it to improve their own portfolio value, and were subsequently charged with various forms of fraud, with some eventually going to jail \cite{Channel4News.2017}; several of them having come from a mathematical (or similar) training background.

Tom Hayes was one such individual. He has a mathematics and engineering degree from the University of Nottingham. He worked in finance from 2001-2010, and in that time was alleged to have rigged LIBOR by offering bribes to others in the industry to slightly change their bank's daily LIBOR submissions.  In 2015 Tom was found guilty of rigging LIBOR, and was sentenced to 14 years in prison, serving 5 years in custody before being released on parole. There have been articles \cite{Enrich.30032025} and books \cite{Enrich.2017} written about Tom. Tom maintained his innocence throughout the process, arguing that he was simply doing what everyone else was doing, and saying ``Well look, I mean, it’s a dishonest scheme, isn’t it? And I was part of the dishonest scheme, so obviously I was being dishonest.'' \cite{Liam_2017}. In 2025, Tom had his conviction overturned as he had been ``deprived of a fair trial'' \cite{Kalyeena_2025}. But even though there was ``ample evidence'' for a jury to find him guilty, ``the jury was not properly directed,'' and so ``The convictions are therefore unsafe and cannot stand'' \cite{Kalyeena_2025}.

Another such individual was Alex Pabon. With degrees in engineering and then financial mathematics, Alex then went to work in finance from 2002-2006. In that time, he was also alleged to have rigged LIBOR. In 2016 he was found guilty, and sentenced to 2 years 9 months in jail \cite{Brinded.04122017}.

\subsubsection{Example: Quantitative finance}

Ke Xu did his undergraduate degree in mathematics at the University of Cambridge, graduating in the top 3 of his cohort in 2007 \cite{Chellel.01122019}. After graduating, he went to work for Goldman Sachs for a few years, and then in 2012 move to the quantitative finance firm Trenchant Ltd (commonly known as ``G-Research'') producing algorithmic trading software. In 2014, he ``only'' received  GBP\pounds400,000 in bonus; a sizeable sum, but far less than the \pounds 1.1million he had been expecting \cite{Fortado.07032015}. Dissatisfied with this, he allegedly spent the next  several months  reverse-engineering the trading code at G-Research, then took it and abruptly flew to Hong Kong \cite{Chellel.19112018}. He was extradited back to the UK shortly after, convicted of various crimes related to his alleged actions, and served several prison sentences until he was eventually deported back to his home country of China in 2019 \cite{Chellel.25012019}.

\subsubsection{Example: 2010 flash crash}

On May 6, 2010, there was a ``flash crash'' across the US stock markets \cite{Kirilenko2017}. They suffered a 9\% drop, and almost identical rebound, in the space of 36 minutes \cite[Section I]{CFTC2010}. It has been argued that Navinder Sarao, a trader, triggered the flash crash. He apparently did so using a technique called ``spoofing'' \cite{Tarm2020}, and profited \$900,000 from his trades that day. Navinder was described as a ``maths wiz'' \cite{TheSeattleTimes.2020} and as having had ``some extraordinary abilities with respect to pattern recognition and certain sorts of mathematical abilities'' \cite{Vaughan.02102017}.  Arrested and charged in 2016; he was initially jailed for 4 months before being extradited to the US for sentencing \cite{Staff.11102016}, and eventually received a sentence of 12 months under house arrest \cite{Tarm2020}.

\subsubsection{Example: Cryptocurrency exchanges}\label{s:ftx}

In late 2022 the cryptocurrency exchange FTX Trading Ltd (an abbreviation of ``Futures Exchange''), and it's sister cryptocurrency trading firm Alemeda Research, collapsed into bankruptcy in spectacular fashion \cite{guardian_ftx_2022, ft_ftx_2022}. This was as a result of alleged illegal behaviour by several of its executives in transferring depositor funds from FTX to Alameda for investment purposes. The executives of Alameda took a ``high risk'' approach to investing \cite{Sier.21112022}, and subsequently many of those funds were lost by Alameda, causing the double-bankruptcy and loss of billions of dollars of depositor funds.

The key executives from FTX and Alameda all received undergraduate degrees in mathematics or related disciplines in the mid-2010's: 
Sam Bankman-Fried (MIT) \cite{Zeidan.31032025}, Caroline Ellison (Stanford) \cite{TodayUsa.24092024}, Gary Wang (MIT) \cite{Larson.20112024}, and Nishad Singh (U.C. Berkley, \textit{computer science}) \cite{Explainers.14112022, Osorio.26122022}. In the mid-2020's, they all received criminal convictions and custodial sentences: Bankman-Fried (25 year \cite{Zeidan.31032025}), Ellison (2 years \cite{TodayUsa.24092024}), Wang (time served and 3 years supervised release \cite{Ciccomascolo.21112024}), Singh (time served \cite{Ligon.30102024}). The prosecutor for Bankman-Fried was particularly disparaging in his sentencing memorandum, stating ``Time and time again the defendant has expressed that his preferred path is the one that maximizes his version of societal value, even if imposes substantial short term harm or carries substantial risks to others'' \cite[p.~88]{USvBankmanFried2023}.

\subsection{Modelling in finance}

We then moved on to explain that mathematical work in finance can harm not just the people doing it, but also wider society, and that this can manifest as systemic harm from the interactions between multiple actors (mathematicians, markets, the public, etc.). We highlighted the idea that, even though no one person might be to blame for these harmful outcomes, each person involved shoulders some of the responsibility, and could have collectively changed their behaviour to avert such harm.

\subsubsection{Example: The Black-Scholes equation and the 1987 stock market crash.}\label{s:blackscholes}

An \textit{option} is a contract between two parties A,B such that A has the `option' of buying some quantity of goods/stocks from B at some fixed future time for some fixed price. The purchase is not compulsory; it is not a contract of sale.
 
In 1973  Fischer Black (a mathematician-turned-economist) and Myron Scholes (an economist)\footnote{Who won the Nobel Prize in economics in 1997, and whose hedge fund, Long-Term Capital Management, went bust the following year \cite{BBC27042012}} developed the \textit{Black-Scholes Equation} (mathematically identical to the \textit{heat equation} in physics) \cite[p.~644]{Black.1973} as a way to price options, \emph{sometimes}.  This equation was derived from the \textit{Black-Scholes Model} \cite[p.~640]{Black.1973}, in which 7 \textit{assumptions} are made; 3 about the assets themselves, and 4 about the behaviour of the market. Like every mathematical model, this is somewhat of a simplification. And some assumptions were later found to be demonstrably (and catastrophically) false, such as the assumption of \textit{constant volatility}.
 
This was the first time a pricing method for options had been widely accepted (by buyers \textit{and} sellers). As options could be priced, it meant they could be traded, creating a large market for options throughout the 1980s. Unfortunately, many traders thought the Black-Scholes formula was a way to price options \emph{all the time}. But  the formula could fail to hold, particularly during times of high market volatility. And with countless financial firms using this formula to trade options (several even doing so via automated trading), when volatility did hit the markets in October 1987, it led to a stock market crash wiping nearly $\$2$trillion off global markets, including a $20\%$ drop in the US stock market \cite[p.~123]{McKeon.2009}; it's largest one-day drop in history.
  
Warren Buffett said of the Black-Scholes formula in 2008 \cite{Buffett.2008, Gupta.2016}: ``If the formula is applied to extended time periods, it can produce absurd results. Black and Scholes almost certainly understood this point. But their devoted followers may be ignoring whatever caveats the two men attached when they first unveiled the formula.'' This leads to a key principle in using mathematical models: 

\begin{center}
\textit{ If we don't understand the derivation, \\then we can't know the limitations.}    
 \end{center}
Nowadays the equation is still used, albeit with several extra (mathematically `ugly') correction terms included. It has become somewhat of a `Frankenstein'-Black-Scholes Equation \cite[p.~855]{Ki.2005}; after all, why should a market necessarily be modelled by a `nice' equation?

\subsubsection{Example: Collateralised debt obligations and the 2008 global financial crisis.}\label{s:copula}

A \textit{Collateralised Debt Obligation} (CDO) is a way to take a collection of interest-bearing assets (such as mortgages), combine them into a pool, and then cut up the pool into pieces (called \textit{tranches}) so that each tranche is an (abstract) interest-bearing asset. The tranches are ordered, and defined by two numbers: a start, and end, percentage. As income comes into the pool each month, the first tranche (call it the $[0\%,N\%]$ tranche) receives the first $N\%$ of payments due for the month. The next tranche, $[N\%,N+K\%]$ receives the next $K\%$ of payments above $N\%$, and so on. The tranches don't `care' \textit{which} of the mortgages are making payments, so long as \textit{some} do to reach their payment interval. When defaults occur, the bottom tranche takes the first losses; the top tranche takes the last. Thus the bottom tranche carries the greatest financial risk, whereas the top tranche carry the least (roughly speaking).

To understand the risk of default of a tranche requires an understanding of how all the mortgages in the pool interact. David Li's paper in 2000 \cite{li1999default, li2000} presented a model explaining the ``joint'' probability of differing mortgage holders defaulting: if one defaults, what is the expected time until the other defaults? This model, called the \textit{One-factor Gaussian Copula Model}, offered high analytical tractability, meaning it can be solved analytically in terms of a closed form expression, and hence can be computed quickly. But it also made many assumptions, which roughly held during relatively stable market conditions, but could break down in times of market instability and decline.

To create the tranches, mathematicians would select some crossover points where the tranches begin and end, and use the Gaussian Copula Model to compute the default risk of each tranche. This would then be sent off to ratings agencies, and each tranche given a rating. Thus, each tranche would have an expected monthly income, and a rating, so could easily be priced (and traded) on financial markets. This would be attempted several  times, with different crossover points, to try and get the highest overall value. Because the top tranches were generally considered low-risk (high rating), and because of the price premium the market placed on highly rated interest bearing assets, the combined values of the tranches always exceeded the values of the mortgages placed into the pool. So this mathematical process generated instant, large, `guaranteed' profit.

The Global Financial Crisis came to a head in 2008. At the heart of the crisis was the use of CDOs. The mathematicians creating and pricing CDOs simply did no have the depth of understanding of the assets, or the model, to use it safely. Unfortunately, some of the mathematical assumptions did not hold, such as the assumption of zero tail dependence; the model assumed that low-probability events (i.e., two different mortgages defaulting) would be uncorrelated. But if both houses were on the same street (or, more generally, both houses were exposed to the same economic factors), then this simply wasn't true.

The boom in CDO sales led to a boom in the (US) mortgage market. Mortgages were a way to make free money, and brokers knew they could sell them to hedge funds for an immediate profit. So they were motivated to issue more mortgages, even to those who could not necessarily afford them (so-called \textit{sub-prime mortgages}). Thus the mathematicians incentivised different market behaviour.
 
In the end, the risk was not properly accounted for, and when house prices declined it led to the write-down of trillions of dollars of losses; the global financial intuitions recorded a collective write-down of USD\$420billion of CDO value (out of a total \$621billion of CDO issuance)  by 2007 \cite[p.~2]{cordell2011collateral}. The system that had been set up and fuelled by mathematicians working in finance came crashing down, an their ``formula had killed Wall Street'' (see also \cite{mackenzie2014formula}).

\subsection{Modelling during COVID-19}

Our focus then shifted slightly to another high-impact area of mathematical modelling: modelling during the COVID-19 pandemic. Given the complicated, novel nature of pandemics and their surrounding circumstances, we wanted to highlight to students that there are times where modelling is viewed as the only way though an otherwise insurmountable situation; an essential approach. However, we chose the following examples to highlight that, if not done carefully, modelling might bring about more harm than it set out to prevent.

\subsubsection{Example: Modelling examination results}

As described in \cite[p.~4]{chiodo2023manifesto}: ``During 2020, in the midst of the COVID-19 pandemic, the UK government cancelled all school examinations. In lieu of exams, the UK Office of Qualifications and Examinations Regulation (Ofqual) was tasked with `calculating grades' for final-year students \cite{BBC.17082020}. These grades were primarily intended to be used to determine which students qualified for university entrance. The grading system asked school teachers to rank the students at their school within each subject (say, final-year chemistry), and then used that ranking to assign grades in a normal distribution derived from the distribution of final-year chemistry students at that school for the past 3 years. That is, all that determined a students grade for a subject was their assigned ranking in \textit{their} year, and how well students \textit{in previous years} had performed. Teachers were also asked to predict grades for each student in each subject, but for most students this data was ignored in the calculation. However, for statistical reasons, those in small subject cohorts of size $<5$ received their teacher predicted grade instead, which were on average much higher than the Ofqual-calculated grades. A total of $39\%$ of the calculated grades were one grade lower than predicted by teachers, and $3\%$ were two grades lower (i.e., A$\rightarrow$C). When the grades were released to students in August 2020, this caused an outrage. Students felt very hard done by, as many received calculated grades much lower than their teacher predicted grades. And in addition, private schools had many more subject cohorts of size $<5$ than public schools, giving them an inherent advantage with more teacher-predicted grades. The public outcry was immense, and after 4 days of immense turmoil the UK government rescinded all Ofqual-calculated grades and replaced them with teacher-predicted grades.
    
Here, mathematicians were asked to go well beyond the limit of what is `predictable' by mathematics, and with minimal evidence calculated unfalsifiable grades for students that the public revolted against.'' \cite[p.~4]{chiodo2023manifesto}

\subsubsection{Example: Modelling the reopening of universities}

In mid-2020, in the midst of the COVID-19 pandemic, universities around the world were trying to understand how they might operate for the upcoming academic year. The pandemic and subsequent government lockdowns and restrictions had led many to close in early 2020, but for many there was more flexibility in how they might operate in the 2020--21 academic year, subject to various rules regarding social distancing, infection testing, contact tracing, etc.

As well as improving the student learning experience, universities were under an additional pressure, in the form of financial losses. With halls of residence laying vacant with students studying from home, universities were facing a vast loss of income from student rental. Thus, several universities around the world, such as the University of Illinois, Urbana-Champaign, and Pennsylvania State University, in the USA \cite{QuantaMagazine.2021}, and teams within the UK \cite{NGTM.07042025}, started to carry out their own COVID modelling, in anticipation of reopening at the start of the 2020--21 academic year.

However, as careful as these modellers were in carrying out their work, many of the models were inaccurate, and grossly under-predicted the spread of COVID on university campuses. For example, at the University of Illinois, Urbana-Champaign, ``on September 3, just one week into its fall semester, the university faced a bleak reality. Nearly 800 of its students had tested positive for the coronavirus -- more than the model had projected by Thanksgiving.'' \cite{QuantaMagazine.2021}. And at Manchester Metropolitan University, thousands of students were forced to self-isolate early into the first term \cite{bbc2020manchester}.

Some of these oversights might, in hindsight, have been avoided with additional consideration. While much attention was paid to the physics of disease transmission, insufficient attention was paid to the more social aspects of human behaviour such as (non)adherence to rules. Indeed, ``The scientists had seemingly included so much room for error, so many contingencies for how students might behave. ``What we didn’t anticipate was that they would break the law,'' Goldenfeld said -- that some students, even after testing positive and being told to quarantine, would attend parties anyway. This turned out to be critical: Given how COVID-19 spreads, even if only a few students went against the rules, the infection rate could explode.'' \cite{QuantaMagazine.2021}

These mathematical models simply hadn't worked, and moreover had given universities a false sense of security about how (un)safe reopening really was. Students were invited back into residence, only to be then forced home, or forced to self-isolate, far away from the (practical, logistical, social, and mental) support networks they had at home.

\subsubsection{Example: Modelling COVID-19 transmission in aged care homes}\label{s:covid}

Early in 2020, in the initial months of the COVID-19 pandemic,  much mathematical modelling was done in the UK to try and predict the spread of COVID, and the effects of such spread on health and on the healthcare system. Initial steps were taken to protect those most vulnerable, i.e., the elderly. In particular, attention was paid to places with high numbers of elderly people: aged care facilities.

Through various models, it was believed that aged care facilities could avoid the worst affects of the outbreak by having a ``protective ring'' placed around them \cite[54:41]{YouTube.07042025} whereby only staff could enter and exit; all other people, such as family visitors, were banned from entry. Modellers suggested that this ``cocooning would reduce deaths'' \cite[21:06]{YouTube.07042025}.  

However, the models were grossly inaccurate, and even with the measures described. In particular, the movement of agency staff \textit{between} care homes propagated the virus throughout the entire aged care sector, and this was not factored into the modelling \cite[21:14]{YouTube.07042025}. The UK government thought they had solved the problem for care homes, but they hadn't. Consequently, there was a large outbreak of COVID across the entire UK aged care sector in Spring 2020, and ``Between 2nd March and 12th June 2020, there were 66,112 deaths of care home residents in England and Wales, of which 19,394  are officially attributed to COVID-19'' \cite[p.~986]{daly2020covid}.

When interviewed about this failure, one of the mathematicians involved in the modelling, Dr Ian Hall of the University of Manchester, said ``The failures of those models, I guess, was that we didn't know how connected the social care settings were with the community. As modellers we didn't know. I'm sure there are lots of academic and policy makers out there that could have told us this, if we'd asked them.'' \cite[21:23]{YouTube.07042025}. Later in the same interview, he said ``We realised that care homes were not isolated and shielded in the way that we thought they were, as modellers. But yeah, there was some sort of external connection that was forcing the disease''. And in response to the follow-up question ``So what made you think they were shielded?'', he paused for a 12 seconds before answering ``That's a good question. We never checked.''\cite[53:45]{YouTube.07042025}.

\subsection{Understanding the limitations of mathematical modelling}\label{s:limitsmodelling}

We concluded this lecture with a summary of some of the pitfalls of modelling. We explained to the students that the purpose of creating a mathematical model is to gain \textit{predictive power}; to be able to use mathematics to predict something about the world without seeing it happen first. Drawing on the examples given above, we suggested that, when creating a mathematical model, mathematicians might be inserting into the mathematics some set of beliefs, values and opinions on how a system functions. 

We explained that in areas such as  physics there is reproducible experimentation to test the hypotheses and robustness of a model, and that models of the laws of physics can be (rather quickly) falsified. However, again drawing from the examples above, we pointed out to students that modelling finance, or a pandemic, is not like modelling physics; it is very difficult to  verify that the model is correct (what is the `right' price, or how `will' the virus spread?). We pointed out that, if a model is wrong, yet being widely-used, then it may only be a system-wide failure (such as a stock market crash, or mass transmission of a virus) that reveals it is wrong, and that those using the model might not understand its derivation and thus fail to see its limitations. 

We extended this concept, by pointing out an additional complication about financial or epidemiological modelling: the model affects the thing being modelled. We highlighted that while the laws of physics do not change regardless of what model is being implemented, the behaviour of a financial market, or a society in a pandemic, can, and will, be affected by the predictive models being used. And so we posed to the students  ``How do you model the impact of your model?'' We then gave the students our interpretation of responsibility in modelling, which is that if a model fails, it is not because the world is wrong, but rather because the model, and by extension modellers, were wrong. We suggested that the world is not obliged to follow the elegant, or `nice', second-order differential equations set up by mathematicians, so they should not be surprised when it doesn't.

To finish, we presented the students with a (crude) classification of modelling into two types: 

\begin{enumerate}
    \item Modelling something that you're going to do anyway. For example, a hospital will distribute money to wards anyway, and a power company will take a hedge on their customer power usage anyway.
    \item Modelling something such that the modelling makes you think it possible to achieve (creating a market or possibility), incentivising people to do something new. You don't go to space unless you've done some modelling. You model so to `decide' it is safe to bring university students back on campus. You model so to algorithmically predict A-level grades. You model so to price and trade options and CDO tranches.
\end{enumerate}
We ended by explaining that any claims of the form `It is always better to model a situation than not', mostly applies to the first case above, but not the second. We suggested that in (1) the model (hopefully) makes an existing (potentially bad) situation better, but that in (2) the model might create a (potentially bad) situation that was otherwise not going to happen, and thus ``inadvertently solve one crisis with a bigger one''  \cite[p.~1]{chiodo2026handing}.

\begingroup
\let\section\subsection
\renewcommand{\bibname}{Bibliography (Lecture 2)}
\renewcommand{\refname}{Bibliography (Lecture 2)}
\addcontentsline{toc}{subsection}{Bibliography (Lecture 2)}
\renewcommand{\addcontentsline}[3]{}
\bibliographystyle{IEEEtran}
\bibliography{bibLec2}
\endgroup

\newpage
\section{Lecture 3: Cryptography, surveillance and privacy}\label{s:cryptoprivacy}

\fbox{
  \parbox{0.98\textwidth}{
Privacy is a critical part of how society functions, protecting people at the times and places where they are most vulnerable. Mathematicians have always played a central role in the creation and implementation of privacy-preserving, and privacy-infringing, technologies, in many cases determining how much privacy people can have.

This lecture looks at how mathematicians can avoid  infringing the privacy of others, as this can happen even completely unintentionally. It shows how they might do so deliberately, by designing tools to break strong encryption, or indirectly, by creating systems and platforms which collect massive amounts of personal data, or accidentally, by being careless or sloppy in the way they store and share data.
}}

\subsection{Introducing dual use of mathematics}\label{s:dual}

In this lecture we took the opportunity to introduce the concept of \textit{dual use} of mathematical technologies. It seemed particularly relevant here, as dual use of cryptographic tools and exploits can be fairly straightforward to see, with the tools themselves needing very little, if any, change to go from `good' to `bad'.

We started off by getting students to think about who the `good guys' are, and whether this is always obvious, and always clear to a mathematician working in areas of cryptography and surveillance. We used a  straightforward worked example from cryptography \cite[p.~30--32]{chiodo2023teaching} to illustrate it. This was done by posing a standard problem of the form `Alice and Bob are communicating using encrypted means, with an eavesdropper Eve listening', where one could do some mathematics to show that Eve can indeed decrypt the messages. The question put to the audience was `Who do you tell, and what do you tell them?' However, by changing only the setting of the actors in the question (i.e., who Alice, Bob, and Eve were), but not any of the mathematics, the students saw very different versions of the same problem. For example, we first presented no information about Alice, Bob, or Eve, and the students suggested to ``Tell Alice and Bob that Eve can decrypt their messages''. We then \textit{added} that Eve was part of a law enforcement agency which suspected Alice and Bob of committing a crime; the students then suggested to ``Tell Eve how to decrypt the messages''. Finally, we added that the law enforcement agency was in a country with a questionable human rights record, and Alice and Bob were suspected of being political dissidents; the students then reverted their suggestion back to ``Tell Alice and Bob that Eve can decrypt their messages''. Full details of this exercise, with worked solutions and further discussion, can be found in \cite[p.~30--32]{chiodo2023teaching}.

This discourse served to get students thinking about who the `good guys' are in a given situation, how that might change with a differing perspective or evolving situation, and ultimately how that might affect how they carry out their mathematical work. We used this to lead into a discussion about privacy.

\subsection{What is privacy?}

Having touched on how mathematicians might be involved in work that affects people's privacy, we presented some ideas about what privacy actually is, and how a mathematician might view it.

We started by presenting privacy as a complicated notion: an idea many people  rely on, with a deep philosophical history,  yet  lacking a rigid mathematised definition. We brought up the maxim \emph{If you've done nothing wrong, you've got nothing to hide} (a play on, but not the contrapositive to, \emph{If you've got nothing to hide, then you've done nothing wrong}), and then asked students to think about whether they would be comfortable publishing their personal details, history of romantic partners, school and university results, banking details, or email logins in public, and how this reconciled with their understanding of the aforementioned maxim.

We described privacy as the tool that people use to protect themselves at their most vulnerable times and in their most vulnerable positions. We pointed out that client-lawyer, and patient-doctor, communication is considered `private' (and legally protected in many jurisdictions), as people are extremely vulnerable in those situations. We pointed out some of the consequences if someone's parents, workplace boss, or any political opponents, happened to have access to all their medical and legal records. We explained that there are times in life when someone might do something which is not necessarily `wrong', but which opens them up to opportunistic attack or judgement from others, highlighting that  privacy is an important tool for people to protect themselves from undesirable consequences of their (legitimate) actions.
 
We then tied this back to the work of mathematicians, and how they can potentially violate the privacy of practically every person on earth. We presented a categorisation of this into three main drivers of how the violation might be carried out: deliberately, indirectly, or just from sheer ignorance. This was to clarify that mathematical work need not be malicious in order to be harmful.

Expanding this categorisation was what the remainder of the lecture covered.

\subsection{Mathematicians deliberately violating privacy}

We first gave several examples of mathematical work where privacy violations were an obvious consequence of the work being carried out, along with some discourse on why such work was conducted, and the potentially harmful effects. We emphasised that there was no `obvious' ethical conclusion to be drawn from any of these, but rather than one might wish to appreciate and evaluate these benefits and drawbacks before embarking on such work. 
  
\subsubsection{Example: Mass surveillance -- Snowden revelations}

In 2013, the US National Security Agency (NSA) contractor-turned-whistleblower Edward Snowden released an immense trove of documents to journalists \cite{greenwald2013snowden}. These were known as the \textit{Snowden revelations}, and detailed how the NSA, and its partner surveillance institutions in the UK, Australia, New Zealand, and Canada (collectively known as the \textit{Five Eyes} nations), had set up a comprehensive digital surveillance system that allegedly intercepted and collected a vast proportion of global digital transmissions (i.e., calls, sms, emails, fax, internet traffic, etc.) \cite{macaskill_borger_2013}. Starting with the revelation of a court order obliging Verizon to continuously hand over ``telephone metadata'' of all calls it processes \cite{greenwald2013}, more and more surveillance systems and processes were revealed. This included:

\begin{enumerate}
     \item The \textit{PRISIM surveillance programme}, in which the NSA was alleged to have tapped into the servers of key US internet giants Microsoft, Yahoo, Google, Facebook, PalTalk, AOL, Skype, YouTube, and Apple \cite{greenwald2013nsa}.
     \item The alleged targeted spying by the NSA on European Union Offices \cite{poitras2013nsa} including the \textit{Dropmire surveillance programme} which intercepted messages sent via encrypted fax machines to the EU mission in Washington DC \cite{macaskill2013}.
     \item The \textit{Tempora surveillance programme}, in which  GCHQ was alleged to have tapped several high-traffic internet fibre-optic cables such as transatlantic cables \cite{macaskill_borger_2013}.
     \item The \textit{TURMOIL surveillance programme}, in which the NSA collected and decrypted significant amounts of (encrypted) VPN traffic in near-realtime \cite{spiegel2014, adrian2015}.
\end{enumerate}
These are just some of the systems and tools the NSA and its other Five Eyes counterparts are alleged to have developed and deployed. The technical and physical infrastructure was immense, with claims that the associated capabilities were violating the privacy of billions of people worldwide.

Whatever the activities were that these surveillance organisations were carrying out -- full details are difficult to ascertain -- mathematicians were at the forefront of these capabilities, as articulated by the former director of research at NSA Michael Wertheimer who used descriptive terms such as ``NSA mathematicians lead efforts'', ``NSA mathematicians are at the forefront'', and ``Mathematicians are leading the way'' \cite{landau2015}. This is also reflected in the scale of the teams of mathematicians partaking in such activities; the NSA is the largest employer of mathematicians in the USA \cite{nature_nsa_2013}.

\subsubsection{Example: Creating a cryptographic exploit -- Diffie-Hellman key exchange}\label{s:DHkey}

In 2015, a team of computer security researchers released a paper on some weaknesses in implementations of Diffie-Hellman key exchange, which they proposed may have been exploited by the NSA \cite{adrian2015}. Diffie-Hellman key exchange involves two parties A and B sharing information over an insecure channel to establish a shared (secret) key. This is done by agreeing on a prime number $p$, and a generator $g$ for the multiplicative group of integers modulo $p$. A and B then choose secret integers $a,b$ respectively, and then send each other $g^{a} \mod p$ and $g^{b} \mod p$ respectively. Each can then compute  $g^{ab} \mod p$, which becomes their shared secret key \cite[p.~6]{adrian2015}. To break this protocol requires taking discrete logarithms, which is known to be a computationally difficult problem. Difficult, but not impossible.

The best known general methods for attacking Diffie-Hellman is to compute discrete logarithms \cite[p.~6]{adrian2015}, and the most efficient method for that is the \textit{number field sieve}; a technique from algebraic number theory. This requires two main steps. The first is the \textit{precomputation}, which entails computing a large lookup table for factorisation of many integers modulo $p$ with only `small' factors up to some given size. In the second step, the lookup table is then used to derive the shared key, making use of $g$, $g^{a} \mod p$, and $g^{b} \mod p$. The vast majority of the computational work is in step 1, depending only on $p$, and with step 2 near-trivial in comparison. Thus, if $p$ were to be re-used, and an attacker had already precomputed the lookup table for $p$, then they could easily compute the shared key each time $p$ was reused \cite[p.~6]{adrian2015}. This is the exploitation that the authors claim the NSA might have carried out.

The authors observed that many websites and protocols implementing  Diffie-Hellman  continuously re-use the same prime $p$. This is an implementation choice, coming from the fact that a prime $p$ is ``safe'' if $p-1=2q$ for some other prime $q$; this is to ensure that the multiplicative group modulo $p$ (which is cyclic of order $p-1$) has the minimal number of subgroups (of orders $2,q,2q$) \cite[p.~7]{adrian2015}. However, generating a safe prime can take a few moments \cite[p.~7]{adrian2015}, and so for user convenience developers instead chose to generate a single safe prime and hard-code it into their implementation, on the assumption that precomputing the lookup table for the prime $p$ would be prohibitively expensive and time consuming, if $p$ was sufficiently large. However, this was not necessarily the case, especially as the precomputation is parallelisable.

The authors of \cite{adrian2015} first observed that there were several `short' 512-bit primes being re-used, and that they themselves could compute a lookup table for such primes in a little over a week on their standard machines \cite[p.~9]{adrian2015} (taking 7,600 +  21,400 + 60,000 core-hours in total; 10 core-years), and then complete step 2 of the attack in ``about 70 seconds'' \cite[p.~9]{adrian2015}. They also estimated that lookup tables for 768-bit primes (which were also being commonly re-used) could be computed in 36,500 core-years; a process that could be reasonably carried out on academic resources at the time \cite[p.~11]{adrian2015}. In addition, they estimated that  lookup tables for 1024-bit primes (again, these were also being commonly re-used, and were widely considered as being totally secure) could be computed in 45million core-years \cite[p.~12]{adrian2015}; a process that they estimated could be carried out in under 1 year with hardware costing of ``the order of hundreds of millions of dollars'' \cite[p.~11]{adrian2015}, and thus within reach of a state actor.

At this point, the authors made use of the Snowden revelations, and in particular the TURMOIL programme mentioned above, to assert that a likely way the NSA was able to decrypt the amounts of VPN traffic mentioned in these revelations would be if they had precomputed \textit{one} particular prime number. This would have allowed them to decrypt traffic from 25.7\% of VPN servers using SSH and about 60\% of those using IKE (\textit{all} of which were re-using the \textit{same} prime) \cite[p.~14]{adrian2015}. Similar statistics were presented for HTTPS (web) connections and SMTP (mail) connections; there, precomputing the next 9 most popular primes (alongside the one above) allowed decryption of approximately 30\% of encrypted traffic \cite[p.~14]{adrian2015}. In addition, the authors pointed to leaked NSA budget reports from 2013 \cite{black_budget_2013} which showed a jump increase of \$100M in a program named ``cryptanalytic IT services'', and another \$100M increase in one named ``cryptanalysis and exploitation services program C'' \cite[p.~12]{adrian2015}, consistent with their estimates on the hardware requirements for precomputing primes in under a year. 

Not only did the authors demonstrate a distinct possibility that the NSA \textit{could} (and may well \textit{is}) implementing such cryptographic attacks, the work of the authors can be used to uncover another concern: the potential for a leak of such lookup tables. Their own work showed that the lookup table for a 512-bit prime ``occupies about 2.5 GB in ASCII format'' \cite[p.~9]{adrian2015}. The estimates in their Table 2 \cite[p.~12]{adrian2015} show a scaling factor of approximately 2500 in disk space for a lookup table for a 1024 bit prime ($\frac{5.2B}{2.1M}$), resulting in an approximate file size of a little over 6TB (that is, 2500 x 2.5GB). Such quantities of data, even at the time, could easily fit on a couple of commercially-available portable hard drives and stored in a backpack. If the conclusions of the authors were true, then not only was the NSA explicitly constructing an exploit to decrypt vast amounts of internet traffic, but they were implicitly creating the possibility of a data leak that could hand such capabilities over to almost any actor, even ones with very limited resources.

\subsubsection{Example: Loss of a surveillance tool -- WannaCry ransomware}

The idea that a large (i.e., state) surveillance organisation could lose or leak some of it surveillance tools is not purely hypothetical. This has transpired in real life, even from the NSA itself.

In 2017, a group called the Shadow Brokers leaked 1GB of NSA software exploits into the public domain, including a tool known as EternalBlue which allows remote code execution on Windows operating systems \cite{gallagher2017}. This was a very powerful set of \textit{zero-day exploits}; ones that the developer (Microsoft) was unaware of. And the NSA kept this knowledge hidden from Microsoft for over 5 years, even though the NSA  knew ``the widespread havoc it could wreak if it ever got loose'' \cite{nakashima2017}. Shortly after this entered the public domain, hackers were able to make use of these exploits and create a \textit{ransomware} attack called WannaCry \cite{nakashima2017}. This encrypted the user's files, and requested they pay a ransom of approximately USD\$300 in Bitcoin to have them decrypted. As well as affecting individual users, this also affected multiple large organisations around the world, including multiple government institutions in Russia, US delivery service FedEx, French car manufacturer Renault, German rail network Deutsche Bahn, Indian state police departments, and hospitals in Indonesia, Ireland, and the UK, with UK ``hospitals and clinics forced to turn away patients after losing access to computers'' \cite{bbc_wannacry_2017}. This release and subsequent disruption was so devastating that Microsoft President Brad Smith compared it to ``the U.S. military having some of its Tomahawk missiles stolen'' \cite{nakashima2017}. It was alleged that the disclosure was part of a bigger release by NSA contractor Harold T. Martin III \cite{nakashima2017}, with claims that he took over 50TB of data, and that he  ``took the documents and computer files and simply walked out the front door'' \cite{politico_martin_2016}; he was later sentenced to 9 years in prison for these offences  \cite{zapotosky_2019}.

By any measure, 50TB is a vast amount of data to exfiltrate. In addition, this is almost exactly the size of the lookup tables for the 10 primes that the NSA may have precomputed discrete log tables for when they allegedly cracked Diffie-Hellman key exchange (Section \ref{s:DHkey}); 10 primes at 6TB per table, giving 60TB. Thus, it is not unreasonable that such a database would also have the potential to leak or be exfiltrated, even by just one individual.

\subsubsection{Example: Mass storage facilities -- Utah Data Center}

In 2013, claims were made that the NSA was finishing off construction of the massive \textit{Utah Data Center} \cite{hill2013blueprints}, at a cost of over USD\$1billion, and with a capacity of around 10 Exabytes ($10^{19}$ bytes). This would have been a significant amount of storage for the time. In comparison, the  capacity required to store ``all words ever spoken'' (in written form) is estimated to be approximately 5 Exabytes \cite{nyt2003information}; half the capacity of the Utah Data Center. And with \textit{Moore's law} (the doubling of computer power every 2 years), such a facility might be reasonably expected to hold close to 100 times that amount in 2025, i.e., 1000 Exabytes (1 Yottabyte). Such data storage by a surveillance agency shifts the bar of \textit{mass surveillance}, to the point where on might refer to it more aptly as \textit{total surveillance}.

\subsubsection{Overview of deliberate privacy violations}

We informed students that debates about whether, and to what extent, state surveillance agencies should be able to surveil the public, are both deep and complicated.  Tying this back to the earlier example, we posed that mathematicians might view helping a government agency to surveil suspected terrorists as being noble work, but such work might be viewed very differently if these tools are then used to surveil the entire population en masse, even those who have `done nothing wrong' (a case of dual use). And we raised the idea that, even if the government of the time is considered `benevolent', a change in leadership could mean that the next government might not be considered so benevolent.

We also highlighted some common concerns throughout the examples given. Namely, that the sorts of tools such agencies use to carry out surveillance can have significant unintended consequences. We made it clear that mathematicians sit at the forefront of enabling such technologies and capabilities, and that different types of tools carry different risks. Students were able to see for themselves that well-resourced actors are perhaps the only ones capable of \textit{creating} such tools, but once created, these tools can then quite feasibly fall into the control of others, with potentially devastating consequences. We highlighted that the \textit{use} of such mathematical tools are not constrained by the \textit{intent} of their creators. In particular, working on tools that are `Hard to make, easy to lose, and devastating if lost' might warrant some careful consideration, as the notion of dual use is multi-dimensional. 

In addition, we demonstrated to students that the potential risks and harms from such tools are not constrained by the capabilities of today. For instance, the 10 Exabyte Utah Data Center contained far more data than any individual or group could ever read \textit{back in 2013} (it was primarily a database to `search through'). However, with the advent of Large Language Models, the NSA now has the tools to read, and infer from, this \textit{entire} database. Dual use considerations not only need to factor in where the tool is used and who is using, but also \textit{when} it is used. In particular, new mathematical work (LLMs) can make old mathematical work (large surveilled and stored datasets) even more risky. We emphasised that every system that mathematicians design, every algorithm that they write, and every exploit that they find that allows them to surveille our intended target, can also be used to surveille unintended targets, perhaps by unintended users. Mathematics does not discriminate and filter between `right' and `wrong' uses; if mathematicians develop it, and it might be harmful, then one day that harm may be realised, in particular if it falls into adversarial hands.

\subsection{Mathematicians indirectly violating privacy}

We then moved on to examples of mathematical work where privacy violations were an indirect consequence of the work being carried out. That is, the intention was not to violate privacy, but the systems created would, by design, result in privacy violations nonetheless. At this point we tried to  emphasise that just because mathematical work is \textit{intended} to do a good thing, does not mean that every aspect or consequence of it is necessarily good, and that mathematicians should think very carefully about all the secondary effects of the systems, tools, and capabilities that they are creating.

\subsubsection{Example: Business models of data collection -- surveillance capitalism}\label{s:surv-cap}

In recent years the term \textit{surveillance capitalism} has been introduced to describe the commodification of (personal) data \cite{zuboff2019surveillance}. This highlights the idea that companies can provide some sort of product or service that is `free' to the end-user (e.g., access to a digital platform for posting content), but then make commercial use of the content that user provides in order to sell various services to other entities (e.g., analytics allowing them to target certain adverts to certain individuals based on perceived preferences inferred from that data). 

The number of companies running such business models are numerous; these are the typical `big tech' platform-based companies that spawned around the turn of the 21st century (Google, Facebook, Snapchat, YouTube, Twitter, etc.); slightly distinct from the software- and/or hardware-based companies originating in the 1970s such as Microsoft and Apple. These new companies run a business model whereby users provide `valuable' data. This might be vast data in the form of behaviour and preferences (e.g., search terms, interaction with platform content such as which items are viewed or clicked on), or large files such as videos or photos, or vast amounts of written text through communications (e.g., messages, emails, etc), or any other `big' data. But this can also extend to data that is not very `big' in terms of number of bytes, but nonetheless valuable. This could be in the form of labels and tagging (e.g., data on when/where a user went, who they interacted with, identifying individuals in a photo, etc.), buying or reading habits (e.g., goods purchased and in what combination, content read or watched, etc.), or even `sensitive' data (e.g., data from dating/matching platforms demonstrating romantic preferences). 

Users of these platforms and systems can often be unaware of the extent to which the data they provide is being used and commercialised, either by the organisation itself, or by third parties. And it is not uncommon for the terms and conditions of such systems to be sufficiently long, complicated, or difficult to understand that most users simply accept them without meaningful understanding of the implications. For example, many people uploaded photos of themselves, and their young children, to the platform Flickr in the mid 2000s. 14 years later, those same photos became part of a large AI training database called \textit{MegaFace}, and have been used to train numerous AI facial recognition algorithms. Those children, now grown up, lament the fact that their photos form part of such databases, and the friends and family who uploaded them are full of regret \cite{hill2019flickr}. However, at the time of upload, no-one was any the wiser.

\subsubsection{Example: Fitness app heatmap -- Strava}\label{s:strava}

Strava is an online fitness data platform, built to allow its users to upload data from their fitness tracking apps or devices, ostensibly to compare with other users and thus `gamify' fitness activities with journey time, leader boards, etc. In 2017, Strava released a (pseudo)-anonymised global \textit{heatmap} of all the journeys taken by its users \cite{burgess2018strava}. This drew a `line' on a map for each journey tracked; the more journeys over a single path, the `hotter' that path appeared on the map. They `anonymised' the data by removing usernames from the heatmap. The database did not contain any usernames associated to the paths, but nonetheless a lot was able to be inferred from it. This included the ability to identify secret military bases and ``black sites'' across the world, potentially placing military personnel at great risk \cite{stewart2018strava}. This not only included the location of the sites \cite{stroud2018fitness}, but moreover some of the internal roads and paths within them; this extended to the internal pathways within top-secret sites such as ``Area 51'' in the USA and the Trident nuclear weapons facility in the UK \cite{burgess2018strava}. This was a case of a platform operating as intended, but with unintended consequences. Indeed,  Strava's CEO James Quarles stated the company was ``reviewing features that were originally designed for athlete motivation and inspiration to ensure they cannot be compromised by people with bad intent'' \cite{stewart2018strava}. 

Unfortunately, the 2017 heatmap release was not an isolated incident of sensitive information being released into the public domain through Strava. In 2024 it was revealed that members of the security teams of various world leaders were uploading their fitness data to Strava, such as jogging routes that started and ended at the hotel they (and the leader they were protecting) was staying at, in advance of arrival of the leader. These included members of the security teams of the Presidents of the USA (Joe Biden),  France (Emmanuel Macron), and Russia (Vladimir Putin) \cite{independent2024strava}. And in 2025, it was revealed that the location and patrol schedule of several French nuclear-armed submarines over the past decade could be determined from user data uploaded to Strava \cite{dailymail2025strava}. 

In all these cases, there was no `malfunction' of the Strava platform; user data was uploaded, shared, and displayed, as intended by Strava, with privacy and sharing settings chosen by the users. However, the construction of such a platform, along with the visibility and sharing settings allowed to users on it, meant that their privacy (and the privacy of others) was compromised, in ways that could have had devastating consequences. And a deferral to complicated terms and conditions, and privacy settings, was made by Strava, who ``basically put the onus on users'' for the data releases, even though the platform had 7 different privacy settings, some of which could only be modified in the browser version of the platform (not the app), and all enhanced privacy modes were ``opt in'' by default \cite{bonazzo2018strava}. Strava response was that users could do things like ``specifically opting out of the global heatmap by unchecking a box in the settings page'' \cite{hern2018strava}. 

It is not clear exactly why Strava chose to produce, and continue to maintain, their global heatmap. However, they did declare that ``Our global heatmap is the largest, richest, and most beautiful dataset of its kind. It is a direct visualisation of Strava's global network of athletes.'' \cite{stroud2018fitness}, suggesting some sort of appeal to `mathematical beauty' that will be covered later in this course in Lecture \ref{s:community}.

\subsubsection{Overview of indirect privacy violations}

We concluded this section by firstly reminding students that, as mathematicians, they may well end up working for such organisations, designing the functionality of these platforms and systems. As such, they should be mindful of what they create, and the harm it might inadvertently cause, even with the best of intentions by those developing it. We emphasised that in creating convenient solutions and tools for users, mathematicians might also be creating problems (in the form of privacy violations) that otherwise would not have existed. Some of these might do a lot of damage, even though the developers didn't plan it that way. And the damage could be very long-lasting. We drew a comparison here to the use of asbestos throughout the 20th century. Asbestos was once considered a very useful and convenient building material. It was cheap, easy to work with, durable, and hailed as a `miracle solution' to lots of different construction problems. However, it turned out to be highly carcinogenic. Moreover, the removal of asbestos is prohibitively expensive, if not impossible. We suggested to students that `big data', and the platforms that encourage, host, and share it, might become the asbestos of the 21st century: convenient for users, providing an affordable and effective service, but potentially causing vast harm in the long run, and near-impossible to remove later, as those whose images are on MegaFace now know.

\subsection{Mathematicians violating privacy out of ignorance}

We started off by highlighting that mathematicians can also design systems that (very unintentionally, but also unsurprisingly) violate privacy, and thus that it is the responsibility of mathematicians to respect and protect the privacy of those who provide their data. We put forward the idea that designing digital systems to collect people's data also requires significant consideration towards how well it is protected, to ensure it is kept safe. And we made students aware that, just because mathematicians can do mathematics and understand cryptographic algorithms, does \textit{not} necessarily mean they understand how to properly encrypt data and keep it safe. We pointed out that data protection is only as strong as its weakest component, and illustrated this by highlighting the importance of salting and hashing passwords \cite{stallings2003network}.

Overall, we used this section to emphasise that data safety is not a \textit{mathematical} problem, but rather a \textit{practical} problem; one for which mathematics can certainly help, but does not provide all the solutions to. To emphasise the `practical' nature of this, we pointed out that the best way to protect data is to not store it in the first place. We exemplified this with some simple examples: if an organisation is running a website registration system, why would they need the users full name, date of birth, address, etc.? If an organisation is never going to post them anything, collecting their address poses unnecessary additional privacy risks, and should be avoided. We acknowledged that collected information has value to advertisers etc.~(as mentioned in Section \ref{s:surv-cap} with surveillance capitalism), but that this value comes at a cost to the end user. These seemingly `non-mathematical' approaches to privacy preservation came as an initial surprise to the students, but they quickly saw the value and importance of them.

From there, we moved on to give examples of where mathematicians had built systems and/or carried out activities where things went horribly wrong, leading to massive privacy violations.

\subsubsection{Example: Large data breaches}

With the advancement of digital systems and online platform, mathematicians now often work on projects that involve the  private data of millions, or even billions, of people. Large-scale hacks and data dumps are a recurring, escalating problem, affecting even some of the largest (private, and government) organisations in the world; this list is dynamic and ever-increasing\footnote{A partial running list can be found at \url{https://en.wikipedia.org/wiki/List_of_data_breaches}}. Significant breaches have occurred from organisations like Equifax (143million data records \cite{cnbc2017equifax}), Ticketmaster (560million customer account records \cite{Lawler2024ticketmaster}), the Shanghai National Police (1billion citizen records \cite{ni2022shanghai}), and Yahoo (3billion user accounts \cite{perlroth2017yahoo}). Even Experian, a large credit reporting company who offer professional services to assist with data breaches \cite{experian2026databreach}, suffered a data breach itself, losing 15million customer records \cite{Thielman2015experian}.

In all these instances, the release of the data set was \textit{not} intended by the organisation. However, with contributing factors such as poor security, misconfigured servers, lost storage media, hacking, inside jobs, or even accidental publication\footnote{Some of the reasons give in \url{https://web.archive.org/web/20260809075657/https://en.wikipedia.org/wiki/List_of_data_breaches}}, such releases inadvertently occurred, leading to privacy violations at scales that were unfathomable before the advent of the internet and digital technology. Unfortunately, even though those affected might ostensibly have the protection of the law, the law cannot always meaningfully protect them. Sometimes the legal system chooses not to, for reasons of limited resources such as when the Leicestershire police chose to only investigate burglaries if there were at even-numbered houses, to save on costs \cite{bbc2015leicestershire}. Other times, there might not be a `police force' at all, such as when the UK Ministry of Defence accidentally leaked the details of thousands of Afghans fleeing the Taliban in 49 separate data \textit{breaches} \cite{kenber2024}. The consequences here are not merely hypothetical or abstract;  recent research suggests that 49 people were killed by the Taliban as a result of these Ministry of Defence breaches \cite{taylor2025afghan}.

\subsubsection{Example: De-anonymisation hackathon data release}\label{s:ptv}

As outlined in \cite{ovic2019myki}, in July 2018 the meetup group \textit{Data Science Melbourne} held the \textit{2018 Melbourne Datathon} in the state of Victoria, Australia. As part of this, they requested, and were granted access to, a large dataset from Public Transport Victoria (PTV). This dataset was from the \textit{Myki Card} system used to pay for public transport journeys in Victoria; similar to the Oyster Card system in London. The dataset contained 1.8billion historical ‘touch on' and ‘touch off' event records of public transport users' activity from 15.1million Myki cards between 1 July 2015 and 30 June 2018 \cite[p.~20]{ovic2019myki}. Before release, PTV had ``substituted each card number with randomly generated card number'', which they claimed  ``anonymised the individuals'' \cite[p.~21]{ovic2019myki}; a claim later disputed. Specifically, a touch on/off event contained the following data (quoted from \cite[p.~20]{ovic2019myki}):

\begin{itemize}
    \item date and time - the date and time of the ‘touch on’ or ‘touch off’ event, to an accuracy of one second;
    \item location information — this data varies depending on the mode of transport, but includes information such as vehicle identifiers, route numbers and stop numbers;
    \item card identifier — the unique number assigned to each Myki card when preparing the dataset for release. The card identifier permits all ‘touch on’ or ‘touch off’ records within the dataset relating to a particular card to be linked; and
    \item card type — a descriptor of the type of card used. There are approximately 70 Myki card types, including card types specific to Victoria Police, Federal Police, State and Federal Parliamentarians, asylum seekers, veterans and pensioners.
\end{itemize}
PTV released the dataset to Data Science Melbourne with absolutely no restrictions, and thus the Datathon website informed participants  ``No NDA [non-disclosure agreement] to sign this year -- you can do what you like with the data'' \cite[p.~13]{ovic2019myki}. It was freely available on the Datathon website for the duration of the Datathon (24 July -- 26 September 2018), and one  participant republished the dataset in
full on their blog for an additional 3.5 months \cite[p.~24]{ovic2019myki}.

During the Datathon, some privacy researchers discovered the data set, and went about investigating if individuals could be de-anonymised \cite{Culnane2019}. By comparing their own travel history (which they had access to in their Myki accounts) to the database, the researchers were quickly able to identify their own Myki cards in the database, simply by searching for 2 distinct touch on/off events \cite[p.~3]{Culnane2019}. They were also able to identify the Myki cards of some of their known co-travellers (friends of theirs) in various journeys, having first sought permission from them to do so \cite[p.~3]{Culnane2019}. Going further, they sought permission from a member of State Parliament -- Anthony Carbine, who tweets  frequently about taking the train -- to locate his Myki card. Having analysed his twitter account, they found 3 tweets describing 3 separate train station visits which, when compared to the database, yielded only one possible Myki card matching those date and (approximate) time visits to those stations; they had, with very high probability, de-anonymised Mr Carbine's Myki, using only publicly available data \cite[p.~4--5]{Culnane2019}. The went on to conclude that it would take about 3 data points to identify the Myki of \textit{any} random stranger; less if that person had co-travelled with the person doing the analysis \cite[p.~5]{Culnane2019}. They also raised concerns around the potential for harm this might present, especially seeing as the (released) database included the card type, which listed sensitive or high-risk demographics such as the police, parliamentarians, or children \cite[p.~4]{Culnane2019}.

PTV carried out a Privacy Impact Assessment (PIA) prior to the data release. It was  ``completed by a PTV data scientist and signed off by the data owner and PTV Chief Information Officer. It was then sent to PTV’s governance and legal team by the data scientist'' \cite[p.~33]{ovic2019myki}. The legal team, relying on the PIA, approved it, stating  ``Given that there is no personal or confidential information involved in this project, privacy or confidentiality laws are not applicable to this project'' \cite[p.~33]{ovic2019myki}. However, the PIA was later found by the Victorian Information Commissioner to be ``flawed in a number of respects'' \cite[p.~34]{ovic2019myki}, that the ``PIA created a false sense of security'' \cite[p.~35]{ovic2019myki}, and that ultimately ``The flawed approach taken by PTV in completing the PIA, and the resultant confidence it provided to PTV and DPC that privacy issues had been considered and addressed, are key contributors to this incident'' \cite[p.~35]{ovic2019myki}. On the back of all this evidence, the Victorian Information Commissioner found that the PTV had carried out a serious disclosure of personal information, and ordered them to reform their data handling processes or face a fine of approximately \$500,000 AUD (3000 penalty units, where a unit at the time was \$165.22 \cite{sentencingcouncil2026fine}).

\subsubsection{Overview of violating privacy out of ignorance}

We used the examples in this section to emphasise the responsibilities mathematicians had to properly protect any data they hold or work with, thus leading into some of the safety aspects of the design and development process of mathematical work. This section helped students see that being able to design a database or system  to \textit{collect and hold} data is not the same as being able to \textit{secure} said database or system, and that if they can't keep it safe, then perhaps they should reconsider how they are doing it, and whether they should even be starting it in the first place. This idea of `If you can't do it safely, you shouldn't be doing it at all' is a core tenent of ethics in mathematics, and one that by now the students had started to appreciate after seeing many, many examples.

\subsection{The importance of privacy}

\subsubsection{The danger of false privacy}\label{s:false}

Drawing together all the examples given in the lecture, we explicitly articulated to students that the only thing more dangerous than no privacy, is false or misplaced privacy. We gave additional specific examples, of how the Germans and the Japanese learnt this lesson (the hard way) in World War II; they transmitted \textit{everything} thinking their cryptosystems were unbreakable -- the Germans with their \textit{Enigma machine} \cite{hinsley2001codebreakers} and Japanese with their \textit{Type B ``Purple'' machine} \cite{clark2011man}. This was of significant detriment to their war effort as the Allies had decrypted most of their communications. We also used a historical example of false privacy, the \textit{Babington Plot}, whereby Mary Queen of Scots was executed for plotting to overthrow and execute Queen Elizabeth I of England; Mary's `encrypted' letters to her co-conspirators (including Anthony Babington) were decrypted as the cipher being used had been provided by Elizabeth's spymaster Francis Walsingham \cite[pp.~21--22]{dooley2013brief}. We conveyed the idea that if someone \textit{knows} they can be (or are being) surveilled, they can (and often do) self-censor, and are more careful with what they communicate. However, if someone thinks their communications are secure, they then tend to communicate more details, some of which may cause them great harm if intercepted. We summarised this by saying that lulling people into a false sense of security is the greatest possible violation of privacy; in this case, people \textit{volunteer} information that they otherwise would have never have revealed. Such `lulling' is an act that mathematicians can, and often do, carry out with their mathematical work (be it intentionally, unintentionally, or out of ignorance), and the students were able to make the connection between that idea and the myriad of examples given in the lecture.

\subsubsection{Mathematicians destroying the notion of privacy}

We ended the lecture by illustrating thy responsibilities mathematicians have to respect privacy. To do this, we articulated that data should not be viewed and treated as `just a bunch of $1$'s and $0$'s on a computer'; rather, we suggested to students that when working with data sets, it may help to consider these not as `data', but as \textit{information about people}, because most of the time that is exactly what is being worked with. For example, mathematicians don't work with `shopping data'; they work with \textit{information about how people buy things}. They don't work with `medical data';  they work with \textit{information about people's health}. And they don't work with `financial data'; they work with \textit{information about how people use their money}. We informed students that the data they work with and manipulate as mathematicians is almost always information about people, and to treat it as such. Overall, we suggested to students that, rather than using the term `data', they replace it with the term `people's personal information', as most of that time that is what it is. This simple change in language 
can have a significant effect on how mathematicians think about, and use, data. Students quickly saw, and felt, the difference between phrases such as `My platform has collected 10 Terabytes of data' and `My platform has collected 10 Terabytes of people's personal information'.

We put it to the students that mathematicians create privacy \textit{preserving} tools such as cryptography, database security, etc. We also put to them that mathematicians create privacy \textit{subversion} tools, by breaking cryptography, setting up vast data collection systems, not keeping data secure, etc. And so, we pointed out that this is not a zero-sum game; the public would be much better off knowing there was no privacy, rather than being fooled into thinking there was, because at least then the public could take measures into their own hands and restrict what the say, what data they provide, and to whom they provide it. We finished by asking students to imagine a world where there's no privacy whatsoever; that's the world mathematicians are creating, and whether or not they personally assign much value to privacy, they would be wise to think about the consequences of such a world.

\begingroup
\let\section\subsection
\renewcommand{\bibname}{Bibliography (Lecture 3)}
\renewcommand{\refname}{Bibliography (Lecture 3)}
\addcontentsline{toc}{subsection}{Bibliography (Lecture 3)}
\renewcommand{\addcontentsline}[3]{}
\bibliographystyle{IEEEtran}
\bibliography{bibLec3}
\endgroup

\newpage
\section{Lecture 4: Fairness and impartiality in algorithms and AI}\label{s:algorithms}

\fbox{
  \parbox{0.98\textwidth}{
Mathematicians design the algorithms that run the world, deciding what people read, who they date, and when then can get a loan, often done under the auspice of increasing efficiency and removing human subjectivity. The scale, and speed, at which such algorithms are deployed means that their impact, and potential for harm, can swell far beyond original intentions.

This lecture looks at how mathematicians can understand the ways in which their algorithms rank, predict, suggest, and rate things for and about people and society, and how this might lead to harm. It shows how they can see the subjectivity of their work as a strength when interpreted correctly, as well as how to avoid any large-scale runaway behaviour and subsequent harm. 
}}

\subsection{Impartiality}

We stared the lecture by informing students that algorithms now run a substantial proportion of the important processes in the world, and that it is mathematicians who design them. From  deciding what people read, to what they buy, and whether they can get a loan, mathematicians often design these systems to \textit{remove} human subjectivity from decision making processes and to make them `impartial', as is being done with predictive policing and prison sentencing algorithms. And so we started to probe students to question how impartial, or fair, a system designed by humans can ever be?

One of the main purposes of this lecture was to dispel the `myth of impartiality' surrounding mathematical work; just because a process is mathematical, does not mean it is impartial.  Throughout the lecture, we used examples to show that mathematical processes can even have \textit{more} bias than human-run processes in certain circumstances. To start understanding why mathematics is not impartial, it was first necessary to see why mathematicians might think it is. This was done with the use of an interactive activity at the start of the lecture.

\subsubsection{Activity: Impartial plagiarism checking algorithm}\label{s:impartial-algo}

We divided up the class into groups of approximately 5, and each group was tasked with the following: \textit{Describe how to make an impartial algorithm to check for plagiarism in written essays}. Strong emphasis was put on the term \textit{impartial} as being a core part of the task. The groups were given 5 minutes to come up with ideas, and then a spokesperson from each group was given 1 minute each to present their algorithm. They key steps or points for each group were recorded on the board. These included steps (aggregated from across the groups) such as:

\begin{itemize}
    \item \textit{Decide} what is meant by plagiarism.
    \item \textit{Choose} a dataset to compare essays to.
    \item \textit{Set} a string length $N$ and look for matches of that length.
    \item \textit{Pick} some $K$ and identify the top $K$ most similar passages between the texts.
    \item Do a \textit{Google search} on the text.
    \item Have a \textit{human} check any matches to make a final determination.
 \end{itemize}
The purpose of this was to tease out aspects of the algorithms, which may \textit{seem} impartial as they are carried out by a mathematical process (e.g., matching strings of a given length), but actually contain some degree of bias from the choices of if, and how, to carry them out in the first place. For example, students were asked follow up questions such as ``Whose definition of plagiarism should be used?'', or ``Why are string matches of length 20 the right length to search for?'' or ``What makes you think the Google algorithm will be impartial?'' This was supplemented by deeper questions such as ``What happens if someone knows how your algorithm works?''; something students initially suggested wouldn't make a difference, until it was pointed out to them that plagiarising from a different language and translating, or being sure to `break' all overlapping strings at length $N-1$ with superfluous words, might completely nullify their algorithms. Students saw from these overall observations that, actually, their algorithms were not impartial, but rather, they were laden with value judgements and human opinions. We used this as a springboard to start giving deeper examples of where mathematics had been used in ways that were purported to be impartial, but in practise were far from it. Having established that mathematical processes \textit{can} be biased, we moved on to show that they are \textit{unavoidably} biased.

\subsection{Ranking algorithms}

\subsubsection{Example: Amazon ``Buy it Now'' algorithm}

The online platform Amazon has a ``Buy it now'' button for most product on its webpage, where Amazon `selects', in some automated/algorithmic way, one seller and places their listing in the ``Buy it now'' box. One might \textit{think} that it simply chooses the cheapest seller, perhaps factoring in postage. But in actual fact, the algorithm is substantially more complicated than that, and ``takes into account many data points and layers of information gathered about the marketplace, current market conditions, and every seller'',  to choose a seller in a way that is ``determined by a complex, machine-learning algorithm that takes millions of data points into consideration'', all designed and computed by Amazon \cite{Feedvisor2025}.

Customers make their purchase decisions base on Amazon algorithms. Thus, these algorithms are determining the price customers pay, and the service they receive. But customers are often want a quick suggestion as to which seller to buy from, and Amazon takes advantage of this fact. Amazon itself claims its ``Buy it now'' algorithm is `neutral', and in a testimony give to the US Congress it stated ``Our algorithm such as the buy box is aimed to predict what customers want to buy, and we apply the same criteria whether you’re a third-party seller or Amazon to that because we want customers to make the right purchase regardless of whether it’s a seller or Amazon'' \cite{Mitchell2025}. However, not only does Amazon run the ``Buy it now'' box itself, in a very opaque way, but they manipulate it in their own favour in a manner that is completely invisible to the average consumer, such as ``artificially giving itself perfect customer scores for every product it sells'', and by ``giving its own warehousing and shipping service, Fulfillment By Amazon (FBA), a perfect score on the metrics related to fulfilling an order, such as the on-time delivery rate and depth of inventory'' \cite{Mitchell2025}. By being able to manipulate the inputs of the algorithm, Amazon was thus able to manipulate the algorithm's outputs, and therefore its overall effect. It is Amazon's unique position as both the \textit{host} of the marketplace, and a \textit{trader} on that marketplace, that makes this possible. If an entity can make the rules for the game that is played, they they are always bound to `win' (in this case, at the expense of consumers). And most customers would be totally unaware that this is happening.

\subsubsection{Example: Google PageRank algorithm}

In 1998 the original Google PageRank algorithm was proposed by Sergey Brin and Lawrence Page \cite{brin1998anatomy} to `rank' search queries on the internet and thus present a web user with an ordered list of webpages matching their search. This algorithm in its original form works as follows: Each webpage $w_{i}$ on the web (of $n$ pages) is assigned a value $v_{i}\geq 0$ such that, if $L_{i}$ is the set of pages that link to $w_{i}$, and $n_{i}$ is the number of outgoing links from page $w_{i}$, then
\[
 v_{i}=\sum_{w_{j}\in L_{i}}\frac{v_{j}}{n_{j}}
\]
That is, each page $w_{i}$ `donates' $\frac{1}{n_{i}}$ of its value $v_{i}$, uniformly, to each page that it links to. This leads to a `normalised' adjacency matrix $A$ with value $1/n_{i}$ in each entry $(i,j)$ where $w_{i}$ links to $w_{j}$; there should be $n_{i}$ such nonzero entries in each column, summing to 1, and this matrix can simply be read off by counting all the page links on the internet. At this point, one simply needs to find a column vector $v=(v_{1}, \ldots, v_{n})$ of `values' such that $Av=v$; that is, an eigenvector of $A$ with eigenvalue $1$. However, such an eigenvector might not be unique, especially if the web has disconnected components. To deal with this, and produce a `unique rank', Brin and Page decided to convert $A$ to a \textit{positive} matrix (one where all entries are $>0$) by use of what is known as the PageRank \textit{damping factor} $d$; this can be any real value $0<d<1$, and is ``usually set to $d=0.85$'' \cite[p.~109]{brin1998anatomy}. One can then define the $n \times n$ matrix $B$ with all entries $1/n$, and then computed the final `damped' adjacency matrix $C:= dA + (1-d)B$. The columns of $C$ also sum to 1, and $C$ is positive. It follows by the Perron-Frobenius Theorem that $C$ has an eigenvalue of 1 corresponding to a 1-dimensional eigenspace, so one can find a unique eigenvector with eigenvalue 1, whose entries sum to 1 \cite{Tanase2009}. This then gives a `value' for each webpage, and these can be ordered to give a rank of search results. Moreover, by starting with a vector $v$ which is a column of 1s, and successively taking powers $Cv, C^{2}v, \ldots$, by the Power Method Convergence Theorem this (quickly) converges to an eigenvector of $C$ with eigenvalue 1 \cite{Tanase2009}. And so Brin and Page had developed a computationally efficient method to give a unique ranking to webpages matching a given search term. In addition, the choice of $d=0.85$ is also somewhat arbitrary, without deep justification.
 
The motivation behind this algorithm was straightforward. Rather than using text based ranking systems (i.e., counting instances of the search terms on the page) which could be easily gamed by thousands of instances of a common search word, PageRank assigned importance of a page by how many other pages linked to it, normalised by how important \textit{those} pages were. Brin and Page interpreted the PageRank value (the normalised eigenvector computed above) as a `probability' that a websurfer randomly clicking links on pages will eventually reach that page; the damping factor $d$ was interpreted as a probability of the surfer's next move being a link click, and $(1-d)$ as a probability of them giving up clicking and choosing a `totally random' page instead \cite[p.~110]{brin1998anatomy}; how this totally random page is chosen is not described.

\subsubsection{Analysis of PageRank}

We highlighted to students that PageRank was not without its drawbacks. Just as text based ranking could be easily gamed, so to could PageRank, by setting up thousands of dud websites that all link to each other, boosting their PageRank `score'. In addition, PageRank did not reflect any aspect of the \textit{content} of  a webpage; an Encyclopaedia Britannica page on `giraffes' may well receive a lower rank than a simple webpage called `Dancing giraffe' containing only an entertaining video of someone dancing in a giraffe costume. The PageRank algorithm might seem reasonable, but is there always a `correct' algorithm to implement in any give case, scenario, or for any given problem such as website ranking? In cases of mathematical calculations, for example in calculating interest payments, this is true. But the choice of which algorithm to use (and how to use it) in a search engine ultimately affects the way billions of people interact with the internet. It might be \textit{effective} in achieving some sort of outcome, but is it \textit{impartial}? Recall that several human decisions were needed to ensure that the eigenvector above was unique; mathematics gives an answer, but only after humans set the question. If mathematicians decide what algorithm to implement (say, PageRank), and furthermore need to then decide how to tie-break certain outputs (to obtain `unique' rankings), then is this an \textit{impartial} way to rank pages on the internet?

One might assume that Google would always strive to implement an `impartial' search engine. However, in 2017 Google announced  that it was `deranking' various news sites \cite{BBC2017derank}; a move that those sites referred to as `a form of censorship'. If a ranking algorithm is objective and impartial, then how can a company \textit{choose} to derank websites? Lowering the ranking of a website sufficiently far means it becomes rather inaccessible to the public. An average Google search yields billions of hits, and the average Google user might go through the top 100 at most. Lowering a website to rank below 100 thus makes it inaccessible to most people; doing so would lead to the suppression of information, hence the accusations of Google's actions as `censorship'. Now perhaps Google has a valid and just reason to undertake such actions; there are many complicated issues at play. But by this point it is clear that such actions (i.e., `deranking'), even if carried out by an algorithm, are in no way impartial.  Indeed, Google was even fined EUR\euro2.4billion  ``for having abused its dominant position in several national online search markets by favouring its own comparison shopping service over those of its competitors'' \cite{ECJ2024}.

Ultimately, a search engine is simply a distortion machine, by design. It looks through trillions of web pages, picks out all the ones that it `thinks' are relevant, and then orders them in some `optimal' way. Therefore, it distorts what users see. This may still be \textit{useful}, but that does not mean it is \textit{impartial}.

\subsubsection{Overview of ranking algorithms}

We concluded this section by outlining to students that those who create ranking algorithms will have their own reasons for wanting to carry out ranking, which could go beyond the `greater good' and include their own self interests. We highlighted that what they see as important, as well as what is important to them, are two big factors that can influence how such algorithms are created.

\subsection{Predictive algorithms}

We then moved on to cover work that had  been done on algorithmic policing and algorithmic sentencing. We outlined how this might be quite appealing for those running the justice system; using a computer to influence where to send police or how long to sentence people for is seen as a cheap and efficient solution that `removes human bias' from decision making. However, the following examples put into question whether these algorithms can genuinely remove such bias?

\subsubsection{Example: COMPAS sentencing algorithm}

As reported in ProPublica \cite{angwin2016machine}: In 1989 statistics professor Tim Brennan and prison manager Dave Wells founded the (private) company \textit{Northpointe}. They decided to ``build a risk assessment score for the corrections industry'', which they named the \textit{Correctional Offender Management Profiling for Alternative Sanctions} (abbreviated to COMPAS). This (proprietary) algorithm purports to give a probability of a person re-offending, known as a \textit{recidivism score}. This was used extensively across the US criminal justice system to assist judges making decisions, from sentencing to parole. In theory, COMPAS was only supposed to be used to ``determine which defendants are eligible for probation or treatment programs'', but not as part of sentencing; even Brennan himself stated ``I wanted to stay away from the courts'' \cite{angwin2016machine}. But nonetheless, this software found its way into courts to influence sentencing decisions.

The COMPAS algorithm computes several types of risk score, including \textit{general recidivism} (any new offence within 2 years), violent recidivism (any new violent offence), and pretrial misconduct (and offence while awaiting trial) \cite[Section~4.1]{northpointe2015}. These take into account factors such as ``prior criminal history, criminal associates, drug involvement, and early indicators of juvenile delinquency problems'' (for general recidivism), ``current charges, pending charges, prior arrest history, previous pretrial failure, residential stability, employment status, community ties, and substance abuse'' (for pretrial misconduct), and ``history of violence, history of non-compliance, vocational/educational problems, the person's age-at-intake and the person's age-at-first-arrest'' (for violent recidivism). Indeed, the Practitioner's Guide to COMPASS \cite[Section~4.1]{northpointe2015} gives that the violent recidivism score is calculated as simple weighted sum of those 5 variables:

\begin{align*}
    \textnormal{Violent Recidivism Risk Score}\ =\  &\textnormal{(age)}* w_{1} 
    \\&+ \textnormal{(age-at-first-arrest)}*w_{2} 
    \\&+ \textnormal{(history of violence)}*w_{3}
    \\&+ \textnormal{(vocation education)}*w_{4}
    \\&+ \textnormal{(history of noncompliance)}*w_{5} 
\end{align*}
How these weights are calculated, or how the other recidivism scores are calculated, is not explicitly given. However, the variables are based only on the answers given in the 137-part questionnaire given to the accused at the onset (where 131 of these questions are multi-choice, and the remaining 4 are non-negative integer responses). The \textit{Practitioner's Guide to COMPAS}, with  these, and additional, details is available online \cite{northpointe2015}, in which it claims that such methodology is ``an objective method of estimating the likelihood of reoffending'' \cite[p.~29]{northpointe2015}.

ProPublica analysed some COMPAS scores, in particular, data from Broward County, Florida \cite{angwin2016machine}. They found that COMPAS risk scores correctly predicted recidivism $61\%$ of the time, similar to the Northpointe claim of $68\%$. However, ProPublica also found two additional properties of that data set, by analysing the risk scores (which were integers between 1 and 10):

\begin{enumerate}
    \item Scores for white defendants were skewed toward lower-risk categories, with almost $30\%$ receiving a low score of 1/10, and under $3\%$ receiving a high score of 10/10. Scores for black defendants were not, and were almost uniformly distributed between 1 and 10.
    \item The false positive rate (i.e., labelled higher risk, but didn’t re-offend) was lower for white defendants ($23.5\%$) than for African-Americans ($44.9\%$). And the false negative rate (i.e., labelled lower risk, yet did re-offend) was higher for white defendants ($47.7\%$) than for African-Americans ($28.0\%$). 
\end{enumerate}
From this, ProPublica concluded that the COMPAS algorithm was ``biased against blacks'' \cite{angwin2016machine}, which Northpointe later disputed \cite{northpointe2016}. BBC News \cite{BBC2016maths} summarised the issue as follows: ``But all sides in this statistical showdown do seem to agree on one thing. If black people in general are more likely to re-offend, then a black defendant is more likely to be given a higher risk score. So while the algorithm itself may not be racially biased, it is reflecting racial biases in the criminal justice system and society more widely.''

\subsubsection{Example: PredPol predictive policing algorithm}

As reported in \cite{johnson2024predictive}, the notion of \textit{predictive policing} entails ``using quantitative techniques to identify likely targets for police intervention and prevent crime'' \cite[Section 1]{johnson2024predictive}; that is, to mathematically predict where to send police \textit{before} crime has occurred. A company called PredPol (an abbreviation of ``Predictive Policing'') created and patented an algorithm of the same name to carry out predictive policing. It sought to optimise resources used by a police force by predicting hotspots of crime, by taking in historical and current crime data, and identifying where police departments should deploy resources \textit{to} (and where to pull resources \textit{from}). The algorithm relies on a mathematical model with two components.

The first component is the \textit{Epidemic-type aftershock (ETAS) model} \cite[Section 2.2]{johnson2024predictive}, treating crime as a ``continuous time, discrete space epidemic-type aftershock sequence point process''; the same process in seismology to predict aftershocks from parent earthquakes. This describes crime as a branching process, with initial crime occurring as a random (i.e., Poisson) point process, but then subsequent crime occurring locally as `aftershocks', which may continue for several iterations. That is to say, the mathematics formulated to describe the (highly physical) movement of tectonic plates is also used to describe the (highly sociological) occurrence of crime.

The second component is the \textit{reaction-diffusion model} \cite[Section 2.1]{johnson2024predictive}, which is used to describe how visible instances of crime \textit{attract} more crime to an area. This model is ``typically used to describe chemical reactions, in which activators and inhibitors move, mix, and interact''; in PredPol, criminals are interpreted as activators, and law enforcement as inhibitors. The specific way this model is used in PredPol relies heavily on a core assumption of the \textit{broken windows effect}; a (disputed) concept in criminology stating that visible instances of crime (such as broken windows, graffiti, etc.) in a geographic area convey a perception that order is breaking down in that area, thus emboldening people to commit more crimes in or around that area \cite{thacher2003order}. This component models people being drawn in to a crime hotspot, to then commit more crime. Mathematically, it is given by solving the following coupled differential equations \cite[Equations 2--3]{johnson2024predictive}:
\[
\frac{\partial B}{\partial t} = \frac{\eta D}{z} \nabla^2 B - \omega B + \epsilon D \rho A \ , \ \ \frac{\partial \rho}{\partial t} = \frac{D}{z} \nabla \cdot \left[ \nabla \rho - \frac{2\rho}{A} \nabla A \right] - \rho A + \gamma
\]
However, as Ensign et al. \cite{pmlrv81ensign18a} describe, such a system can introduce feedback loops. The parameters used in the PredPol algorithm are determined by the crime data from the preceding 180 days; this is constantly updated. And crime data comes in two forms: \textit{reported incidents} (directly reported by the public), and \textit{discovered incidents} (those observed by police) \cite[Section 1]{pmlrv81ensign18a}. Unlike the former, the latter of these is heavily influenced by the police presence in an area; more police = more observations of crime, even if the crime rate itself is unchanged. As PredPol included discovered incidents in the crime dataset, Ensign et al. \cite{pmlrv81ensign18a} show that as soon as there exists non-uniform crime distribution across an area, the use of PredPol leads to crime hotspots become `hotter' over time, and cold spots become colder. And as the initial datasets used for PredPol may contain some bias towards certain groups, Johnson et al. \cite[Section 4.1]{johnson2024predictive} argue for example that ``due to the focus of predictive policing on crime data that is skewed towards communities of color in poor neighborhoods, the biases that we have discussed above in prediction will lead to overpolicing in those areas''. That is, rather than remove bias in policing, an algorithm such as PredPol not only calcifies such bias, but also enhances it. 

Further claims of such feedback loops creating ``self-fulfilling prophecies'' are made in Haskins \cite{haskins2019academics}, observing that the reason such a mathematical model works so well for earthquakes, but reinforces bias and creates feedback loops in policing, is that ``you have seismographs everywhere -- wherever an earthquake happens, you’ll find it''. By comparison, the measurement of crime is highly depended on the placement of (limited numbers of) police officers; the model thus impacts how crime is measured. On the back of that report (and other works), a petition of over 1500 mathematicians \cite{aougab2020web} was put together to ``call on the mathematics community to boycott working with police departments'' \cite{aougab2020boycott}.

\subsubsection{Overview of predictive algorithms}

We went on to highlight that there is much work being done now on algorithmic policing and algorithmic sentencing. This is, of course, quite appealing for those running the justice system; using a computer to decide where to send police  or how long to sentence people for is seen as a cheap and efficient solution, and moreover is believed to remove human bias from decision making. However, students saw from these examples that bias can certainly remain, or even increase.

We pointed out that all these seemingly-useful algorithms can be designed with the best of intentions, but then easily re-assigned for another more serious use that the developers may not have intended. Moreover, the perceived \textit{authority} of these algorithms is almost unquestionable; a judge may interrogate a prosecuting or defending lawyer, or a witness, but is highly unlikely to interrogate a computer. And even if they wanted to, how would they do so? The algorithmic decision assistance processes can easily become a decision \textit{making} process, as decision makers themselves are unlikely to challenge the algorithm; this is known as \textit{automation bias} \cite[p.~5]{chiodo2025formalising}.

We highlighted to students that, when designing and training these sort of predictive and decision-making algorithms, they might make use of large existing data sets to do machine learning and determine parameters. In doing so, this absorbs the existing cultural and societal bias that contributed to the creation of these data sets (e.g., who was caught committing crime, who was sent to jail, etc.) into these algorithms via the machine learning. Making these processes algorithmic does not wash away the bias, quite the contrary; it calcifies the bias, carving it into stone in the form of algorithms that are much harder to challenge and scrutinise than human decision makers.
 
We suggested to students that if they, as mathematicians, might be involved in the development of algorithms that determine things as important as police deployment and criminal sentencing, they need to be aware of the drawbacks of such systems. Just because a computer does it, does not make it impartial or free from bias.

\subsection{Suggestion algorithms}

Next, we moved on to suggestion and recommendation algorithms. This was perhaps the most counter-intuitive section of the lecture, with students seeing algorithms that simply \textit{suggest} options that users \textit{might} want to view, visit, or meet, as being unmitigated tools for good. However, as the following examples demonstrate, such `recommenders' can still have harmful effects.

\subsubsection{Example: Content recommender algorithms -- YouTube} 
 
Content-sharing platforms run (large, complicated, and opaque) algorithms to suggest content alongside whatever content  a user is currently watching (e.g., videos). These suggestion algorithms are not random, nor do they solely look for whatever content is `most similar' to the one currently being viewed; instead, they look to suggest content that the viewer is most likely to view next. These algorithms are often not optimised over what is most relevant, but rather over what will keep users viewing for the longest. Of course, this seems like a sensible strategy for such platforms, as many make their money from adverts. However, it has had the unfortunate side effect that the sorts of content that keeps people watching the longest sometimes tend to be unpleasant. Indeed, such recommendation algorithms have been shown to suck people ``down the rabbit hole'' and in to ``vortices'' of often rather radical content, such as far-right extremism, as has been reported with YouTube \cite{Roose2019}. Content recommendation algorithms, ever-present on social media platforms, ``can quickly spiral and take young people down rabbit holes of harmful and distressing material through no fault of their own''  \cite{Boyd2025}. But, through \textit{morbid curiosity}, it is precisely this distressing content that evokes the most human interest \cite{Oosterwijk2017, Scrivner2021}. Thus it is unsurprising -- nay, to be expected -- that these algorithms recommended greater quantities of distressing content \cite{Michel2024, Milli2025}. So if an algorithm keeps recommending unsavoury content to people, some (or indeed many) of them might just watch it regardless.

\subsubsection{Example: Dating app recommender algorithms -- Tinder}
 
Dating match-up apps also runs suggestion algorithms. Platforms such as Tinder use complicated, machine-learned algorithms which \textit{choose} who to suggest to a user as a potential match \cite{Rolle2019}. Again, such algorithms are far from random, and moreover take in to account much more than just  partner preferences the user submits (gender, age,  distance from them, etc.); in the case of Tinder, these algorithms factor in past swipe history, the way other people have swiped on the user, how the user interacts with other Tinder users, etc., and gives all users a ``rating'' in order to make matching suggestions with others on or near the same rating \cite{Harrison2023}. And thought it might not seem like it, such apps keep a vast amount of data on their users; one user requested to be send a complete copy of all the data Tinder had on them, and received 800 pages in response \cite{Duportail2017}, all of which is internally available to compute matching suggestions.

Platforms such as Tinder have extremely large userbases (60million users a month as of 2025 \cite{Curry2026}, with  actual impact far wider as people tend to stop using the platform when they are in a relationship) across many countries, and with most users (80\%) being of child-rearing age (18-44) \cite{Curry2026}. However, such vast and targeted reach, for platforms that are in the business of recommending romantic partners, could carry certain risks. By influencing the way so many people of child-rearing age meet romantic partners, these platforms might be altering the reproductive patterns of a substantial proportion of the human population, and thus (inadvertently) reshaping the human genome. This cuts both ways, with the recommender algorithms affecting who a user \textit{will} be shown (and thus selecting the pool of people the might choose to parter with), as well as who a user will \textit{not} be shown (thus excluding them from a large group of potential partners). Here, the recommender system ends up tailoring (and restricting) the set of actions the user can do, for perhaps one of the most important decisions of their life: who to parter and have children with.

Thus, the question now remains: what do these (predominantly technology) companies know and understand about relationship psychology, medical genetics, or any other disciplines relevant to romantic human relationships? And how do their choices and decisions about such matching algorithms impact how people reproduce?

\subsubsection{Example: Targeted (political) advertising algorithms -- Cambridge Analytica}\label{s:CA}

The idea of targeted advertising is nothing new; advertisers have always tried to be efficient in who they advertised their product to, in order to get the most return on investment (ROI). But usually this has been quite `coarse'; placing adverts for high-performance car tyres on car racing shows, or advertising fast food during children's TV shows (which in the US was subsequently regulated in 1990 with the \textit{Children's Television Act} \cite{story2004food}). But of course these adverts were targeted in a rudimentary, yet still effective, way: at `car enthusiasts', or `young children', and so on. Platforms such as Facebook have enabled targeted advertising at a whole new level. Now, advertisers can send an advert based on dozens of criteria, to the point where the advert may only be seen by as little as 20 people. This already causes a problem; in the past, a lot of people (including parents) saw fast food companies advertising to children, and could complain. Now, only the target audience sees the targeted adverts; no-one else. This removes society's ability to critique such adverts, and such small-group targeting can have serious and disturbing consequences for those being targeted. And such actions are simple, and cheap; one individual was able to send targeted adverts specifically to their housemate (and no-one else) over Facebook for a prolonged period, for a total cost of only \$1.70 \cite{Swichkow2014}.
 
Various entities have taken the idea of targeted advertising even further, and implemented targeted \textit{political} advertising on Facebook. Most seriously, this happened with the data science company Cambridge Analytica. They first had someone create a Facebook app called ``thisisyourdigitallife'' for which 270,000 users subscribed \cite{Facebook2018suspending}. Due to a peculiarity of the options available to app developers, the app was able to harvest not only the personal information of the 270,000 users who subscribed, but also the (approximately) 87million ``friends'' \cite{Isaak2018} of all these users \cite{Facebook2018DataAccess}. Cambridge Analytica then applied machine learning algorithms to this data set of 87million people and created a predictive tool that could profile people and deduce things about their voting preferences, life ambitions and preferences, etc., all through their Facebook profile. Cambridge Analytica then created many different targeted political adverts \cite{Isaak2018}, and directly advertised to all Facebook users in the US during the 2016 US presidential election \cite{LewisHilder2018}. This advertising was low-cost, accurate, and highly impactful. Thus, from a single, allegedly exploitative, app, and the lax data concerns of 270,000 people, the entire adult Facebook user base in the US of 167million people (68\% \cite{Greenwood2016} of the 245million US adults in 2016 \cite[Table 1]{USCensus2016}) was targeted; a vast proportion of voters. Facebook was eventually fined USD\$5billion by the US Federal Trade Commission for their role in the privacy breach \cite{Wong2019}, which also ruled that Cambridge Analytica had engaged in ``deceptive practices'' regarding data collection \cite{FTC2019CA}.

Being able to target political adverts with such granularity means they can be tuned to fit very specific demographics. But being able to specify who sees an advert makes it possible to specify who does not see it. Thus, a new opportunity to advertise in a selective, inconsistent, and self-contradictory way has arisen, with little possibility for scrutiny. Naively one would send adverts for politician X to people who show preference to politician Y, to change their mind. But advertisers can now use this in very creative, and arguably manipulative, ways. Hypothetically, one can send an advert saying ``Voting is important; make sure you vote'' only to people who might be inclined to vote for your party but not likely to vote (which voters may interpret as coming from some electoral agency), thus appearing partisan and appealing to their civic duty. Or sending ``You're party has a strong lead in the polls, and is guaranteed to win'' only to people inclined to vote for your opponent (which voters may interpret as coming from their preferred party), thus disincentivising them to vote.

Indeed, Cambridge Analytica sent out adverts in almost exactly this fashion, with reports that ``The Trump campaign seized the opportunity, showing two different ads to different categories of voters according to the detailed geographical information of visitors to the YouTube home page. Voters in areas where people were likely to be Trump supporters were shown a triumphant-looking image of the nominee, and help finding their nearest polling station. Those whose geographical information suggested they were not fervent Trump supporters, such as swing voters, were shown photos of his high-profile supporters, including his daughter Ivanka Trump, a celebrity from the reality TV show Duck Dynasty, and Dana White, the president of the Ultimate Fighting Championship'' \cite{LewisHilder2018}. The purpose of this was not to convince every voter how to vote, but rather to ``tip the scales'' of the election \cite{Isaak2018}. And so ``As described by Cambridge Analytica's CEO, the key was to identify those who might be enticed to vote for their client or be discouraged to vote for their opponent. Every vote added or disrupted (in the intended way) tips the election results'' \cite[p.~57]{Isaak2018}. And it does not take much voter movement to have a significant impact, as ``Tipping the scale in a few states—with as few as 100,000 voters—using individualized, high-impact messages is sufficient to impact election results. This might not be the only reason for the specific 2016 US election outcome, but there is every indication that it was a useful if not a critical contribution'' \cite[p.~58]{Isaak2018}.

\subsubsection{Overview of suggestion algorithms}

The above examples were presented to illustrate that the decisions that developers make can have critical (and potentially harmful) impact on people, even in the context of \textit{recommendations} which no-one is obliged to follow. Whether it be recommending more content (leading people into radicalisation vortices), or recommending romantic partners (leading to how people will reproduce), or recommending political adverts (leading people to vote in a certain way), we made it very clear that none of these domains are in any way `impartial'; by their very nature, they are there to influence people in a specific way, for a specific (behavioural, economic, or political) outcome. Mathematicians might be being used here. Used to exploit unassuming citizenry. Used to generate profit in harmful ways. Used to gain political power. And they might not even realise it.

\subsection{Rating algorithms}

We moved on to another type of algorithm: \textit{rating algorithms}. These take collections of things (usually humans), and give them a (usually linear) rating. Such algorithms might \textit{seem} to be impartial processes, but, as we demonstrated to students, that need not be the case.

\subsubsection{Example: Credit scoring algorithms}

Credit scoring has existed for many years, and comes in three basic types \cite[p.~5]{Upturn2014}. The first are \textit{traditional models}, which take in baseline credit data such as identity, address history, past and existing loans with their amounts and payment histories, and default or bankruptcy proceedings. Then there are also \textit{mainstream alternative models} which take in alternate mainstream data such as utility payments and other regular payments (mobile phone contracts, gym memberships,  etc.).  And finally (and more recently) there are \textit{fringe alternative models} which take in fringe alternative data such as shopping habits, social media profiles, web tracking and location data, etc. The first two types of models have been around for many decades; the last of these is relatively new.

There are several companies doing fringe alternative credit scoring, and they make heavy use of mathematicians in the process. One of them, ZestFinance, claimed on its website in 2018 \cite{ZestFinance2018} that ``the team has grown to more than a hundred crazy smart data geeks, mathematicians and computer scientists''. Another, Lenddo \cite{Lenddo2013}, sought  direct access to applicants' social media accounts such a Facebook, via login credentials (i.e., their Facebook username and password provided to Lenddo) \cite[Column 9, Line 45]{LenddoPatent2014}, to harvest and analyse their social media data in order to compute a ``credit score''. This includes who the applicants friends are, the places they visit, biometric data, etc., which are all explicitly specified in the Lenddo patent \cite[Column 9, Line 45]{LenddoPatent2014}.

The potential for bias here is significant, as there are so many pairs of datapoints that could be correlated (which the mathematics in the machine learning will `see') but have no causal link (which the mathematics, and indeed the ``crazy smart mathematicians'', have no mechanism for distinguishing). Credit scoring is a crucial part of many people's lives, from borrowing money for recreation (holidays), to essential items (cars), to crucial assets (housing). This practice of using social media data to determine credit worthiness has been  criticised \cite{Eichelberger2013}, as ``It’s an unprecedented practice that consumer advocates say can be unfair or discriminatory'' \cite{Eichelberger2013}. They specifically single out Lenndo which they say ``denies credit to applicants who are Facebook friends with someone who was late repaying a Lenddo loan'' \cite{Eichelberger2013}.

\subsubsection{Overview of rating algorithms}

With these, we were able to help students see that, just because a \textit{computer} has carried out an assessment using \textit{data} and a new \textit{algorithm} -- rather than, say, a \textit{person} using an \textit{application form} and a standard credit model -- does not necessarily make the process any better; more data does not mean better results. Indeed, given how little ``crazy smart data geeks, mathematicians and computer scientists'' might know about complicated social and financial issues that people face, there may be a greater propensity for bias to appear in this type of credit scoring.

\subsection{Machine-learned algorithms and AI}

We then moved on to cover machine-learned algorithms such as AI. This is a vast, and ever-expanding field; AI ethics is already a well-established area of study, and of concern. And moreover, the pace of AI development, and scale of deployment, is rapidly accelerating, especially in specific areas such as frontier AI models (e.g., Large Language Models) whose reach and capabilities are evolving almost daily. So as to not overwhelm the students with this particular sub-area of ethics in mathematics, and building on the fact that it has so much public attention, we only gave it a brief coverage.

We explained that creating machine-learned algorithms or AI requires designing the learning incentives and cost functions of the algorithmic training, and then choosing training data sets to input, and that both of these steps are a human construct. Students saw that, by the nature of machine learning, it is possible to take the same data set and produce two different algorithms which purport to achieve the same overall aim, simply by changing the incentives and cost functions. We pointed out that, when developers use phrases like ``the data says X'', it is somewhat nonsensical, as those developers could make the same data set `say' something different, or choose a different data set altogether. 
 
We gave an example illustrating just how important choice of data set is: the AI soap dispenser that would only work on white hands, and not work on black hands (complete with video) \cite{Gizmodo2017}. In that instance, there was a workaround: placing a white piece of towel on a black palm made the machine dispense soap, suggesting the AI had learned to dispense to a white, hand-sized `blob', rather than an actual hand. But already this demonstrated unintended bias: by not making sure to choose a set of test data that sampled from the whole population -- in this case, they probably sampled mostly, or even only, white hands -- the resulting AI's operation was biased against black hands, as it simply lacked the training data to recognise them. We posed to students how convincing it might be for a mathematician developing such software to think ``I have 1million images of hands; \textit{surely} this machine learning is going to be extremely accurate!'', and they saw that just increasing the \textit{quantity} of training data does not necessarily increase its \textit{diversity}, and thus \textit{overall quality}.

We brought up the idea of `scrutinising the AI code, perhaps by making it open-source' as a mechanism by which to avert or resolve such issues, and highlighted that in general scrutinising code is very difficult. An AI algorithm is often just a machine-learned parameter matrix, and in the case of modern (i.e., 2026) LLMs this can be of the order of 1trillion entries (1million x 1million matrix) \cite{Rohilla2026}. Unfortunately, automated evaluation cannot work \textit{in general}, as it is algorithmically impossible to \textit{always} determine if a piece of computer code is just the function $f(n)=n+1$ or not\footnote{Indeed, for \textit{any} Turing machine $T$ which is \textit{total} (i.e., halts on all inputs),  the set \textbf{Tot} (of Turing machines which are total) many-one reduces to the set of  Turing machines that compute the same function as $T$ (call this set $S$). Since \textbf{Tot} is $\Pi_{2}$-complete in Kleene's arithmetic hierarchy, then $S$ is at least $\Pi_{2}$-hard, so is not computable (and, more strongly, neither $S$ nor its complement can be computably enumerated). As $f(n)=n+1$ is a total computable function, this result applies to $f$. See \cite[Theorem 4.3.2]{soare2016turing}.}. And even if a piece of code is comprehensible by humans, having it as open source does not always mean problems will be identified rapidly, as in the case of \textit{Heartbleed} which was one of the most catastrophic bugs ever to appear in open source software (OpenSSL); it  went undetected for 2 years \cite[p.~476]{Durumeric2014}. We highlighted that even very well resourced organisations struggle with this; Google struggled to understand their own machine-learned ad-placement algorithms, as demonstrated when they placed official UK government adverts next to extremist content \cite{Independent2017}, thereby causing extremists to receive derived ad revenue from the advertiser (the UK taxpayer). Of course, the ad placement algorithms were doing precisely as they were trained, but with undesirable consequences that Google did not detect.

\subsection{When algorithms meet algorithms -- Algorithmic collisions}\label{s:algo_collisions}

We then introduced to students a new, and under-reported,  danger of mathematics: when algorithms meet algorithms, also known as \textit{algorithmic collisions} \cite{chiodo2025problem}. This is well-illustrated by the (somewhat silly) example from 2011, when two algorithms from different Amazon sellers interacted to raise the price of a new book (which was selling second-hand for USD\$35) to \$23million \cite{Eisen2011}. There, the process was simple: seller A would always try and be the cheapest by setting their price to 0.9983x the next cheapest price. Seller B would always set their price to 1.270589x the cheapest price (presumably they did not have the book, but had such a good seller rating that they could buy it from the cheapest seller and onsell it for a profit \cite{Eisen2011}). Sellers A and B met on the Amazon marketplace, and were the only sellers there. Thus their two algorithms `met', and after many iterations the price skyrocketed to \$23million. Each seller seemed to have a `clever' algorithm to out-wit other sellers and potential buyers, but neither considered the possibility that the other algorithm might exist.

We then posed to students what might  happen in a more serious instance, such as self-driving cars? At the moment, self-driving cars are trained on real-road data sets of cars. But these datasets only have \textit{human} drivers on them. Thus, self-driving cars know how to deal with human drivers. We  explored what might  happen if a self driving car (say from manufacturer X) met another self driving car (from manufacturer Y)? There is simply no training data for driving amongst self driving cars, so neither car is very well prepared for such an encounter. If the Amazon book example from before reappears in a much worse fashion, then (for some unexpected reason) both cars might `decide' to avoid a head-on collision, with one veering left, and the other veering right, still leading to collision. This illustrated a very deep problem, whereby such `algorithmic' collisions might occur in ways that cannot be initially detected or measured by those designing each algorithm, and that there is a significant responsibility on developers of algorithms to put in place mitigations to detect, dampen, and ultimately prevent such disasters \cite{chiodo2025problem}.

\subsection{Algorithmic fairness}

We highlighted that many of the cases and algorithms in this lecture (in particular, COMPAS) have generated much interest among researchers about \textit{algorithmic fairness}. That is, understanding when an algorithm is, or is not, treating people fairly. We noted that wider society is still debating the notion of `fair', and even philosophers (who have been thinking about this for thousands of years) are yet to come up with a succinct, or even mutually-agreed, notion of what it means for something to be `fair' or `just' \cite{Sandel2009}. 
 
We pointed out that computer scientists and mathematicians have started to take it upon \textit{themselves} to single-handedly define, mathematically, what is meant by `fairness' of an algorithm. Such discussions and mathematical definitions go back to the 2011 work of Dwork et al.~\cite{Dwork2012}, with a later paper by Friedler et al.~\cite[Definition 3.1]{Friedleretal2016} summarising and extending this definition as:

 \begin{quote}
 A mapping $f : CS \to DS$ between metric spaces $CS=(P, d_{P})$ (Construct Space) and $DS=(O, D_{O})$ (Decision Space) is said to be \textit{fair} if objects that are close in CS are also close in DS. Specifically, fix two thresholds $\epsilon, \epsilon'$. Then $f$ is defined as $(\epsilon, \epsilon')$-fair if for any $x, y \in  P$, we have $d_{P}(x, y) \leq  \epsilon \Rightarrow d_{O}( f(x), f(y)) \leq \epsilon'$. Note that the definition of fairness does not require any particular outcome for entities that are far apart
in CS.
\end{quote}
Students initially thought this was a great advancement in understanding, until we revealed that all is not what it might seem. For example, in a later (2018) paper \cite{dwork2019}, Dwork et al.~then show that this definition of fairness is not closed under composition. That is, applying a fair algorithm to a fair algorithm need not yield a fair algorithm; a somewhat counter-intuitive property of how `fairness' might (intuitively) work. Moreover, the paper by Friedler et al.~\cite{Friedleretal2016}, as well as some other papers \cite{Miconi2017}, then go on to study various impossibility results for fairness, akin to Arrow's impossibility theorem for voting \cite{arrow1950difficulty}. By taking different `intuitive' mathematical definitions of fairness (known as \textit{demographic parity} -- the definition by Dwork et al., \textit{equalized odds}, and \textit{predictive rate parity} \cite{Hsu2022}), one can show that they are mutually inconsistent, in the sense that there is no algorithm that is simultaneously satisfies all these `fairness' conditions. 

We used this discourse to demonstrate to students that, while such research is academically interesting, and probably done with the best of intentions, it is not necessarily the case that problems brought about by mathematical processes (e.g., COMPAS) have clear-cut mathematical solutions, and that, given how complicated a social concept such as fairness is, it might be beyond the abilities of technical researchers to `solve it' with purely technical tools. Indeed, when Dwork was asked ``Whose responsibility is it to ensure that algorithms or software are not discriminatory?'', she responded with ``This is better answered by an ethicist. I'm interested in how theoretical computer science and other disciplines can contribute to an understanding of what might be viable options. The goal of my work is to put fairness on a firm mathematical foundation, but even I have just begun to scratch the surface. This entails finding a mathematically rigorous definition of fairness and developing computational methods - algorithms - that guarantee fairness'' \cite{nyt_dwork_2015}. We suggested that such a response defers responsibility and judgement to ethicists, but that such people are  unlikely to comprehend the mathematical definitions researchers such as Dwork are producing, once again drawing on the notion from Section \ref{s:blackscholes} that  ``If we don't understand the derivation, then we can't know the limitations''.

To illustrate how such work might lead to tangible harm, we posed the following hypothetical scenario to students: Suppose company X creates an algorithm (say, for assessing job applicants) and uses it to filter applicants. Then suppose a rejected applicant takes them to court for unfair discrimination in the selection process. That company might then present to a judge the argument ``We developed our algorithm, and we have mathematically \textit{proven} it is fair'', making reference to the algorithmic fairness papers above. And assume that the algorithm is indeed fair by (one of, but of course not all three of) the main three mathematical definitions of fairness. How, then, is a judge to consider this evidence? Students saw that an algorithm (which already possesses all the authority of mathematics), being \textit{verified} by yet another mathematical (published, peer-reviewed by experts) process, gains a `Teflon coating', and is almost impossible to challenge.

\subsection{Accountability when creating algorithms}\label{s:insider}

The lecture so far challenged the notion that algorithmic processes are tools of unmitigated good, and we took this one step further by exploring whether mathematicians should be allowed to develop them without restriction, or whether there are limits that should be, or are, in place. We illustrated this with a specific, direct question to students: ``Are there any data sets we should not, or are not allowed to, train machine-learned algorithms on?'' Students were initially perplexed by the question, asking why this would be the case. So we presented existing examples in society which might serve to guide mathematicians, such as rules and laws surrounding insider trading \cite{bainbridge1998insider}. In many jurisdictions, it is illegal to go in to a company, `learn' about that company, and then walk out and (without taking any data) act on what was learned by buying/selling shares in that company. That is, one can't train \textit{oneself} to act on confidential company data. And such principles now have analogues in machine learning, with frontier AI companies finding themselves in all sorts of complex legal (and moral) quandaries due to their use of vast amounts of copyrighted data (text, images, videos, etc.) scraped from the internet and used in LLM training.
 
We then explored the (mis)conception about algorithmic processes such as AI that humans will always be able to ``maintain control'' \cite{chin2026reframing}, showing that is no longer the case. Common algorithmic safety mechanisms such as having a `human-in-the-loop' are now being show to suffer from serious shortcomings and numerous failure modes, to the point where they cannot be considered a universal infallible safety method \cite[p.~5]{chiodo2025formalising}. We revealed that there are now ways to implement (machine-learned) algorithms that cannot be `turned off', as `unstoppable transactions' are being deployed on the Etherium blockchain in the form of DAOs \cite{coindesk_dao_2019}. Humanity cannot turn these off, as the blockchain is a completely decentralised system and no one entity possesses the computational power to take it over, or the physical means to shut it down. And finally, we showed that super-intelligent AI which can evade human oversight is now being studied as a very serious and plausible threat to society, and plausible mechanisms for the `Skynet scenario' of AI taking over, such as by being handed the keys to the city to resolve a (different) crisis, are now being examined \cite{chiodo2026handing}.

And thus, we asked students: ``When a mathematical process makes a moral decision, who do we hold morally accountable?'' Suppose a self driving car hits and kills someone; who is at fault? Suppose COMPAS gives someone a disproportionately high risk score leading to an unnecessarily long jail term; who is at fault? Suppose someone is denied a mortgage because an alternate credit scoring system gives them an undeservedly-low credit score; who is at fault? Suppose someone is repeatedly recommended matches on Tinder which become abusive partners; who is at fault? We left the students with the following open question to consider: If mathematicians can claim that the algorithmic processes they create are fair and impartial, and these are used to make  decisions in a way that causes harm, how much accountability sits on them as developers? And furthermore, if a mathematical process makes a moral decision, who should be held morally accountable? All these issues were discussed later in Lecture \ref{s:law}.

\subsection{The myth of impartiality}\label{s:impartiality-myth}

By now, the students had seen the purpose of the lecture: to dispel the myth of impartiality. No mathematical process is completely impartial, because those who develop it are not. Developers' own internal bias -- coming from their own experiences, lessons and cultural influences -- impacts how they design and create systems, and often leads to the creation of bad systems (usually unintentionally, but bad nonetheless).  Into every algorithm  or model they create, they pour all their prejudices, their ignorance, their cultural norms, their moral beliefs, and their subconscious aspirations -- where they want to go, and where they want the world to go). Their algorithms are loyal, but only to them. They create tools of great power, but how that power is wielded is influenced by its creator. 

We suggested to students that often, too much `honour', or `authority', or `impartiality' might be assigned to mathematics, and mathematicians might defer their own personal moral judgement (which they may claim can be flawed) to the `purity' of mathematics (discussed later in Lecture \ref{s:community}), which they may claim is impartial. But mathematics cannot be completely impartial, because its creator never is.  Yet this doesn't seem to be widely-understood by mathematicians, and one can perhaps see evidence of this in \cite{siam_terrorism_2016}, where the author writes ``Modelling data is an economically cheap option that lets the data speak for itself, extricating the researchers' emotions from the equation''. We put a different view to students: namely, that modelling data might calcify the researchers' beliefs, emotions, desires and opinions into the equation, in a way that often cannot be seen, cannot be scrutinised, cannot be understood, and cannot be undone.

We finished the lecture by asking the students a simple question: ``Do you think that the algorithm you described at the start of this lecture was impartial?'' None have ever said yes.

\begingroup
\let\section\subsection
\renewcommand{\bibname}{Bibliography (Lecture 4)}
\renewcommand{\refname}{Bibliography (Lecture 4)}
\addcontentsline{toc}{subsection}{Bibliography (Lecture 4)}
\renewcommand{\addcontentsline}[3]{}
\bibliographystyle{IEEEtran}
\bibliography{bibLec4}
\endgroup

\newpage
\section{Lecture 5: Regulation, accountability, and the law}\label{s:law}

\fbox{
  \parbox{0.98\textwidth}{
Mathematical work in industry is now very close to its tangible impact, and such outputs are now extremely quick and easy to implement, with impacts that are amplified by modern computational processes. But with such great power comes great responsibility, and so various aspects of criminal, civil, and moral accountability can quickly come in to play for those carrying out such work.

This lecture looks at how mathematicians can learn just how close they now sit to the applications of their work, and associated accountability mechanisms they may face even in circumstances where they were `just following orders'. It shows how there are things mathematicians are, and aren't, legally (or morally) allowed to do, and how they can change course to avoid consequences in the courts of justice and of public opinion for falling foul of such expectations.
}}

\subsection{Revisiting past case studies -- moving from examples to theory}

This was the first lecture where we transitioned from the examples-based `part 2' content (showing \textit{there exists ethics in mathematics}), to more theory-based `part 3' content (showing that \textit{for all mathematics, there are ethical issues}), reminding students that universally-quantified (i.e., `for all') statements cannot be proven \textit{by example}, but rather, require the development of \textit{theory}. We began by revisiting some past case studies where mathematicians went to jail, such as those involved in the collapse of FTX and Alameda Research (Section \ref{s:ftx}). But rather than viewing this as an example of \textit{what} happened, we started to look at \textit{why} it happened, and in particular, the existence of laws and regulations that everyone in society (including mathematicians) are expected to follow.  In this case, those FTX executives moved billions of dollars of depositor funds from FTX to Alameda, and then proceeded to lose vast sums on risky investments. Those jailed from FTX or Alameda were convicted of actions such as fraud and conspiracy, reflecting the fact that just because depositors give money to a company to hold, does not give those within the company full license to do whatever they want with it.

We also revisited several of our earlier `data releases' examples, whereby organisations simply chose to release data `because they could'. This included the Strava heatmap data release from Section \ref{s:strava}, and the various security issues that it caused, at both the individual, and state, level. It also included the Public Transport Victoria (PTV) hackathon travel data release from Section \ref{s:ptv}, and the significant individual security issues it raised, along with the stiff legal consequences for PTV. This illustrated the idea that, just because individuals give data to a corporation, does not mean that corporation has the right to circulate that data in whatever way they want. Indeed, there are now very strict data rules in various jurisdictions, such as in the EU with the General Data Protection Regulation (GDPR) which came into effect in 2016 and (very strictly) governs the use of (personal) information within the EU \cite{EU2016GDPR}.

And finally, we revisited the notions from Section \ref{s:insider} of insider trading, and  from Section \ref{s:CA} of the use of data in AI training for Cambridge Analytica. Both of these illustrate the same principle: namely, that there are laws in place which govern what sort of algorithmic training one can do with datasets. Such laws exist to prevent those who possess and train on the data from exploiting wider society. In the case of insider trading, those with the data can make stock transactions that are extremely favourable to them (i.e., at the direct expense of the other parties in such transactions). And with AI training on improperly-obtained data, a similar thing occurs, whereby those training on the data are able to make gains from it in ways that exploit wider society (i.e., by influencing elections). And we added to this the case of the company \textit{Stability AI}, which was recently found to be committing limited trade mark infringement in the UK High Court because its image generation software (Stable Diffusion) was inadvertently generating images with the Getty Images watermark, because they (Stability AI) had trained on so many images from Getty Images, all with that watermark \cite{pfitzer2025getty, smith2025getty}.\footnote{The lawsuit by Getty Images also alleged that Stability AI had engaged in a massive copyright violation of their images. But this part of the lawsuit was ultimately unsuccessful, as the courts found that the AI models were not an ``infringing copy'' of the images \textit{per se}, in the context of copyright law \cite{pfitzer2025getty, smith2025getty}. However, this case has been appealed to the Court of Appeal by Getty Images, so another assessment on the application of copyright law to this case is yet to come. Nonetheless, this has been a significant legal battle for Stability AI to fight; a smaller organisation may well have  capitulated given the vast effort and legal resources and costs needed to defend against a company as large as Getty Images. Reputationally, Stability AI have been the face of widespread media backlash resulting from the high-profile litigation.}

And so we explained our first (of many) \textit{universal} observations to students: just because a mathematician has access to something (e.g., other people's money, data, or trademarked logos), does not mean they can do whatever they want with them (e.g., invest them, circulate them, train AI on them, or reproduce them). We explained that the law can, and will, always apply here, even to mathematicians doing mathematical work. Laws – especially in \textit{common law} systems\footnote{\textit{Common law} is ``law developed through judicial decisions rather than enacted statutes'' \cite{wex_common_law}.} – do not explicitly refer to `AI' or `mathematics', but apply by adapting general principles and previously established reasoning to new fact patterns, e.g., existing rules on privacy, discrimination and fairness often still hold.

\subsection{What is the purpose of regulation and law?}

We picked up again on a previous example; this time of self-driving cars. In Section \ref{s:algo_collisions} we discussed \textit{how} self-driving cars might be involved in a collision, and potentially hit at kill someone. We now moved the discourse on to \textit{who} should be legally responsible if a self-driving car hits and kills someone? To do this, we explored some of the reasons for why society imposes laws. Criminal law (loosely) reduces to the following key ideas \cite[p.~6]{Wells_Quick_2010}:
  
 \begin{enumerate}
  \item \textit{Retributive justice}. Society wants to exert some sort of retribution on those who have done wrong.
  \item \textit{Deterrence}. Society wants to disincentivise certain actions or behaviour.
 \end{enumerate}
That is, society want someone to blame, to make the victims feel equality has been restored (justice), and wants to make an example of someone so that others won't repeat the action (deterrence).

We looked at how society might reconcile these principles with the case of a self-driving car hitting and killing someone. In particular, how might society visit retribution on, or disincentivise, a \textit{mathematical} process? Students saw that throwing a computer in jail might not make the victim or their family feel any better, nor would it likely change the behaviour of the next self-driving car. And so, it became clear that, in a functioning legal framework, some humans would need to suffer legal consequences, if the law is to have its desired effect.

\subsection{Self-driving cars and their legal ramifications}

\subsubsection{Example: The first pedestrian fatality from a self-driving car}

We explained that, the first few times this idea of a self-driving car fatality was discussed in the lectures (2016-17), it was a hypothetical question. Unfortunately, in 2018, it became a real example, when Elaine Herzberg was the first pedestrian hit and killed by a self-driving car. For ease of reference, the following two paragraphs are taken almost exactly from \cite{chiodo2025formalising}; the bracketed mentions of how the setup failed refer to some of the common human-in-the-loop failure modes given in \cite[Section 3.1]{chiodo2025formalising}:

\begin{quote}
``In March 2018 the first recorded incident of a pedestrian fatality from a collision with a self-driving car occurred; an Uber self-driving test vehicle in the US hit and killed Elaine Herzberg. The car had a human-in-the-loop setup created by Uber, with a human driver (an Uber employee) at the wheel poised to intervene and `take over' the driving at any point if they observed a problem with the self-driving AI system. The human operator monitored car and had the ability to intervene by braking and/or taking the wheel, but they did not have the ability to change the autonomous decision-making in more complex ways. The human monitor was expected to perpetually `hover their hands above the steering wheel and foot above the brake pedal' (somewhat unrealistic expectations), without ever touching them as that would disable the self-driving mode. The car had been in self-driving mode for 19 minutes before the crash, with no tasks required of the human driver (inducing fatigue and automation bias). The self-driving AI first identified the pedestrian 5.6 seconds before collision, as they were jaywalking across the road. However, as it had not been programmed to classify a jaywalker, it spent the next 4.1 seconds misclassifying the pedestrian as `another vehicle', `bicycle', and `other' (the AI had unexpected inputs). The AI first predicted a collision 1.5 seconds before impact and began evasive action; at 0.2 seconds before impact the AI established that impact was inevitable, instigated a controlled slowdown, and gave the human driver the first warning notification of upcoming impact (this was a delayed notification, giving the driver insufficient reaction time). The human driver had, from 6 seconds before impact, been gazing down at the control panel, allegedly watching their cellphone (clearly incongruous intentions). With a warning lead-in time of 0.2 seconds---approximately human reaction time---the human driver had virtually no time to take control; they took the wheel 0.02 seconds before impact, and the car then hit the pedestrian. The entire human-in-the-loop setup was affected by a poor safety culture (with insufficient support, and other external pressures). In this example, the failure cascade spans across the human-in-the-loop failure taxonomy numbering in \cite[Section 3]{chiodo2025formalising}: (1) the AI components failed by not recognising the pedestrian, (2) the workflow failed as the human allegedly was not embedded in a supporting safety-culture and had insufficient reaction time and unrealistic expectations, (3) the human--machine interface setup failed by allowing an environment whereby the human could watch a video, and (4) the human failed by being inattentive. However, one key further failure of this setup arose post-incident, when the law failed to hold the company accountable in any way or form, corresponding to (5).''  \cite[Section 3.3]{chiodo2025formalising}.

``Here, we see that both criminal and civil liability fell on the Uber \textit{employee}; the human operating the car. This was despite allegations that Uber's technology was flawed (i.e., not recognising the jaywalking pedestrian) and the company's poor safety culture (for case details, see~\cite{stamp2024reckless}).'' 
\cite[Section 4.1]{chiodo2025formalising}.
\end{quote}

\subsubsection{Analysis of self-driving car fatality example}

After giving the self-driving car example, we took some time to discuss what actually went wrong in the process, and who (or what) might have `made a mistake' worthy of punishment. Further details of these shortcomings, and what could have been done differently or better, can be found in \cite[Appendix C.2]{chiodo2025formalising}.

We covered the fact that the self-driving AI had learned a very strict definition of ``pedestrian''; it could only recognise a pedestrian on the road if they were crossing at a pedestrian crossing, not a jaywalker, hence why Elaine was categorised as an ``unclassified object'' on the road. But as an `AI mistake', and in the knowledge that one cannot throw a computer in jail, who was ultimately responsible for it?

We then explored the safety mechanism of a human-in-the-loop (the backup driver), and whether this provided meaningful safety. While such solutions might sound ideal to a mathematician or other tech developer, it was clear in this instance that the normal failure modes of humans in such setups was not properly considered. For example, can a human be reasonably expected to hold their hands and feet over the steering wheel and brake of a car for prolonged periods? Or to react to a situation with 0.2 seconds notice? And so, if theses expectations fail, who is `at fault'; the human driver (who was expected to be \textit{superhuman} \cite[Section 4.2]{chiodo2025formalising}), or those who designed and developed the setup (those who created these expectations, and put a human in place to carry them out)?

Turning to the fact that the human driver was able to watch videos on their phone, we suggested to students that a system that allowed the driver to bring a device capable of playing videos into the car was itself a safety oversight, and therefore point of fault and liability for whoever designed the system. And, more generally, the allegedly ``poor safety culture'' at Uber could have contributed to the litany of failures that led to the crash, leading to Uber shouldering some (corporate) liability,

But ultimately, in this case, it was the driver who was found \textit{fully} at fault by the courts \cite{stamp2024reckless}, not Uber, which was completely exonerated. So we left the students with two questions to consider. First, we posed whether the law here succeeded in achieving both retributive justice and deterrence. That is, whether all parties that contributed to the harm caused were `punished' in some way, and whether the punishments given out would succeed in deterring future parties in such scenarios from carrying out the same actions. Second, we posed whether, overall, they thought that the law had done a good job of interpreting the situation and acting in society's best interests, or whether there were any shortcomings they could see in how the law treated what had happened? This naturally led on to the next activity, on what the law \textit{should} do.

\subsubsection{Activity: Self-driving car optimisation process}

We highlighted to students that self-driving cars are not solely produced by the usual people who build human-driven cars; they involve several much more technical processes, in particular, machine-learning. Thus, mathematicians and computer scientists are now at the forefront of car manufacturing. So we made the following simple activity for students, asking them:  ``In a potential collision situation, what should a self driving car optimise over?''  They saw this as a reasonable question to consider, given that they (as mathematicians) may one day need to implement such an optimisation strategy. We presented the following list of options, but also solicited further comments and thoughts:

\begin{enumerate}
\item Total number of lives.
\item Lives of the occupants of the car.
\item Total insurance payout.
\item The net benefit/harm to the company that built or sold the car.
\item Something else.
\end{enumerate}
Students began to see the utilitarian perspective of each of the listed options, and realised that each optimisation was in the interest of some particular group or organisation. We followed up by asking if the answer to this was obvious, and some students said it was, but with different `obvious' answers. We pointed out that this was similar to the trolley problem introduced in Section \ref{s:disease}, but much more `real' (i.e., will happen much more often), and requiring much faster decision-making (within about 1 second). 

We then modified the question slightly, posing \textit{who should decide} what to optimise over? We queried whether this was something that mathematicians should necessarily decide, as they'll be the ones implementing it in practice? We also put it to students that a (team of) ethicists or lawmakers may decide on one sort of optimisation, but then the mathematicians need to figure out a way to realise that optimisation. For example, if it is decided to optimise (by minimising) the total insurance payout, how might that be calculated, or even estimated? What factors should the mathematicians consider, and do they even have access to the information they seek? Moving the conversation towards regulation and accountability again, if for some reason this implementation is not effective, and is in some way `wrong' (or indeed, if the decision on what to optimise over is later deemed `wrong'), who \textit{should} (and who \textit{is}) legally \textit{accountable}? Would this be the ethicists? The mathematicians? The car company? The owner of the car? Or other individuals or groups entirely? This is different to who is (morally) \textit{responsible} to prevent any such harm from occurring; the person whose `job it is' to stop something going wrong is not always the one who is accountable if it does go wrong -- responsibility and accountability are different things.

We left the students thinking about what the current laws might (or should) be surrounding this, and how might regulators and lawmakers begin to go about considering these very new technologies and the associated potential for harm, especially seeing as such laws aren't always drafted by those who are trained in the mathematical sciences. Furthermore, we highlighted that existing cases may involve different fact patterns, for example, a case about causation of cancerous diseases from asbestos fibres in workplaces may inform how the courts respond to similar responsibility gaps in the case of decisions made by machine learning models \cite[Section 4]{chiodo2025formalising}.

\subsection{What sorts of accountability are there?}\label{s:sorts}
  
This section covered the idea that accountability in society is not simply a binary state; there are different ways in which we can do wrong, and thus different ways in which people can be held accountable. Even in the strict eyes of the law, there are different forms of accountability, which have different judicial processes. To summarise, we covered three main types of accountability:
   
\begin{enumerate}
    \item \textbf{Criminal accountability} (i.e., who is going to jail?). This is where society has laid down some very clear rules (criminal laws) which, if violated, lead to criminal prosecution (usually by the state), and then some form of state-enforced punishment (often incarceration).
    \item \textbf{Civil accountability} (i.e., who is going bankrupt?). This is where individuals or companies in society have a private dispute with others, and appeal to a state-run judicial process to seek compensation or redress. This form of accountability (usually) cannot lead to incarceration, but can lead to (vast) sums of money changing hands, and corrective orders.
   \item \textbf{Moral accountability} (i.e., who is going to hell?) This form of accountability is not written in any form of law, though often underpins the legal principles which make up criminal and civil law. But it can still have a significant effect on people, especially in terms of reputation. 
\end{enumerate}
Here, moral accountability is not part of the legal process, but still has `law-type' effects. We explained that, for example, if a large car manufacturer make cars with a propensity to skid off the road, then even if they can absolve themselves in the courts (perhaps with many high-priced lawyers), their public image is tarnished, and they may find it difficult to keep selling cars. Or if someone is a serial adulterer, and this is publicised over social media, getting a date in the future might be tough. But more seriously than that, moral accountability is humanity's own understanding of right and wrong; if someone does something bad -- such as contribute to corporate activities that make profit but cause significant harm -- then that is an action they need to live with for the rest of their life. And it may have seemed like a `good idea' at the time, but if their perception of life changes in the future, that's something that they need to face every morning when they look in the mirror. 

We tied this back the Uber self-driving car case, by putting to the students that if they were part of the coding team who developed the algorithm, and implemented the safety setup, that led to the death of Elaine Herzberg,  even if they and their company were absolved in criminal and civil courts, how would they \textit{feel}? That is moral accountability in it's very deepest form; one can never escape one's own self-judgement.

\subsection{Following orders -- the Nuremberg Defence}\label{s:Nuremberg} 
 
A common, though not necessarily valid,  legal defence people often try and use when implicated in some wrongdoing as part of a command structure (workplace hierarchy, military command, etc.) is the argument  `I was just doing my job' or `I was only following orders'; more commonly referred to as \textit{Superior Orders} \cite{dinstein2012defence}.  That is, if someone's `boss' tells them to do something, does that absolve them of all accountability? And moreover, what sort of accountability does this refer to here:  Moral accountability? Criminal accountability? Civil accountability? All of these? None of these?

The defence `I was only following orders' is often used during war crime trials, where it is usually referred to as the \textit{Nuremberg Defence} in light of its use at the Nuremberg trials following WWII \cite{dinstein2012defence}. However, it has has a long history, with inconsistent and mixed levels of success. In the 15th century, Peter Von Hagenbach was found guilty of ``crimes against the laws of nature and God'' \cite[p.~7]{bassiouni2006miller} (i.e., sack and pillage) and executed, even though he ``acted on the orders of the Duke of Burgundy'' \cite[p.~7]{bassiouni2006miller}.  At the Nuremberg trials themselves, it was established that the Superior Orders defence could only be used to \textit{lessen} a punishment, but could not absolve someone of illegal acts simply based on the fact that they were following orders \cite{dinstein2012defence}.

The 1998 Rome Statute of the International Criminal Court is a rather modern interpretation of this defence \cite{RomeStatute1998}. It states in \cite[Article 33]{RomeStatute1998} that:

\begin{quote}
\begin{enumerate}
    \item The fact that a crime within the jurisdiction of the Court has been committed by a person pursuant to an order of a Government or of a superior, whether military or civilian, shall not relieve that person of criminal responsibility unless:
   \begin{enumerate}
        \item The person was under a legal obligation to obey orders of the Government or the superior in question;
        \item The person did not know that the order was unlawful; and
        \item The order was not manifestly unlawful.
    \end{enumerate}
\item For the purposes of this article, orders to commit genocide or crimes against humanity are manifestly unlawful.
\end{enumerate}
\end{quote}
And so we finished by querying whether the activities of a mathematician would be absolved by the Rome Statute, or any other principles in law? We asked students to carefully consider, if their boss or other superior tells them to do mathematical work, whether they may be \textit{personally} accountable, regardless of the outcomes, and whether they can always make effective use of the Nuremberg defence?

\subsection{Designing an algorithm}\label{s:designinig}
 
To make the discourse in this lecture as relevant as possible, we turned our attention to  the act of designing and/or coding up an algorithm, and considered the question ``What sort of accountability is there for those who merely design or code an algorithm, without partaking in it implementation or use?'' Is there any legal or moral accountability for such actions, and are there any examples of this? The initial responses we received from students were usually along the lines of ``Only those \textit{using} the algorithm could be accountable'', thus implicitly absolving mathematicians and other technical developers only involved in the design stage. We then proceeded to give the following examples refuting this idea.

\subsubsection{Example: CIA torture algorithms}

In  July 2002, less than a year after the 9/11 terrorist attacks, two psychologists -- Dr James Mitchell and Dr John Jessen -- were approached by the CIA to develop ``enhanced interrogation techniques'' to used on terror suspects. They founded a company,  Mitchell, Jessen \& Associates, and signed a contract with the CIA to the value of USD\$180million, at least \$81million of which eventually paid out to them by 2009 \cite{Windrem2014CIATortureTeachers}. They soon developed a ``menu of 20 enhanced techniques'', which included ``waterboarding, sleep deprivation and forcing prisoners to assume ``stress positions'', being referred to as ``Torture teachers'' \cite{Windrem2014CIATortureTeachers} due to the nature of these techniques and the fact that some are ``widely viewed as torture'' \cite{Fink2017CIATorture}.

To use clearer language: Mitchell and Jessen designed various ways to torture inmates, and then provided these designs to the CIA. Indeed, the CIA themselves later referred to Mitchell and Jessen as the ``architects of the interrogation program'' \cite{Fink2017CIATorture}; they had, in essence, designed the `algorithms' for the interrogation program. The New York Times reported that two of these inmates tortured via such techniques -- Mohamed Ben Soud, and Suleiman Salim -- were later released; another -- Gul Rahman -- died in custody \cite{Fink2017CIATorture}. In 2015, the American Civil Liberties Union filed a lawsuit against Mitchell and Jessen, on behalf of these three former inmates (or their estates) for damages as a result of the alleged torture \cite{Apuzzo2015PsychologistsSuedCIA}.  It is important to note that  ``Mitchell and Jessen didn't directly interrogate the 2 surviving men who brought this case, but they methods they proposed were used on them, and at least 37 other men'' \cite[Video1 from 2:00]{Fink2017CIATorture}. 

In July 2017, the lawsuit was permitted to proceed to trial, which the judge said would focus on whether the plaintiffs ``were subjected to torture, and if so, did the defendants aid and abet — legally aid and abet and factually aid and abet — that torture?'' \cite{Fink2017Psychologists}. The following month, 
the plaintiffs settled  for an undisclosed sum \cite{Chappell2017Settlement}; as civil proceedings, this would only have ever involved monetary payment, rather than jail for them. Fortunately for Mitchell and Jessen, ``the CIA provided a multi-year indemnification agreement to protect the company and its employees from legal liability arising out of the enhanced interrogation program'' \cite{Windrem2014CIATortureTeachers}, and so the CIA would have paid the settlement. Had they not been indemnified, Mitchell and Jessen would have been \textit{personally} liable for any settlement amount.

The New York Time described this situation well by saying  ``This entire case really gets to the question of personal responsibility within a larger system''. That is, an instance (and legal test) of Superior Orders. It is worth noting that, had Mitchell and Jessen tried to use the Nuremberg defence, they may have run in to difficulty; they were alleged to have been involved in or abetted torture, and  according to Article 7 Item 1(f) of the Rome Statute, ``Torture'' is a crime against humanity \cite[Article 7]{RomeStatute1998}. This would have likely prevented any such use of  a Nuremberg defence, as described  in Section \ref{s:Nuremberg}.

\subsubsection{Example: The EU AI Act}

The EU AI Act \cite{EU2024AIAct} came in to force in the EU on 1 August 2024, setting out ``clear requirements and obligations regarding specific uses of AI'' \cite{EuropeanCommission2024AIActEnters}. With 113 `Articles' (i.e., broad items) and 13 Annexes \cite{EU2024AIAct}, this is a substantial piece of legislation, spanning hundreds of pages, and covering almost all aspects of AI development and use in the EU. Adherence to the AI Act is not optional; it is law, across all EU member states.

The AI Act makes a very clear  distinction early on in items (3) and (4) of Article 3, where it defines \textit{providers}, and \textit{deployers}, of AI systems \cite[Article 3]{EU2024AIAct}:

\begin{quote}
\begin{enumerate}
    \item[(3)] `provider' means a natural or legal person, public authority, agency or other body that develops an AI system or a general-purpose AI model or that has an AI system or a general-purpose AI model developed and places it on the market or puts the AI system into service under its own name or trademark, whether for payment or free of charge.
    
    \item[(4)] `deployer' means a natural or legal person, public authority, agency or other body using an AI system under its authority except where the AI system is used in the course of a personal non-professional activity.
\end{enumerate}
\end{quote}
Loosely speaking, AI providers are those who \textit{make} AI systems, and AI deployers are those who \textit{use} AI systems. An entity (e.g., a company) may be acting as both a provider and deployer if it makes and uses its own AI. But it can also act as just one of the two. Hence, the AI Act places numerous specific requirements and constraints on AI providers, and AI deployers, separately (though with some overlap). The salient point here is that AI deployers (i.e., those who design, code, train, and sell AI systems) face significant regulation under the AI Act, even if they only sell their software and never act in the capacity of a deployer, i.e., never actually use the AI themselves. And if they violate these regulations, there are severe penalties, e.g., fines of up to EUR\euro35million or up to 7 per cent of a company's total worldwide annual turnover.

And so, we again see that those who simply design and code algorithms can (and will) be subject to stringent regulation, and significant legal penalties. They will be legally accountable for their actions, even though they only `did the mathematics and/or wrote the code'.

\subsection{Legal consequences of data scraping and machine learning}

\subsubsection{Example: Scraping Tinder data}\label{s:tinder-scrape}

In 2017 a machine learning developer, Stuart Colianni, was alleged to have been involved in a mass data scrape of Tinder \cite{Lomas2017Tinder}. According to TechCrunch \cite{Lomas2017Tinder}, he scraped 40,000 images from Tinder users\footnote{Most users have multiple pictures on their profile, so this may have affected approximately 10,000 users, as an initial estimate.} in the Bay Area (California, USA) via a ``simple script'' created by ``exploiting Tinder’s API'' \cite{Lomas2017Tinder}. He then uploaded the entire data set to Kaggle, which is ``a platform for machine learning and data science competitions'' (and had been recently acquired by Google), as well as uploading the scraping script to GitHub \cite{Lomas2017Tinder}. 

At the time of writing the TechCrunch article, the data set (called ``People of Tinder'') had been downloaded over 300 times \cite{Lomas2017Tinder}. His aim was to ``use the data set with Google’s TensorFlow’s Inception (for training image classifiers) to try to create a convolutional neural network capable of distinguishing between men and women'' \cite{Lomas2017Tinder}. i.e., to create an AI for genderising. He released the dataset on Kaggle under a CC0: Public Domain License (i.e., practically unrestricted) \cite{Lomas2017Tinder}. The Terms and Conditions of Tinder at the time meant it users granted it ``worldwide, transferable, sub-licensable, royalty-free, right and license to host, store, use, copy, display, reproduce, adapt, edit, publish, modify and distribute'' user content (e.g., pictures) \cite{Lomas2017Tinder}. The platform had already suffered numerous exploits of its API in the past, and in this instance asserted that Colianni violated Sec.~11 of its Terms of Service \cite{Lomas2017Tinder}.  TechCrunch was able to locate and interview at least one of the people whose data was released by Colianni, asked whether she was ``happy at her data being repurposed to feed an AI model''; she replied with ``I don’t like the idea of people using my pictures for some sad `researches' ''.

\subsubsection{Activity: Who might be accountable for the Tinder data scrape?}

Having presented as many details as we could find about the Tinder data scrape (Section \ref{s:tinder-scrape}), we then ran an interactive activity with the students to try and understand who might be accountable? This required a bit of `setting up' of the problem.
 
First, we highlighted that some users caught up in the scrape (and many will remain blissfully unaware) were understandably aggrieved by this series of events, and may choose to take civil action in the courts for damages. Moreover, there \textit{may} have been some laws broken in the process, regarding copyright, cyber security, and data protection. And it stands to reason that some may consider the actions of various actors\footnote{A common term used in law to describe a participating entity in an event or circumstance is \textit{actor}.} in this situation somewhat immoral. We pointed out that the situation certainly seemed complicated, and that it might help to look at the three accountability principles from Section \ref{s:sorts} again: criminal, civil, and moral accountability. That is, who might go to jail, who might sue whom, and who has acted immorally and be judged by society.

Next, we though it helpful to list out all the actors involved that we could identify; there may have been more:
 
 \begin{enumerate}
  \item Stuart Colianni: The developer who did the scraping, posted the data set, and then did the AI training for a genderising algorithm.
  \item Tinder: The platform on which users uploaded their images, and whose API was used to scrape the 40,000 images.
  \item Kaggle: The platform on which the data set was later uploaded.
  \item Google: The company that provided the algorithm to do the machine learning on the data set.
  \item The 300: The people who downloaded the dataset from Kaggle.
  \item The 10,000: The Tinder users whose data (images) were scraped.
 \end{enumerate}
We noted that this list did \textit{not} include things like `The data set'; it only included \textit{actors}, which are individuals or groups with some sense of autonomy who took actions relevant to the situation.

Then we stressed that `the law' is very complicated, and while individuals may make their own decisions about moral accountability, one cannot always \textit{a priori} figure out what the position of the law is, partly because more details are often needed\footnote{A very common response given by legal experts when asked a somewhat general question about the law is `it depends''; the specific details of the circumstance can make a very big difference to how the law might treat each actor.}. To get the discussion going, we led the students with the following (non-exhaustive) list of questions that  would likely need to be resolved  before any sort of legal analysis could be done of this situation:

 \begin{enumerate}
  \item What is/are the jurisdictions? Where are the companies based (Tinder, Kaggle, Google, GitHub)? Where are the individuals based (Colianni, the 300, the 10,000)? Where is the data stored? Where did the images originate from? Tinder may be located in California, but if any of the data (images) originated in the EU then the EU GDPR will apply, which `follows' the data wherever they go. And if the AI is further developed and/or deployed in the EU, then the EU AI Act will apply.
  \item For criminal liability, what are the relevant laws in those jurisdictions, and is there overlap? If any of the data could constitute Child Sexual Abuse Material (CSAM) there are criminal offences in Europe. There are also various criminal laws which punish unauthorised access to computer material.  
  \item For civil liability, what are the damages suffered? To sue, one needs to prove damages -- unlike criminal law, which is simply `did you break the law?'. And what knock-on effects are there? If the 10,000 sue Tinder, perhaps as a class action, does that then create damages for the platform, allowing it a trigger to sue Kaggle\footnote{In one year, a student asked if Kaggle could sue Google. Seeing as Google now owns Kaggle, we suggested this might be possible from a legal standpoint, but unlikely from a practical one.}, Google, Colianni, or any other party?
  \item Are there any `hidden' parties or actors? For example, did Colianni do any work or store any data on computers at his place of work, or place of study, even if not part of his job/studies and/or done outside working hours. Can his employer/university be sued vicariously\footnote{This is when an employer can be held legally responsible for the actions of an employee carried out as part of, or within, their job. This can happen even if the employee itself is `blameless'.}?
  \item Are the 300 protected from whatever they do with the dataset because they downloaded it with a CC0 Licence? That is, can Colianni `wash' the dataset and make it permissible to use? What onus is there on the 300 to verify the source of the dataset? Can they be liable in some way?
  \end{enumerate}
From this list, students were able to make follow-up observations and ask further questions. We emphasised that the point of this example and exercise was not to establish exact legal liability in the confines of a classroom; that is for a court to decide, and a room full of mathematicians is unlikely to posses the right knowledge or insight to even know where to begin\footnote{In some years, we were able to bring in a legal scholar to the lecture. This was always quite productive, especially as students got first-hand experience of a lawyer answering `it depends' to practically every question they asked.}. Instead, it was to illustrate that the law can be very complicated, and that what mathematicians think is \textit{surely} permissible may instead land them, and those around them,  in quite a lot of legal trouble.

We encouraged the students to, as a follow-up exercise, go home and have a think about where they thought all the liability (criminal, civil, moral) lay in this example. And, if possible, we suggested they then go and talk about it with people who have the relevant expertise -- perhaps a lawyer, a philosopher, or an ethicist -- to get their opinion. Indeed, many of the students admitted they had never spoken to a lawyer (or philosopher) before.

\begingroup
\let\section\subsection
\renewcommand{\bibname}{Bibliography (Lecture 5)}
\renewcommand{\refname}{Bibliography (Lecture 5)}
\addcontentsline{toc}{subsection}{Bibliography (Lecture 5)}
\renewcommand{\addcontentsline}[3]{}
\bibliographystyle{IEEEtran}
\bibliography{bibLec5}
\endgroup

\newpage
\section{Lecture 6: Understanding the behaviour of the mathematical community}\label{s:community}

\fbox{
  \parbox{0.98\textwidth}{
Mathematicians, just like every other academic field,  form a community, maintaining and handing down their own conventions, common beliefs, and schools of thought through the generations.  This `mathematical' approach is extremely effective when working on abstract mathematics problems, but doesn't always work as well when considering broader industrial problems, even ones that involve mathematics.

This lecture looks at how mathematicians are affected by the way they are trained and taught, and how the habits, norms, and perspectives they pick up can go on to affect and inhibit the way they carry out mathematical work, potentially leading to harm. It shows how they can identify and bridge any shortcomings in their skill sets, and how to release the notion that mathematics is `pure and neutral' in order to produce good solutions to industrial or social problems.
}}

\subsection{The oil pipe problem revisited}

\subsubsection{Observing how mathematicians view the problem}

We again started the lecture by revisiting a past case study; this time the oil pipe problem from Section \ref{s:oil-pipe}. We recalled that that problem asked for the \textit{optimal} pipe route, without specifying what optimal meant, and that it was `incomplete' as posed; there were some `mathematical' omissions (some distances/costs were omitted), and some `non-mathematical' omissions (some other details about the surroundings were omitted). So we began by discussing how, over the years, students would almost always point out the missing `mathematical' details (e.g., missing distances), and say that the problem could not be solved without them. But they would seldom, if ever, ask for more details about the surrounding environment, or about other factors relevant to the problem.

We also brought up that, in a previous year, a student remarked, at the end of the activity and upon hearing the discussion about what might be in the ocean (perhaps a coral reef) or on the shore (perhaps a protected habitat, or a school), ``But then this is not a mathematics problem''. We highlighted that at no point had the problem ever explicitly been posed as a \textit{mathematics} problem.  Rather, the problem was to plan the route of a pipe; the students took it upon themselves to \textit{use} mathematics, and furthermore some specifically \textit{interpreted} it as a mathematics problem.
 
Additionally, we mentioned to students that in past years, upon seeing the problem, some students tried to race each other to see who could solve it the fastest, even though it was never presented as a race or competition, and no prize was ever mentioned or awarded. And when asked why they were racing, they often responded along the lines of `mathematics is a competition'.

\subsubsection{Analysing how mathematicians view the problem}

We moved on to discuss what is it about mathematical \textit{training} that precipitates the sorts of responses to the oil pipe problem students have demonstrated over the years. We put it to students that it is perhaps not  surprising that mathematics students responded in these ways. One key observation we presented was that students in previous years seldom asked about the surrounds, but always asked about the missing mathematical information, and that this comes from the way mathematics is taught and the way students train themselves in mathematics. Students can spend thousands of hours,  across high-school and university, completing thousands of mathematical questions on problem sheets. These questions usually have a very particular format: \textit{either} the students are given a complete, abstract, distilled piece of mathematics to solve or prove, \textit{or} they are given a (usually contrived) `physical-world' setting from which they are implicitly expected to extract an abstract mathematical problem from and produce a mathematical solution of some form -- hence, the oil pipe problem becomes an exercise in calculus. They therefore perhaps come to expect that this is the way mathematics is done, and train how to `mathematise' \textit{any} physical-world problem into an abstract mathematical one. We put it to students to think about whether they felt that at any point in their mathematical training they were ever taught to \textit{ask questions} about what they were solving, as opposed to simply taking the setting exactly as given and extracting an (abstract) mathematical problem out of it. Or moreover, whether they expect that the information provided is the `minimal' required to solve the problem, and that there are no missing details and all presented details are necessary. And we suggested that students reflect on how they might feel if they were given an over-determined problem, even an abstract mathematical problem,  with stronger hypotheses than are needed to solve it, and to think about how often they have seen this in their own mathematics training. We highlighted that, in the physical world, such over-determined problems are very common.
 
We also asked students what they thought about the statement that ``mathematics is a competition''. Several mentioned that they had participated in mathematics competitions over the years, such as various forms of mathematics olympiads. We suggested that this viewpoint, that mathematics is a competition, perhaps comes from the competitive mathematics scene that has become popular in recent decades with local, national, and international mathematics olympiads where `winners' are ordained and rewarded in mathematical contests. Of course, it is reasonable to expect anyone in a time-limited contest to  race. But we suggested that students reflect on how much time and energy they (and others) spend training for, and focusing on, competitive mathematics tournaments, and how this might lead them to believe that \textit{all} mathematics is a competition.

\subsection{How mathematical training affects the way we think}

\subsubsection{Learning more than just mathematics from its teaching}\label{s:subservience}

We used the analysis from the previous section to then lead students into the next of our theory-based ideas: that they learn more than just technical mathematics from the training they receive, and that some of these `lessons' that they learn may have a negative, harmful, and detrimental impact on how they carry out mathematical work.

We first presented some simple numbers to them. A typical one-term university mathematics subject module might have upwards of 100 `exercises' for students to do, spread across practice exercises, assessed exercises, and written exam questions (practice and assessed). In each academic year, a student on a full-time mathematics degree might study 8-10 modules, so attempt approximately 1000 exercises. Across a 3 year degree that would be 3000 exercises (4000 if they add on a 1-year masters), and that says nothing about whatever training was done in the latter years of high school; perhaps another 1000 spread across 2 years. So to a first approximation, students completing a mathematics degree attempt 5000 higher-level mathematics exercises by the time they graduate. As such, any repeating patterns, themes, or implicit assumptions present in these questions, such as stripping away background information, or always having the `exact' quantity of hypotheses needed to solve the problem, will become ingrained after 5000 repetitions. Indeed, we pointed out that unquestioningly `finding $x$' 5000 times might instil in a mathematician the notion that they should never ask `Why?' before doing so, and thus instil some subconscious degree of intellectual subservience in them when it comes to carrying out mathematical work. Such traits might then carry over to their professional career; always solving mathematics problems put in front of them because they never learned how to say no. This argument is further covered by Ernest in \cite{Ernest2018}, where he specifically discussed the impact of mathematical training on the way mathematicians think.

And so we put it to students that what they are (explicitly, and implicitly) taught, as well as what might be implicitly absent from their teaching, affects they way they view and carry out mathematics. That is, the way students are trained to become mathematicians, and the duration and intensity of this training, can change the way they think about, and view, the world. They may start to see some things as mathematical objects that others don't see, potentially `over-mathematising' situations, and their focus can move more and more towards technical abstractions that strip away subtle and/or unquantifiable aspects. In addition, they may fail to question important parts of, or indeed the underlying legitimacy of, the mathematical problems they work on; taking all such problems and tasks `as given', without question. We reminded students that they had seen this before, in Section \ref{s:covid} where we had discussed COVID modelling in care homes where the journalist asked the mathematician ``So what made you think they were shielded?'', which he answered with ``That's a good question. We never checked''.

We observed that, moreover, this intense training does not necessarily come with an `externally-facing' moral code of conduct; and that unlike most other disciplines (law, medicine, physical sciences, etc.) where students \textit{are} taught about acting professionally and ethically, mathematics training not only \textit{lacks} ethical training \cite{muller2022hippocratic}, but can (implicitly) \textit{strip away} any idea that mathematical work can have moral implications, by failing to ever mention that it does during the many years of training. Of course, mathematicians do have their own (internal) code of conduct. But unlike  codes of conduct from other disciplines which focus heavily on how the discipline affects wider society, the (unwritten) `code' for mathematicians is one that's geared towards producing and preserving `good' mathematics \cite{muller2022hippocratic, tao2007good}. Such ``socio-mathematical norms'', and the way mathematics students ``develop beliefs and values'' \cite[p.~458]{Yackel1996}, come up again later in our discussion.

But before going in to the full details about how mathematics training affects mathematicians, we first turned to a somewhat related example from history: that of the longbow archers of Henry V and their actions at the Battle of Agincourt.

\subsubsection{Example: The longbow archers at Agincourt}\label{s:longbow}

This phenomenon, of an elite group that is highly-trained in advanced technology having an immense impact on society, but without its own moral code of conduct to guide such impact, has been seen before in history. One need only consider the English longbow archers under Henry V. Archers were so important to war strategies in the middle ages that, from as early as 1363 (over 50 years before Henry V took the throne), England had enacted laws which ``commanded that on feast days every able-bodied man should in his sports use bows and arrows \ldots and shall learn and practice the art of shooting'' and that ``men were forbidden to either play or watch `games of no value' '', put in place by Edward III \cite[p.~199]{strickland2005great}. These were followed with additional similar laws by Richard II in 1388 (12 Rich.~II, c.~6) and Henry IV in 1409 (11 Hen.~IV, c.~4) \cite[p.53]{gunn2010archery}. And this was primarily to encourage the (decentralised) training of archers, to be used to ``advantage [to] the king in his actions of war'' \cite[p.~199]{strickland2005great}. Whereas it was the nobility who served as knights and other men-at-arms, archers would ``be drawn from a noble's estates'', with laws such as the 1388 Statue of Cambridge  specifically commanding that ``artisans and labourers shoot at the butts every Sunday and of feast days'' \cite[p.~199]{strickland2005great}. With archers being drawn from the lower classes, they therefore required additional encouragement and compulsion to train as archers. And so with few (legal) sports entertainment options, all men from adolescence to old age had little choice but to oblige. And indeed, they would regularly participate in competitions, with ``shooting meets [were] held every Sunday'' \cite[Ch.~6]{Mondschein2021}.
 
Years of such training didn't just equip archers of the time with new skills; it fundamentally changed them. Most obvious was the change to their physique; those who carried out the most extensive archery training, particularly with longbow and its incredible draw weight of over 80kg \cite[p.~17]{strickland2005great}, showed significant physical changes. This ranged from ``changes to their shoulder blades'' \cite[p.~6]{gunn2010archery} such as enlargement of their right shoulder as well as their left elbow \cite[p.~109]{skeletal2000}, to having a bowing of the upper arm bone (humerus) \cite[p.~109]{skeletal2000}, to having ``twisted spines'' \cite[p.~17]{strickland2005great}. And this training (and thus changes) started very early in life, with some ``display[ing] alterations \ldots that derive from strenuous activities commenced in youth'' \cite[p.~116]{skeletal2000}.

But this extends beyond physical characteristics, as ``Draw weights such as those estimated would require considerable physical skill and training, as well as mental control'' \cite[p.~109]{skeletal2000}; training not only altered their bodies, but also their minds, in the way they had to have carefully-timed control of their bodies. So how else might their minds have been altered by this? Well, unlike knights and other men-at-arms, archers trained in relative isolation; theirs was not a combat that directly engaged with their opponent -- be it a sparring partner, or enemy -- that men-at-arms would be regularly exposed to \cite[p.~290]{strickland2005great}. Their targets could be hundreds of meters away (up to 350 yards \cite[p.~295]{strickland2005great}, or even 400 yards for ``good professionals'' \cite[p.~297]{oakeshott1996archaeology}); merely a faint speck in the distance. And so these soldiers were, as we shall soon see, deadly on the battlefield, but somewhat distant and detached from the effects of battle. And, as shall become relevant later, they did not subscribe to the same chivalric code of conduct that knights and other men-at-arms did (in particular with regards to the treatment of prisoners) \cite[Ch.~11]{oakeshott1996archaeology}, as theirs was a relatively `new' branch of warfare with a shorter history and less tight-knit community.

On 25 October 1415, Henry V led his troops into battle against the French at Agincourt \cite[Ch.~17]{strickland2005great}. Coming in severely outnumbered, with only about 900 men-at-arms and 5000 longbow archers, against a French force of around 25,000 with upwards of 14,000 men-at-arms \cite[p.~325]{strickland2005great}, the odds heavily favoured the French. However, with well-trained archers who could shoot longbow in excess of 300m, and against heavily-armoured (and thus weighed-down and slow) French soldiers marching across soft damp soil \cite[p.~333]{strickland2005great} and sinking up to their knees in the mud \cite[p.~333]{strickland2005great}, Henry V was able to use this ``technological innovation'' of longbow \cite{thecollector_archery}, and more importantly those who could wield it, to his advantage. The onslaught of English arrows on the advancing French men-at-arms forced them to ``bunch inward'', and ``hampered their ability to fight'' \cite[p.~335]{strickland2005great}, giving the English an enduring tactical advantage which they capitalised on. By the end of the battle, the English has made a stunning and surprise victory, with approximately 6,000 French troops killed, and many more captured, compared to only the loss of around 600 English troops \cite[p.~338]{strickland2005great}. Superior technology in the form of longbow, coupled with superior training in the form of an entire culture built around archery, decisively won the day, where it was said that ``God and our archers made them stumble'' \cite[p.~325]{strickland2005great}.

But then something unexpected happened. Henry V, worried that the vast number of French prisoners might ``seize the weapons strewn about the field and recommence fighting'' \cite[p.~338]{strickland2005great}, ordered his men to execute their prisoners. This was unusual for the time. Prisoners of war were usually kept alive and ransomed back; this was part of the chivalric code of knights. And as such, the English knights ``balk[ed] at the order'' \cite[p.~338]{strickland2005great}, and did not carry it out as it would have broken their code of conduct. However, the longbow archers followed no such code, and so carried out the order \cite[p.~338]{strickland2005great}, executing their prisoners. For all their years of specialised training, the archers had mastered a tool of combat, without ever mastering their own conduct in combat. Unlike the knights, they lacked any extensive established community notion of what was `right' and `wrong'. This is perhaps unsurprising; when someone spend years of their life shooting something very very far away (targets, other soldiers), they establish a sense of detachment from the actions they undertake to the impact of those actions. Archers are, almost by definition, detached from where their arrows land, compared to knights who fight at very close proximity. Hence, a code of chivalry never really developed for archers. There were probably other contributing factors, such as the different social class that archers originated from. But this detachment meant there was no way such a code might naturally evolve.

\subsubsection{Some of the impacts of our mathematics training}

And so, we started to draw connections between the training and use of mathematicians, with the training and use of the longbow archers of Henry V. Like archery, mathematical training begins from an early age, and there are often nationally-mandated rules on doing such training in schools, that all pupils are obliged to participate in. This often then evolves such that, down the line, the more ``gifted and talented'' students are taken out and given more specialised teaching and training, in settings with their own particular ``socio-mathematical norms'' \cite[p.~79]{ccakir2024socio}. This training and practising of mathematics is usually done in relative isolation, save for times when they come together at various competitions. We also pointed out that doing mathematics over a prolonged period can change the way a mathematician's mind functions, by changing the way they see, interpret, and evaluate the world, as demonstrated with the oil pipe problem. And mathematicians are often situated very far from the impact of their work; through its abstraction it can be difficult to see where, and how, it is having impact, thus there is a distinct detachment from the work they do, and its effects \cite[p.~20]{muller2022hippocratic}. And the `technological innovation' created and used by mathematicians, such as the developments at Bletchley Park in breaking the Enigma code, can give countries tactical advantage in war \cite{thecollector_archery}. And finally, we highlighted that there is no agreed, or even acknowledged, externally-facing code of conduct among mathematicians governing their behaviour in the context of affecting others \cite[p.~19--24]{muller2022hippocratic}; when asked to do some mathematics, so long as that mathematics is `good' in the technical sense \cite{tao2007good}, they may be inclined to do it without deeper questioning (not unlike the archers following all of Henry V's requests), especially if they haven't received any formal or informal training in any form of critical mathematics education that prepares them for it, as discussed in  \cite{skovsmose2024critical}.

\subsection{Mathematical strengths become ethical weaknesses}\label{s:strengths}

Building on this idea that mathematical training could instil harmful characteristics in mathematicians, we started to drill down into some of the specific outcomes of such training, and how certain skills and propensities might be helpful for carrying out abstract mathematical work, but be detrimental in seeing and acting on the `bigger picture' such work might be part of. We broke this up into 4 main attributes which, in our experience, make one `good at (abstract) mathematics':

\begin{enumerate}
  \item Precision.
  \item Focused, undisturbed thought.
  \item Filtering through information to select only the relevant parts.
  \item Going extremely deep into a thought or problem.
 \end{enumerate}
None of these are specific mathematical facts or techniques that can be explicitly examined and assessed; they all form part of the implicit process of doing mathematics, and each one contributes to a student being more efficient at solving abstract mathematical problems. We pointed out that we have seen all these traits being valued by mathematicians over the years, which are trained in them from a young age, and that generally contribute to them being better at producing abstract mathematical results and solving abstract mathematical problems. We explained that mathematicians see these as their `mathematical strengths', and such \textit{skills} are what make them valuable as mathematically-trained experts, far beyond any individual theorems or proofs that they happen to know.

However, the ability to carry out abstract mathematical reasoning is not the only aspect of being a mathematician. So we articulated all these attributes to students again, but this time with a slightly different phrasing and emphasis, to see some of the other consequences of these attributes that we have observed, and our own conclusions about them.

\begin{enumerate}
  \item \underline{Precision} may lead to \textbf{\textit{inflexibility}}.
  \\If everything mathematicians do is absolutely precise, then they leave themselves no room for meaningful interpretation. Indeed, we have heard a mathematician answer the question ``Do you know the time?'' with ``Yes''. This approach is acceptable when solving abstract mathematical problems, but not every problem  solved by mathematicians is abstract mathematics. If a mathematician begins \textit{some} of their sentences with `If, then, else', then their thinking is somewhat conducive to doing abstract mathematics; if however they begin \textit{all} their sentences begin with `If, then, else', they have no way of interpreting the world in a compassionate, human way. Mathematics training might lead them to think in the `If, then, else' sort of way; the extent to which this occurs might have a negative affect on the way they interpret the world.
  \item \underline{Focused, undisturbed thought} may lead to \textit{\textbf{limited adaptability}}.
  \\If mathematicians only ever take a problem `as given', without seeking further information or letting changes in the scenario affect their thinking, they may quickly start solving a problem that is at best irrelevant to the situation at hand, and at worst quite separate and detached from its reality. Such rigidity, based perhaps on an assumption that whoever is presenting the problem is also able (and willing) to provide all relevant details, can lead mathematicians to providing unsuitable and potentially harmful solutions that don't take into account all necessary factors.
  \item \underline{Filtering through information to select only the relevant parts} may lead to \textit{\textbf{blocking out important factors}}.
  \\When mathematicians `mathematise' a problem, or indeed think about any problem, they often choose the parts most relevant to their mathematical interpretation of it. But this means they might be \textit{intentionally} ignoring other parts, some of which may have been extremely relevant to interpreting both the problem, and how to use solution. Choosing (and omitting) various aspects may make the mathematics `nicer' and more tractable, but may again lead to solutions and mathematical outputs that do not sufficiently match the situation at hand. 
  \item \underline{Going extremely deep into a thought or problem} may lead to \textit{\textbf{no breadth of perspective}}.
  \\This sort of thinking can lead to the `over-optimisation' of one aspect of a situation, possibly at the expense of others. If mathematicians do not maintain a broad perspective, and cannot see the way their work or role fits in to the larger picture, it may significantly distort the suitability of their output. Going too deep into a logical chain of reasoning can lead to a `logical' output that might be absurd, as even the smallest oversight or misunderstanding of what is truly important can shift such logical reasoning in an inappropriate direction. As a toy example, being tasked to `water the flowers' and setting up a regular and thorough irrigation system is great, unless of course one is in a drought; mathematicians may do very well watering the flowers, but fail to see the drought.   
 \end{enumerate}
And so we suggested to students that all the traits, skills and training that make mathematicians very good at doing abstract mathematics, may well be the very things which make them easily exploitable, make them somewhat blind to world around them, and make them produce results and analysis that is out of context or out of touch with reality. The greatest strengths of mathematicians may simultaneously be their greatest weaknesses.

We connected this to earlier discussions, pointing out that  this, in turn, helps to explain \textit{why} several of the examples in earlier sections occurred, giving a  \textit{source} of the behaviour and actions of mathematicians in those scenarios. Those earlier sections went into great detail about \textit{what} mathematicians did that caused harm, and now students started to see \textit{why} mathematicians may have made those choices, and subsequent mistakes. Ultimately, the viewpoint we put forward here is not just a diagnostic one, but moreover preventative; by understanding what mathematicians did that caused harm, and why they did it, students have an effective mechanism to scrutinise their own actions in a way that might prevent \textit{them} from going down the same path and causing harm.

\subsection{Motives for teaching mathematics}\label{s:teaching-motives}

\subsubsection{Who do we teach mathematics for?}

So how does this happen? How is it that the teaching of mathematics often focuses on maximising the ability of the students to carry out abstract mathematics, with minimal consideration for what effect that might have on their outlook on the world, or what sort of skills are \textit{not} trained in to them, or indeed trained \textit{out} of them? 

We suggested to students that, to understand  the teaching process they are currently in, they must first look at what motivates the teaching methods, and curriculum selection, of the academics carrying it out. So we gave the following anecdotal observation, which is argument which we have heard (in various forms) from academic mathematician many times, especially those from `prestigious' institutions\footnote{And these are not just our observations; we have spoken to other mathematicians who have heard the same idea articulated by academics at institutions we have never visited.}: 
\begin{center}
\textit{We teach mathematics for those  students who might come and replace us as academics. We happily take on many other students who won't go in to academia, and they are welcome to participate, provided they can keep up with the material. But they are not our target audience. We design and deliver our courses, we choose what to teach, what not to teach, and how to teach it, purely with the interests of those  future academic mathematicians in mind.}
\end{center}
Of course, we stressed that this is not the opinion or approach of all academic mathematicians; far from it. One can find academics from all backgrounds, at all institutions, who are not of this mindset, and who have various different motives and reasons for teaching mathematics \cite{ernest2010teach}. But this somewhat approximates a commonplace opinion, especially for those mathematics courses which are not service courses (i.e., courses for `mathematicians', rather than courses for `mathematicians + engineers + biologists + \ldots'). And we highlighted that this approach starts to show \textit{why} much mathematics teaching aims to  maximise abstract mathematical skills, while giving minimal regard for other skills, or the effects of such over-optimisation; a mathematician in academia will encounter fewer (but certainly not none) of the safety and ethical issues that one in industry might. But perhaps this teaching approach is somewhat justified \textit{if} the vast majority of mathematics students do go on to a full career in academia? So we proceeded to `explore the numbers' with the students, to see if this was indeed the case.

\subsubsection{Example: Who are we teaching mathematics to? Attrition in academia}

One can start to look at various statistics regarding how many students progress to careers in academia. In the UK, the Royal Society carried out such a study in 2010 \cite{RoyalSociety2010}, taking the cohort of all \textit{PhD graduates} (i.e., those who complete their PhD) across all science, technology, engineering, and mathematics (STEM) subjects, and recording their career trajectories. So, of the full cohort of STEM PhD graduates within the UK, their careers progress roughly as follows \cite[p.~14]{RoyalSociety2010}:

\begin{itemize}
  \item 53\% leave academia at the end of their PhD; 47\% continue on to early-career academic research (e.g., find a first postdoctoral position).
  \item 17\% leave academia for industry after doing some number of years of early-career research; 30\% secure a permanent academic position.
  \item 26.5\% leave academia for industry at some point in their career after finding a permanent academic position.
  \item And thus 3.5\% spend their entire career in academia.
  \item Furthermore, (only) 0.45\% become university professors at some point in their career.  
\end{itemize}
Keep in mind that these are the numbers for students who \textit{complete a PhD}; the PhD programmes also have some level of attrition. And so if one traces this back to those who start as an undergraduate, the proportion of those who go in to industry are probably an order of magnitude larger than the figures stated above. So as a first approximation, under 1\% of undergraduate mathematics students spend their entire career in academia. Of course, these numbers are not exact; they only deal with UK PhD graduates, and cover STEM subjects in their entirety (so not specific to mathematics). But they point to the fact that only a (very) small proportion of mathematics students spend their entire career in academia; the (vast) majority will at some point end up in industry, and thus face the full swathe of safety and ethical issues present therein.

There is another way to (partly) verify the numbers given above, in particular with regards to academic mathematics careers. It is possible to make use of the Mathematics Genealogy Project \cite{genealogy} to look up academic mathematicians and see their ``academic descendants'' (i.e., their PhD graduates, the graduates of their graduates, etc.). This is the `family tree' of academic mathematics, and the website has quite thorough data going back to the middle ages. The key assumption here is to use the fact that a mathematics PhD graduate having had \textit{some} PhD graduates of their own as an indicator or proxy for having gained a permanent academic position. And the more such students they have, or the longer the time window the students are spread out over, the longer that person's career in academia. Of course, this is not a mathematical rigorous process, as not all PhD graduates are recorded there, not all academics take PhD students, and no full statistical analysis of all the available data has been done. But some random sampling of mathematicians can still start to provide further insight.

If the `attrition rate' of mathematics PhD graduates -- from graduation to securing a permanent academic role -- was very low, what one might expect to see is that, if an academic mathematician has $N$ PhD graduates, and $k$ descendants in total, then each of their PhD graduates would have approximately $k/N$ descendants. But actually, what we found in our (ad-hoc) random sampling of postwar-mathematicians on Mathematics Genealogy is that most (approximately $3/4$) of their PhD graduates have $0$ descendants each, and the remaining $1/4$  have approximately $4k/N$ descendants each. That is, $3/4$ of their graduates left academia for industry rather early in their career. This is fairly close to the statistic given by the Royal Society of 30\% of STEM PhD graduates obtaining a permanent academic position. Even if that full cohort of mathematicians -- the `30\%' who gain a permanent academic position -- were to stay in academia for their entire career\footnote{We have not studied the numbers further on Mathematics Genealogy to analyse which mathematicians had a long history of PhD graduates and thus presumably spent their entire career in academia.}, rather than dropping to 3.5\% as is the case in general for STEM, that still implies that less that 10\% of mathematics \textit{undergraduates} spend their entire career in academia. And so the conclusion remains that the vast majority of mathematics students will at some point end up in industry.

As an aside, when presented to students, this particular example becomes an instance where they need to use intuitive mathematical reasoning to analyse a moral question before they know the necessary mathematics; in this case, branching processes and statistical hypothesis testing. It is perfectly possible to reason about such a scenario, even when one doesn't have all the mathematical details at hand.

\subsubsection{The under-serviced majority of mathematics students}\label{s:underserviced}

Combining these observations about \textit{why} mathematics is taught (i.e., in significant part for future academics), with the analysis of \textit{where} students end up in their career (i.e., overwhelmingly unlikely to be academics), we were able to demonstrate that many students of mathematics might actually be \textit{under-serviced} by the teaching they receive; the focus on abstract mathematics comes about in part from the disregard for careers that deal with factors other than abstract mathematics. The trade-offs described in Section \ref{s:strengths}, whereby strengths in abstract mathematics can become weaknesses when using mathematics, are only minimally recognised in academic work. As such, it becomes easy to ignore these trade-offs when directing teaching to `future academic mathematicians'. We explained that, for now, students should anticipate that mathematics teaching will continue in this way; focusing on preparing future academics, while (perhaps accidentally) quashing or training-out skills and perspectives that are otherwise very useful for the majority that proceed to industry.

Of course, just because a mathematician leaves academia, doesn't mean they're lost to society; far from it. We explained to students that the mathematics education system \textit{as a whole} is set up to produce many mathematical experts, who deploy their knowledge and skills elsewhere outside of academia. And such support and encouragement often comes from the very highest echelons of government. For example, in 2010 the Australian Federal Government declared mathematics and statistics areas of ``national priority'', and thus (alongside science) reduced university fees in these areas to the lowest rate across all fields of study, 19\% cheaper than the next highest ``band'' \cite{goingtouni2011whatyoupay}. And that has continued in Australia to the present day, with 2027 mathematics/statistics degrees still on the lowest possible fee level; now 50\% cheaper than the next highest ``band'' \cite{studyassist2026studentcontribution}. This goes some way to illustrate that the training of people in higher-level mathematics is funded for a variety of reasons, including economics and technological progress. 

And so we highlighted that having a mathematically literate populace, and a lot of mathematical experts, is a way to generate economic, and  political, power within society \cite{Keitel2015}. But this is already well-known by people other than mathematicians; governments are well-aware. Thus, there are  arguments suggesting that the ``curriculum is a mechanism of governmental control'', that ``teachers are slaves of the textbook'', and that `` assessment
is still primarily and predominantly a mechanism for selecting the mathematical
elite only'' \cite[p.~167]{Keitel2008}. Its teaching ``(re)produces techniques of power which privilege some children while disadvantaging others'', and mathematics itself is a ``knowledge which allows power,
represents a specific world view and serves the interests of certain groups in society'' \cite[p.~1061]{kollosche2014mathematics}. We clarified to students that such arguments are suggesting the same common conclusion: that mathematics creates and serves power, hence why it is so well supported by those in power.

We concluded by arguing that, in these places outside academia, mathematicians deploy mathematics as a means towards different types of goals to those in academia, and hence the trade-offs of Section \ref{s:strengths} can become quite problematic if the mathematicians are not aware of them. Such mathematical skills, and knowledge, combine to create tools of (vast) power, not unlike what was seen at Agincourt as described in Section \ref{s:longbow}. And so, while students might be under-serviced \textit{overall} in the mathematics teaching they receive, at the same time society has (historically) been very happy to embrace their (technical) mathematical expertise, and found a lot of ways to make use of them, sometimes to the point of exploitation. However, such mathematical work is sought for these purposes; economic gain, power, etc. So mathematicians need to properly inform themselves of \textit{who} is paying for their work, and \textit{why} they're paying for it \cite[p.~16]{chiodo2023manifesto}. Ultimately, when someone else is making use of mathematical skills and knowledge, mathematicians need to know how to navigate that world, and unfortunately they still remain under-serviced with regards to this sort of learning during the teaching process.

\subsubsection{Who teaches mathematics?}

At this point, we shifted the discourse from \textit{what} is taught, to \textit{who} teaches it. And the above analysis of mathematicians on the Mathematics Genealogy Project demonstrates something else: From there, we see that mathematicians are taught by a very small proportion of mathematicians, in a very tight-knit `family'. As described by Davis and Hersch in ``The ideal mathematician'' \cite{davis1981mathematical}, they write that 

\begin{quote}
    ``The ideal mathematician's work is intelligible only to a small group of specialists, numbering a few dozen or at most a few hundred. This group has existed only for a few decades, and there is every possibility that it may become extinct in another few decades.'' \cite[p.~34]{davis1981mathematical}
\end{quote}
This demonstrates that such a family often writes for, and indeed looks out for, itself, and itself alone. 

We then explained to students that this small proportion is very particular, in that those in it will have very little (if any) experience outside academia. Unlike other disciplines, it can be extremely difficult for a mathematician to leave academia, spend some time in industry (or elsewhere), and then return to an academic career. Of course it is not impossible, and there are many cases when it does happen. But for the most part, academic mathematicians have spent their entire career inside academia. Thus, the perspective of the world that they do, or even can, teach, is one that is developed from within academia and with minimal influence from the sorts of industries and workplaces where in fact most of their students end up. We clarified that, of course, this is a similar situation in many academic fields, but that those with an explicitly more applied focus such as engineering, management, computer science, etc.~often have a more porous boundary. Nonetheless,  there's a lot that can be learned from such adjacent (mathematically-powered) disciplines, as outlined extensively in \cite{muller2022hippocratic}.

\subsection{Handing down knowledge and culture}\label{s:hand-down}

\subsubsection{The way mathematicians hand down knowledge}\label{s:handing}
 
We then turned our attention to what exactly it is that academic lecturers, teachers, and instructors `hand down' or impart to their students, and the particular type and style of knowledge this is. We started with  Isaac Newton's letter to Robert Hooke \cite{Solr-dc-9792} in which he wrote ``If I have seen further it is by standing on the shoulders of giants'', which partly explains how mathematics is done; mathematicians build on the work, and knowledge, of their predecessors. Of course, all disciplines do this, to an extent. But we stressed that mathematics is particular in that there is rarely a  replacement or deletion process in the knowledge transfer, and it is often strictly cumulative. Other disciplines tend to throw out, or amend, ideas, whereas the mathematics that is provably true today is provably true tomorrow (or until our standards of proof change, which when they do can create significant levels of ``anxiety'' in the mathematical community, as happened in the 19th century \cite{Gray2004}). We pointed out that what generally changes are techniques, and moreover communication and writing style, but not so much the abstract technical knowledge. And by looking at the Mathematics Genealogy Project \cite{genealogy}, one can trace precisely the way in which such (deep) knowledge is `passed down' in mathematics through the generations, from supervisor to PhD student in a way not too dissimilar to a baton relay. 

However, we stressed that this `passing down' can cut both ways; productive, as well as harmful, insight can be passed down. We gave the example of Andrew Wiles, who is well-known for proving Fermat's Last Theorem. His PhD supervisors (and theirs, and theirs, etc.) can be traced back 6 generations on the Mathematics Genealogy Project \cite{genealogy_wiles} to William Rouse Ball\footnote{The full chain of supervisors is as follows: William Rouse Ball $\to$ Ernest Barnes $\to$ John Littlewood (close collaborator of Godfrey Hardy) $\to$ Harold Davenport $\to$ John Coates $\to$ Andrew Wiles.}; indeed, all of them were at the same university (Cambridge). At the opening of the new mathematics building at the University of Oxford in 2013 -- aptly named the ``Andrew Wiles Building'' -- Wiles gave a speech in which he commented on the harmful use of mathematics in finance, and said ``One has to be aware now that mathematics can be misused and that we have to protect its good name'' \cite{silvera2013ibtimes}. We pointed out that while it was helpful that Wiles publicly pointed out this `harm' of mathematical work, he didn't mention the \textit{social harm} that was caused; just the reputation harm to mathematics as a discipline; it's ``good name'' seemed to be his primary concern. We suggested to students that someone living in a world of academic mathematicians might end up over-valuing the ``good name'' of mathematics compared to, say, wider societal harm. In addition, Wiles also said that ``Mathematics has lost its moral purity'' \cite{devlin2013financial}, and so we asked students to consider where he might have picked up that notion from, and whether they had see such a notion -- `moral purity' -- written out and argued \textit{explicitly} in any books, papers, or other resources they have come across. Again, we suggested that such a concept may have been `passed down' through the generations; not necessarily written down, but nonetheless conveyed during the teaching and mentoring process.

And so, we put it to the students that it is not just \textit{mathematical understanding} that is passed down through the generations, but also an understanding of \textit{what mathematics is} (and is not). A theorem can be true, and mathematicians can provide a proof for it; that is universal. But for a theorem to be ``beautiful'', or ``elegant'', or even ``trivial'', is a different (and subjective) notion \cite[p.~625]{tao2007good}. Yet, we posited that such ideas can also be `passed down', alongside mathematical truths and techniques. So where might mathematicians get notions such as `mathematical beauty' from? Well, we see this in ``A mathematician's apology'' by Godfrey Hardy \cite{hardy1992mathematician}. He uses the word beauty/beautiful over 20 times in his short 50 page book; it appears on almost every open double-page. Hardy is clearly expressing this very human notion of what mathematics `is'. Thus, not only has Hardy handed down knowledge, but also (his) understanding of how one might view, appreciate, and be motivated to do, mathematics\footnote{Recall, moreover, that Hardy was a close collaborator of Littlewood, who was in the `line' of mathematicians leading to Wiles.}. Indeed, Terrance Tao does likewise, in his work titled ``What is good mathematics?'' \cite{tao2007good} where he also uses the term ``beauty'' and its variants 3 times as descriptors, and even starts to provide an explicit definition of it \cite[p.~624]{tao2007good}. We therefore suggested to students that all this is far more than just knowledge transfer; this is a transfer of culture. Of course, the setting of `boundaries' of what counts and doesn't count as science (or mathematics) is normal for all scientific disciplines \cite{gieryn1999cultural}. Mathematics is not special here, and so the framing of what makes for `good' mathematics is just part of the process of doing scientific work. Our point here to students was no that this \textit{shouldn't} happen, but rather that it \textit{does} happen, and they should be aware of it and simply be mindful of the potential consequences of how it is done.

\subsubsection{The community and culture of mathematicians}

We used this discussion of mathematical teaching to then raise the idea with students that mathematicians form a community, and have their own culture, with all sorts of aspects that cultures typically demonstrate, including norms, language, heroes, relics, history, and even an appreciation of aesthetics and what is `beautiful'. Not only did we assert that such a culture exists, but also that it has various impacts on who does mathematics, how they do it \cite{Burton2009}, and indeed who even counts as a `mathematician'; the definition of which is far from clear \cite{buckmire2023definitions}. We presented to students some of the specific ways in which we saw these aspects manifest in mathematical culture, as follows:

\begin{itemize}
 \item \underline{Norms:}
 \\Mathematicians often write their documents in \LaTeX, usually in the same font and style, and have even been known to suggest that a piece of mathematics not being written \LaTeX\ is the \textit{first} of ``10 signs that [it] is wrong'' \cite{aaronson2008tensigns}. They have particular ways of laying out proofs, and often place a small square in the bottom-right corner to signify the end. They often insist on using blackboards \cite{zenou2024magic,uchicago2021boardwork}, even though  they rebelled against them in the early 19th century when they were first introduced; the famed ``Conic section rebellion'' \cite{jackson2002communications}.
 \item \underline{Language:} 
 \\Mathematicians have their own particular methods and style of written technical language, which has been studied extensively \cite{ganesalingam2013language}. In addition to this, they export some of these terms into their semi-technical and social vernaculars. For example, they might use a phrase like ``That proof is trivial'' to mean it is (in some sense) easy, or say ``Those two words are a small Hamming distance apart'' to mean they are spelled similarly. This can even leach into their everyday conversation, when they may say things like ``The intersection of their friendship circles is non-empty'' to mean two people have a mutual friend, or ``$5$ of the people were late, and $N-5$ were on-time'' to mean that five were late and the rest were on time.
 \item \underline{Heroes:}
 \\Mathematicians have among themselves certain `heroes' \cite[p.~11,14]{doxiadis2003}; those (past and present) who have reached near-idol status in the community. These are prize-winners such as Fields Medallists \cite{imufieldsmedal}, as well as those who have proven remarkable mathematical results before the times of such awards. Similarly, they have their own `sporting heroes' in the form of Mathematics Olympiad medallists \cite{imoresults}. And they even \textit{name} (important) theorems after whoever first proved them as a sign of respect; even though a more `logical' (i.e., `mathematical') approach might be to use a numbering or descriptive naming system, they choose not to.
 \item \underline{History:} 
 \\Mathematicians have a long and detailed history of mathematics and of mathematicians, going back several millennia to the Babylonians and Egyptians \cite[p.~5--10]{cajori1893history}. So long is its history, that history books such as Cajori's \cite{cajori1893history} were written about it as far back as 1893, with older books such as ``The Life of Mr.~Paschal'' coming from 1744 \cite[p.~xi]{cajori1893history}. And it is seen as important, with William Glaisher saying ``I am sure that no subject loses more than mathematics by any attempt to dissociate it from its history.'' \cite[p.~iii]{cajori1893history}; perhaps a peculiarity, given mathematical \textit{truth} could naively be considered independent from its \textit{historical origins} (but of course, is not).
 \item \underline{Relics:}
\\Mathematicians keep notable relics in museums and libraries. Things such as Ramanujan's notebooks (now in the library of Trinity College, Cambridge \cite{papers_ramanujan}), original printings of Newton's Principia (various locations), and ``Plimpton 322''; a 4000 year old Babylonian clay tablet with Pythagorean triples on it \cite{abdulaziz2010} now in the Rare Books and Manuscript Library of Columbia University \cite[p.~11]{vinciguerra2017}. Mathematically-speaking, there is no need to keep these; we have these results written elsewhere now. But we do, because they form part of the cultural history of mathematics. And mathematicians find great excitement in visiting relics such as Principia and Plimpton 322 in-person, calling them the ``slam-dunk'' of a museum visit, and commenting ``we were all a bit shocked that we’d just been centimeters away from this mysterious, irreplaceable historical document'' \cite[p.~12]{vinciguerra2017}.
\end{itemize}

\subsubsection{What exactly are mathematicians handing down?}

We drew on the previous points, regarding how mathematicians hand down knowledge, and the cultural aspects of mathematics, to suggest that the handing down of (technical) mathematical knowledge also comes with it a handing down of the cultural aspects of the mathematics community. We used these earlier points to highlight that this is inevitable; the fact that mathematical teaching is a human endeavour means that the cultural aspects of the discipline cannot be detached from the technical aspects during the teaching process. Learning is not limited to what is explicitly taught, and so even that which is mentioned only implicitly, and that which is omitted, is also picked up in the learning process. And so we hinted to students that not only does the teaching process contribute to them knowing what their lecturers know, but also to them \textit{being} how their lecturers \textit{are}.

But just as this emulation and cultural propagation can be very useful in terms of knowing how to do abstract mathematics, we emphasised to students that it can also inhibit their understanding and perception of their own role in society. Mathematicians hand down an understanding not only of what mathematics is, but also of what mathematics might, or might not, be used for; this may not always be correct, as it does not follow the same formal proof structure as abstract mathematics. For example in his ``Mathematician's Apology'' \cite{hardy1992mathematician}, Hardy writes:

\begin{quote}    
``The great modern achievements of applied mathematics have been in relativity and quantum mechanics, and these subjects are, at present at any rate, almost as ``useless'' as the theory of numbers. It is the dull and elementary parts of applied mathematics, as it is the dull and elementary parts of pure mathematics, that work for good or ill. Time may change all this. No one foresaw the applications of matrices and groups and other purely mathematical theories to modern physics, and it may be that some of the ``highbrow'' applied mathematics will become ``useful'' in as unexpected a way; but the evidence so far points to the conclusion that, in one subject as in the other, it is what is commonplace and dull that counts for practical life.'' \cite[p.~39]{hardy1992mathematician}.
\end{quote}
Hardy wrote his Apology in 1940, only a couple of years before the world-changing wartime projects at Bletchley Park (using number theory to break the German Enigma code \cite{gannon2014colossus}) and the Los Alamos Laboratory (using quantum mechanics to develop the first nuclear weapon \cite{szasz1992british}) were in full swing. We pointed out that of course one would not \textit{expect} Hardy to have foreseen such developments; however, to categorically rule them out was perhaps quite telling of a particular example of cultural understanding among mathematicians at the time. And it is \textit{that} understanding, of what mathematics is, and what it is not (in the eyes of the community of mathematicians), that also gets handed down to future generations, both of which can significantly impact the way mathematics is thought about, and done, in highly contextualised (and potentially high-impact) situations.

\subsection{Motives for doing mathematics}

\subsubsection{The `beauty' of mathematics}

Having already established that mathematicians use the term `beautiful' and its variants to describe mathematics (Section \ref{s:handing}), we started to explore how this might actually be a motivator for mathematicians to do mathematics. But first, we looked into what `beauty' actually means in the context of mathematics. Tao \cite{tao2007good} begins to describe `beauty' in mathematics as:

\begin{quote}
    ``Beautiful mathematics (e.g. the amazing identities of Ramanujan, results which are easy (and pretty) to state but not to quote).'' \cite[p.~624]{tao2007good}
\end{quote}
Already one notes that this description includes the words `amazing' and `pretty', none of which are described or appear elsewhere in that work. Indeed, Tao's description of beauty is but one of a ``non-exhaustive'' list of 21 qualities of what makes ``good mathematics'' \cite[p.~623--624]{tao2007good}. Moreover, he introduces other aesthetic and somewhat hedonic\footnote{i.e., devoted to pleasure.} terms in that list such as taste, elegance, creativity, depth, and insight \cite[p.~35]{muller2026ethical}. Indeed, Tao also comments on one quality in his list, rigour, by saying ``mathematical rigour, while highly important, is only one component of what determines a quality piece of mathematics'' \cite[p.~624]{tao2007good}. And so we suggested that students consider the juxtaposition between the idea that mathematics is valued and revered for its rigour, and Tao's list which suggests that rigour is but one of many elements of what makes ``good mathematics'', and that all but one (i.e., `$N-1$') of the qualities on that list (rigour) are ``somewhat subjective and contain some inherent imprecision or uncertainty'' \cite[p.~625]{tao2007good}. We used this to demonstrate to students that the main selling point of mathematics -- its rigour -- might not actually be the main reason people \textit{do} mathematics, and that almost every other reason is actually ``somewhat subjective''.

We then told students that our focus on mathematicians' pursuit of beauty was to demonstrate that such pursuits may actually lead to the production of inappropriate or harmful mathematical work; indeed, a mathematician's definitions of beautiful and good mathematics can often be quite disconnected from ethics and morality \cite{Müller2024}. As an initial example, we presented the idea of being ``Lost in math: How beauty leads physics astray'' \cite{hossenfelder2018lost}; a book by physicist Sabine Hossenfelder that explores how the use of mathematics in physics, and how physicists seek `beauty' in the mathematical expressions they derive to describe the physical universe, might actually impede the production of new physics insight by limiting what physicists are willing or able to consider. She describes how ``aesthetic judgment drives contemporary research'' and that ``When asked to judge the promise of a newly invented but untested theory, physicists draw upon the concepts of naturalness, simplicity or elegance, and beauty'' \cite[p.~9--10]{hossenfelder2018lost}. She quotes the head of theoretical physics at CSER, Gian Francesco Giudice, as saying ``The sense of beauty of a physical theory must be something hardwired in our brain and not a social construct'' \cite[p.~12]{hossenfelder2018lost}, and one of the founders of string theory, John Schwarz, as saying ``The mathematical structure of string theory was so beautiful and had so many miraculous properties that it had to be pointing toward something deep'' \cite[p.~12]{hossenfelder2018lost}. Yet she then laments that ``physics isn’t math'' and that ``The main reason we use math in physics, however, is because we can'' \cite[p.~14--15]{hossenfelder2018lost}, having given this pursuit of mathematical beauty as a reason for why ``for more than thirty years now we have not been able to improve the foundations of physics'' \cite[p.~13]{hossenfelder2018lost}. And so we used this to demonstrate to students that just because some mathematics is (subjectively) beautiful, does not mean it fulfils any higher, let alone useful, purpose.

We stressed to students that this propensity to use `beautiful mathematics' in theoretical physics was, although harmful to the progress of the discipline itself, not outwardly harmful to wider society, and that this was because physics still followed the scientific method of experimentation for verification and refutation. Such misleading (but beautiful) theories in physics could either be disproved, or were lacking in predictive power; they didn't do anything. But we then reminded students of some earlier examples from finance where mathematical beauty was sought, and the results then used, in an unscientific, unfalsifiable manner (that is, until grave harm was caused). The first of these was the Black-Scholes Equation for pricing options (Section \ref{s:blackscholes}), with its `nice' assumptions about the market, and subsequent `nice' equation (the heat equation). We pointed out that here, the beauty made the equation analytically solvable, though (unfortunately for many) highly ineffective and ultimately harmful. The second of these was Li's One-factor Gaussian Copula Model, for pricing CDOs (Section \ref{s:copula}). Again, this was a `beautiful' model that offered  high analytical tractability, but again without sufficient physical-world applicability to avert economic disaster.

Building on all these examples, we made the point that while the notion of beauty might be interpretable by mathematicians within the realm of abstract mathematics, it has no bearing, relevance, or credence in the physical world. We stressed that those who pursue `mathematical beauty' when trying to use mathematics to explain and understand the physical world, risk being ineffective, or even causing immense harm. We finished by articulating our own rule of thumb for mathematicians:
\begin{center}
\textit{If you can explain the behaviour of a system containing over 1 million parts, with an equation that fits in 1 line, then regardless of how beautiful your formula is, you've probably (deliberately, or accidentally) missed something.}\footnote{Even in very predictable physical systems from thermodynamics or statistical physics, such formulae might capture \textit{something} about average-case behaviour, but certainly not everything.}
\end{center}

\subsubsection{Why we do mathematics}\label{s:whymath}

Building on this idea of `beauty in mathematics' and all its apparent consequences, we suggested to students that perhaps a large proportion (not all, but many) of mathematicians study and develop mathematics primarily for its intrinsic beauty. Indeed, they may be doing so even though there might be very few other people who could possibly appreciate the `beauty' that they are interpreting, as even within the mathematical community there are only a small number of people specialised enough to understand certain areas; from outside the community of the mathematically-trained, practically no-one. On top of this, we pointed out that in our experience many mathematicians  claim that they derive a particular sense of enjoyment and excitement out of doing mathematics, producing theorems, and solving problems, rather than being  motivated to do mathematics for the direct uses, applications, or societal benefit of their work (again, some do, but many do not). We acknowledged that all these are perfectly fine, legitimate, and honest reasons for doing mathematics, and there is no immediate problem in being \textit{motivated} by these. However, as illustrated in the previous section, if these motivators become the main -- and indeed only -- driving force for doing mathematics, then other important factors can start to be overlooked, and end up making the mathematics much less helpful, or indeed harmful, for the situation at hand.

We further suggested that these sort of motivations for doing (abstract) mathematics -- revealing beauty, enjoying the process, etc.~-- sound very similar to another area that people invest a lot of time and effort in to: (abstract) art. Those making art, be it painting, sculpting, music, dance, etc., do so for several of the same reasons that many mathematicians seek to do mathematics. That is, they see it as a form of beauty, and they enjoy the process. And just like mathematics, some art forms, such as ballet, are difficult to appreciate if one does not have the right understanding and training. There may be very few, if any, applications of the art people produce. But all in all, people still do, and enjoy, art. And indeed, the question of ``whether mathematics is an art or a science'' has already been raised and explored by others such as Rowe \cite{rowe2002years}.

However, this raises an interesting question: who pays for people to do art, or research in the arts? In the UK, the government body tasked to provide research funding for the arts is the Arts and Humanities Research Council, with an annual budget of GBP\pounds 70million in 2025-26 \cite{DSIT2026RD}. So we put it to students that if the act of doing mathematics is so similar to the act of making art, why do mathematicians not seek funding from government arts councils? Why, instead, do they (primarily) seek funding from science research councils, with substantially larger budgets than arts councils, or look for work in industries such as technology and finance (with much higher salaries than most artistic endeavours)? We suggested this seems like a logical decision by mathematicians: follow the money, to do the most mathematics. And indeed there is much more money, as the UK government body for funding the sciences, the Engineering and Physical Sciences Research Council, had an annual budget of \pounds 640million in 2025-26 \cite{DSIT2026RD}; over 9 times the Arts and Humanities Research Council budget.

We asked students to think about why, if what mathematicians do seems so much like an art form, detached from physical-world application in both its formulation an execution, it still gets funding form science bodies and big industries? We suggested that actually, the answer could be rather simple: because these bodies, and companies, believe that the mathematics done will, over some time frame (short, or long term), eventually be helpful and of some use. That is, such entities fund mathematicians because they're making things, doing things, producing things, that are, or might be, useful to society, often through connections to the (other) sciences. Society, in the form of government agencies as well as private industry, are most probably paying mathematicians to do mathematics because it is useful, rather than to entertain their desires for enjoyment and for the production of mathematical beauty, even though many in the community have adopted the belief that their work is `artistic' in the sense that it is done for the sake of beauty and its formulation detached from physical-world application, with any `use' that comes from it being a secondary effect.

This convenient coupling between mathematics and sciences has been known for decades, going as far back as at least 1960 with Nobel laureate Eugene Wigner's paper ``The unreasonable effectiveness of mathematics in the natural sciences'' \cite{Wigner1960}. There, he argues that mathematics is ``the science of skillful operations with concepts and rules invented just for this purpose'' \cite[p.~2]{Wigner1960} and that ``Most more advanced mathematical concepts \ldots were so devised that they are apt subjects on which the mathematician can demonstrate his ingenuity and sense of formal beauty'' \cite[p.~3]{Wigner1960}; again the term `beauty' comes up. Conversely he argues that ``The physicist is interested in discovering the laws of inanimate nature'' \cite[p.~14]{Wigner1960}, and eventually concludes by saying ``The miracle of the appropriateness of the language of mathematics for the formulation of the laws of physics is a wonderful gift which we neither understand nor deserve'' \cite[p.~14]{Wigner1960}. This connects back to Section \ref{s:limitsmodelling}, where we discussed the falsifiability of mathematics when used in physics (but not in other areas). And Karaali summarised this overarching idea in 2019, writing ``[mathematicians] have claimed to be useful mathematicians doing useless mathematics.'' \cite[p.~12]{karaali2019doing}.

We put it to students that ultimately, mathematicians' own beliefs of what mathematics \textit{is}, do not dictate the way in which it is used, and exploited, in the world. Moreover, such beliefs need not even affect how those from outside the mathematical community think about, or view, mathematics and what it fundamentally \textit{is}; what mathematicians think does not always need to match what other people think or do.  Mathematicians may think that their mathematics has no potential for harm because their motivations for doing it are detached from reality, but those paying them certainly have good reason to pay for it, as explained further in Section \ref{s:underserviced}.

\subsection{The (fake) objectivity, purity and neutrality of mathematics}

While objectivity, neutrality, purity and the universality of mathematics often represent the starting point of the academic literature and ``grown up debates'' among philosophers \cite[p.~22--23]{muller2026ethical}, our exposition in this lecture was deliberately structured slightly differently. We wished to discuss such concepts in light of the experience and `intuition' in ethics in mathematics that our students have gained up until now. Our goal, over this and the preceding lectures was to get students to \textit{reach} the idea that there is no objectivity, neutrality, or purity in mathematics, by way of many examples and reasoned arguments. And so, this final section of Lecture \ref{s:community} represented the conclusion of such reasoning  for the students, bringing them (gently) to the point that most philosophers would instead have started at\footnote{`Started', modulo their many years of training in philosophy, much like a trained mathematician might say a mathematical proof is `trivial' modulo their own extensive training.}.

\subsubsection{Truth vs meaning}

So, with the previous discussion of mathematics in the sciences in mind, we then turned to the concepts of \textit{truth} and \textit{meaning}. This was based on numerous experiences of hearing mathematicians argue that mathematics gives `absolute truth'. So we put these two ideas to students, with the following maxim:
\begin{center}
    Mathematics on its own gives \textit{absolute truth}, but with \textit{absolutely no meaning}.
\end{center}
That is, meaning must be \textit{inferred} from mathematics, done as a human endeavour, and with all the subjectivity, shortcomings, and drawbacks that any other human endeavour must endure and work with. That is not to say that meaning cannot be inferred from mathematics in a useful and responsible manner; rather, it says that the \textit{absoluteness} present in abstract mathematics melts away as soon as it is `touched by human hands'.

As an initial example, we asked the students a simple question: ``What is 2+3?'' After some pondering\footnote{After 6 lectures on Ethics in Mathematics, the students were somehow convinced that this was a trick question.}, we eventually received the answer: ``3''. So we then posed a physical-world problem to try and infer meaning from this absolute truth: ``If I have $\$2$, and I go out and earn another $\$3$, how much money do I have?'' The students, now more confident in their mathematical prowess, were happy to answer ``$\$5$''. When asked why, we usually received responses of the form ``Because $2+3=5$''. At this point we reminded the students of the oil pipe problem from Section \ref{s:oil-pipe}, and suggested that they had been misled by their (unquestioning) interpretation of the problem, and by their faith in mathematical truth; we pointed out that we are subject to a marginal tax rate of about $30\%$, and so would have to pay tax on that extra $\$3$, leaving us with about an extra $\$2.10$ after tax. So in this situation, $2+3$ ``$=$'' $4.10$,\footnote{After 6 lectures on Ethics in Mathematics, the students still got tricked.} thus demonstrating that the absoluteness of abstract mathematics quickly degrades in the hands of humans.

We used this maxim, and toy example, to encapsulate what we had been professing throughout the entire course: that mathematical work involves vast human input and subjective judgement. The absoluteness of whatever truth it contains is lost when it is taken \textit{by} humans to create insight \textit{for} humans.

\subsubsection{Objectivity and purity}

Having discussed the unavoidable subjectivity in certain uses of mathematics, we moved on to look at whether \textit{any} part of mathematics is objective or `pure'. This so-called `purity of mathematics' which, though seldom expressed explicitly in the curriculum, comes up frequently in an implicit way throughout the teaching of mathematics and the community of mathematicians (Section \ref{s:hand-down}), as well as in the way (prominent) mathematicians speak about mathematics (Section \ref{s:handing}). And such a viewpoint can have serious consequences for how mathematical work is carried out, as covered in Section \ref{s:impartiality-myth}.  But, turning back to the earlier discourse on truth vs meaning, we explained to students that while there might be objective truths in abstract mathematics, there is no objectivity in the pursuit of mathematics. Indeed, the very act of choosing which mathematical truths to seek out and work on is a highly subjective choice by mathematicians, taking into account surrounding factors such as personal tastes and interests (what they like to work on), sources of funding (what they will be paid to work on), and aspects of the wider community interest (what they will be congratulated for working on). None of this is objective or pure, and this idea is extensively examined by Ernest \cite{ernest2021mathematics}, commenting on how mathematics can (among other things) even be viewed as a ``democratic'' process due to proofs needing to be ``open to scrutiny by others''  \cite[p.~3143]{ernest2021mathematics}, and observing ``the irony that adopting the position that mathematics is largely value free is indeed itself a values position \ldots leading to the `choosing not to choose' fallacy''  \cite[p.~3144]{ernest2021mathematics}. And he explicitly argues that mathematics is value-laden, and not objective or value-free \cite{Ernest2016}. Indeed, even an `objectively-true mathematical statement' can be false; one need only look at the  myriad of mathematics papers in the literature that have been withdrawn or corrected, and the fact that most mathematical `proofs' are not done with predicate logic or other formal tools such as \textit{Lean} \cite{leanlang}. The way mathematicians (i.e., humans) \textit{choose} to verify mathematical statements affects what they \textit{consider} to be true. What \textit{is} true does not change, but what matters is only what humans consider to be true.

Indeed, we pointed out that the very term `pure mathematics' is not an ideological one, but more a descriptor of the motives for doing mathematics (for its own sake), rather than `applied mathematics' which is mathematics which is done to solve some sort of problem (see Wigner's distinction between mathematics and physics \cite{Wigner1960} as outlined in Section \ref{s:whymath}). This is in contrast to the common institutional partitioning of pure/applied mathematics by topic alone \cite[p.~41]{muller2022hippocratic}, and in this sense, ``number theory can be used in applications like cryptography while methods of calculation in numerical analysis can be studied for their own sake without regard to any uses they may have'' \cite[p.~41]{muller2022hippocratic}.

But although this objectivity is somewhat of a myth, can mathematicians perhaps still do \textit{some} things in an objective way? Can the truth of a mathematical statement be something that can stand up to absolute scrutiny, and moreover if mathematicians see a verifiable error in a piece of mathematical work that they can point out in a clear, objective manner, is that not an objective process? Well, we put it to students that this may not always be the case. We asked them to consider what they would do if they found an error in a piece of mathematics somewhere, perhaps in a research paper or text book? We suggested they may choose to inform the author(s), or communicate that error to the community, or write/publish something pointing out the error and a mathematical proof of why it is wrong and/or how to correct it. But we then asked them to reconsider their choice of action in the case that the error was in a piece of writing by a senior member of the mathematical community; perhaps someone who holds a prestigious award such as a Fields Medal, or a high-ranking position such as a head of a mathematics department. We suggested that this may affect their `objective' analysis of the error that they found, and queried whether they would  still be willing to point out the error to the author, or the wider community, or write/publish about it, especially if they were a young research mathematician seeking an (academic) career in mathematics. We suggested that, perhaps, the hierarchy of the mathematical community might instil a combined sense of fear (i.e., ``What if I never get a job because of this?'') and self-doubt (i.e., ``I'm probably wrong, as they're a Fields Medallist'') in them. We asked them to think very carefully about whether they would be brave enough to (openly) go up against a Fields Medallist, or chief scientific advisor of the government, or head of a department? Or whether their self-preservation instinct might kick in, and make them behave in a way that is aligned with the mathematics community and the hierarchy that exist within it.

We also raised the question of whether  mathematicians always have the unrestricted right to do as they please with their mathematical output, even if their piece of mathematics is somehow (objectively) shown to be true. We pointed out the existence of various notions of \textit{responsible disclosure}, which may be quite relevant to mathematicians working in areas such as number theory (and associated cryptographic implications such as the  breaking of critical cryptoschemes), and AI (and the consequences of developing AI tools that are ultra-powerful, and making these widely available). We suggested that just because the mathematicians presenting the mathematics are not doing anything harmful with the work \textit{themselves}, does not absolve them of the responsibility of considering in whose hands they place such capabilities and what \textit{they} might do with it, as already covered in Section \ref{s:designinig}. Mathematical work, like many other technological developments, carries with it the potential for dual use (Section \ref{s:dual}) beyond what is envisaged by those developing it; an ideal already well-studied in many other technical disciplines (physics, computer science, etc.), as well as considered more broadly by philosophers and ethicists.

\subsubsection{Neutrality -- One cannot hide in academia}

We pointed out that, in our experience, mathematicians often like to revel not just in the `objectivity' of abstract mathematics itself, but also in the `objectivity' of their own community. We have heard mathematicians assert that they are able to judge mathematical work, and judge mathematicians, in a totally objective way, given the precision and absoluteness of (abstract) mathematics. Such arguments might seek to exploit the fact that `mathematics is either right or wrong', and so mathematicians can simply `add up' all the mathematics someone has done to evaluate their mathematical worth. Of course, any such argument would need to deal with the fact that, as previously discussed, verifying mathematical proofs is not a flawless process, and moreover some sort of (human-chosen) metric would be needed to `add up' how much mathematics someone has done; akin to the `impartial plagiarism testing algorithm' from Section \ref{s:impartial-algo}. But even if overcome, one then needs to deal with aspects of evaluation that touch on `what is good mathematics?', including attributes such as `beauty'; concepts that mathematicians are far from having an agreed-upon (let alone objective) measure for.

And so, we put it to students that any such perception of objectivity in judging mathematics and mathematicians is simply an illusion, and that as such, mathematicians make subjective decisions on a daily basis, in all aspects of their work. In particular, they judge other mathematicians, and their \textit{mathematical} work, in a subjective way all the time\footnote{To argue that ``this is not part of doing (technical) mathematics itself, and therefore not part of ethics in mathematics'' would be akin to arguing that issues around billable hours would not be part of ethics in law; that may not be part of technical legal arguments, but certainly forms an integral part of `being a lawyer' and `doing law'.}. They (subjectively) judge them when awarding part-marks for exam solutions; for who should decide what makes up part of a solution, and what each part is worth?  They judge them when they write a letter of reference; for if evaluating a mathematician is openly objective, what is the point of a referee? They judge them when they review, or are an editor for, their research paper and its perceived relevance or depth; for beyond correctness of a result, who decides whether it is worthy of publication?  They judge them when they sit on a selection committee, promotions committee, award committee, or employment committee; for if it's just their mathematical work being weighed up, why does this need a whole committee, as surely one person would suffice? All of these processes involve the personal, subjective judgements of mathematicians, which go beyond strict mathematical proof and accuracy. We used all this to demonstrate to students that mathematicians have formed a community, and expect other mathematicians to abide and be judged by the human-made rules of that community.

At this point we broadened the argument even further, highlighting that it is not just  a mathematician's own actions in isolation that matter, but also the actions of those they choose to associate themselves with in various ways, even in academia. In terms of the product of their work, mathematicians might consider where, and what, their students go on to do with their mathematical knowledge, and how they as an instructor affect that. This includes not just what they teach and how (see Section \ref{s:teaching-motives} on teaching mathematics), but also how they do or don't provide support like suggestions on where to work, or letters of reference. We asserted that it is a deep ethical choice by each mathematician to decide if and when to provide a letter of reference. Such an act is a choice, not an obligation, on a teacher, lecturer, or instructor, and while there are some who provide reasons to write such letters, there are also others who provide reasons \textit{not to} \cite{jmc_letter_writing}, and each mathematician needs to weigh up when, and how, to do so\footnote{The particular example we gave for students to consider was whether they would provide a letter of reference for a student seeking to work in nuclear research for a country currently engaged in military conflict with their (the letter writer's own) nation.}. Here, claiming to `be impartial' by always agreeing to provide a letter simply leads to the ``choosing not to choose fallacy'' yet again  \cite[p.~3144]{ernest2021mathematics},  so this this is still an \textit{ethical} choice by the writer.

And furthermore, in the \textit{source} of their work, we suggested that mathematicians might wish to consider not just the work that they choose to do, and why they choose to do it, but also \textit{who} is funding their work, and why \textit{they} (the funder) have chosen to fund it \cite[p.~16]{chiodo2023manifesto}. Money, like mathematics, is not objectively clean; it's source matters. Indeed, one can consider the work by Joi Ito in his time as director of the M.I.T.~Media Lab, and how he chose to accept funding from Jeffrey Epstein who was, even at the time, a convicted sex offender \cite{farrow2019elite}. Ito accepted USD$\$525,000$ of donations from Epstein to the Media Lab, and another USD\$$1.2$million ``for investment funds under [Ito's] control'' \cite{farrow2019elite}. Epstein also assisted in soliciting an additional USD$\$7.5$million for the Media Lab from other prominent donors such as Bill Gates and Leon Black \cite{farrow2019elite}. Furthermore, it was reported that Epstein made use of these institutional relationships for his own gain, and in particular ``Epstein used the status and prestige afforded him by his relationships with élite institutions to shield himself from accountability and continue his alleged predation'' \cite{farrow2019elite}. And so, while no questions about the validity or safety of the academic work at M.I.T.~funded by Epstein were raised, this was nonetheless a scandalous affair, and all but destroyed Ito's (otherwise rather stellar) academic career\footnote{If you roll with pigs, don't be surprised when you get mud on you.}; he resigned from the Media Lab the day after full details regarding his dealings with Epstein were revealed publicly \cite{farrow2019elite}.

\begingroup
\let\section\subsection
\renewcommand{\bibname}{Bibliography (Lecture 6)}
\renewcommand{\refname}{Bibliography (Lecture 6)}
\addcontentsline{toc}{subsection}{Bibliography (Lecture 6)}
\renewcommand{\addcontentsline}[3]{}
\bibliographystyle{IEEEtran}
\bibliography{bibLec6}
\endgroup

\newpage
\section{Lecture 7: Psychology 101: How to survive as a mathematician at work}

\fbox{
  \parbox{0.98\textwidth}{
Mathematicians eventually form part of the workforce, facing the same workplace interactions, issues, conflicts and pressures that arise in all other professions, and needing to know how to deal with these. Yet their focused and dedicated nature, enhanced by the training they receive, might make it easy to overlook instances of others trying to exploit or manipulate them into producing harmful work.

This lecture looks at how mathematicians can be susceptible to factors such as false choice, workplace pressure, intentional manipulation, and social conformity. It shows how they can start to identify such situations, and provides insight and empowerment to protect themselves, and their work, against these harmful influences.
}}

\subsection{Choice}\label{s:choice}

This was the final of the theory-based lectures, looking at principles from psychology and at how mathematicians can be `led astray' in the process of carrying out their work, and how even those with ethical intentions might struggle to enact them. That is, we explored how mathematicians might find themselves in circumstances where they lack (sufficient) autonomy to \textit{be} ethical. The lecture broadly aligns with Raz's interpretation of autonomy \cite{raz1988autonomy}, where he argues that for an agent to be autonomous, it must meet three conditions \cite[p.~369]{raz1988autonomy}:

\begin{enumerate}
    \item It must possess certain mental capacities.
    \item It must have an adequate range of valuable options.
    \item It must enjoy independence from coercion and manipulation.
\end{enumerate}
The preceding lectures extensively covered the ``mental capacities'' needed by a mathematician to act ethically. So we then moved on in this lecture to cover the ``range of valuable options'' (i.e., choice) and ``independence from coercion and manipulation'' (i.e., avoiding pressure, manipulation, and conformity) that mathematicians require to act as autonomous ethical agents. With this final lecture, an ethically-minded mathematician can start to learn some of the tools they will need to maintain their ethical autonomy, in professional and workplace settings that might otherwise erode that autonomy.

\subsubsection{Activity: Determining choice}

We started the lecture with a simple activity, surrounding the concept of choice. We invited three volunteers from the audience to come up to the front and partake in the following simple exercise\footnote{We would always aim to choose students who had attended in a previous year, so that the audience would be almost entirely attendees seeing this for the first time, and moreover so that the volunteers would fully understand what was coming.}:

Each person ($A,B,C$) was asked to choose \textit{one} pen, in different ways:

\begin{enumerate}
    \item $A$ was offered a selection of \textit{one} black pen.
  \item $B$ was offered a selection of \textit{two} black pens.
  \item $C$ was offered a selection of \textit{one} blue pen and \textit{one} red pen.
\end{enumerate}
None of $A,B,C$ could speak; they simply chose one pen each. We then invited a rather open discussion among the remaining attendees (not the three volunteers), to ascertain which of $A,B,C$ had `choice'. This was done in sequence, with a `vote' held at the end of each, asking whether that volunteer ``had choice'' or not.

We first let the students discuss whether $A$ had `choice'. After some deliberation, often with a few alternate proposals such as ``$A$ could have chosen not to take a pen, or even not to come up to the front'' (which we acknowledged as valid proposals), we took a vote, and the students usually (mostly) agreed that $A$ had ``no choice''.

We then let students discuss whether $B$ had `choice'. Here, there was usually heavy debate; students would quickly converge to the notion that ``$B$ had \textit{more} choice than $A$''\footnote{Which we always found odd, seeing as it was clear students were struggling with the very notion of choice, yet somehow had ways to quantify it.}, but we stressed that the question was simply ``Did $B$ have choice?''. Students seemed to get stuck on the fact that there were two pens that did the same thing. The subsequent vote seldom gave clear consensus, across all the years this was done; to a first approximation, it was 50:50 each time.

And finally, we let students discuss whether $C$ had `choice'. Here, students quickly agreed that yes, $C$ had \textit{some} choice, but often brought up additional arguments, usually along the lines of ``What if $C$ didn't actually want either colour, and preferred green? Or wanted chalk instead of a pen?'', concluding that while $C$ certainly had some choice, they could have had more. And the vote usually reflected this, with the overwhelming majority voting that yes, $C$ ``had choice''.

We then gave the students a final opportunity to change any of their votes; few ever did. We then invited each volunteer to write their initials in very large text on a blank A4 page, with the pen they had chosen, and hold it up to the room; all three pairs of initials were written in black. This is because, before the lecture, we had swapped out the blue and red ink cartridges for black cartridges, thus although the red and blue pens \textit{looked} like they would write in red and blue, they both wrote in black, just like the black-looking pens\footnote{After 7 lectures on ethics in mathematics, the students were still surprised that they were able to be tricked so easily.}.

While the attendees would, in their `deliberations', debate whether $A,B,C$ have choice, and say sensible things like ``Perhaps one of the black pens has slightly more ink than the other?'' or ``Maybe the volunteers don't even want a pen?'',  each time we ran this,  none have every asked about, or tested, the writing colour of the pens. They all \textit{assumed} the external colour of the pen was the writing colour. This connects back to the discussion in Section \ref{s:subservience} about mathematicians being subservient and never questioning aspects of the problem put before them (already \textit{implicitly} witnessed in the oil pipe problem of Section \ref{s:oil-pipe}), and seeing this play out again so \textit{explicitly} in real life becomes a useful learning experience for the students in helping them see the value of `questioning the question'.

\subsubsection{Trying to understand choice}

In order to break down the previous activity into more tangible concepts, we introduced the following 4 terms to students, to help them see what sort of choice volunteers $A,B,C$ did, and didn't, have. Each reflects a type of choice which, while non-exhaustive, still illustrates that there are indeed different types.

\begin{enumerate}
    \item \textit{No choice}.
    \\Here, there is only one option available.
    \item \textit{Numerical choice}.
    \\Here, multiple options are available. However, there may not be anything that sets any of them apart from the others.
    \item \textit{Distinguishable choice}.
    \\Here, multiple options are available, and not all identical. However, the differences between them might be minor, and unrelated to the high-level purpose of the selection.
    \item \textit{Meaningful choice}.
    \\Here, multiple, distinguishable options are available, some of which lead to genuinely different outcomes in the context of the reason for making the choice in the first place.
\end{enumerate}
This connects back to Raz's \cite{raz1988autonomy} notion of having an ``adequate range of valuable options'' in order to have autonomy. And so, in our activity, $A$ had no choice. $B$ had numerical choice (2 pens), but they were identical, so no distinguishable choice. $C$ had distinguishable choice (a pen that appeared red, and one that appeared blue), but they performed exactly the same when used to write (both had black cartridges), so no meaningful choice.

Of course, as some students rightly pointed out, they might choose to \textit{engineer} additional options, with the aim of expanding how much choice they have, by refusing any of the options offered and going and seeking out another option themselves. This sort of \textit{ethical} thinking is indeed quite valuable, as it allows an exploration of additional, better options to carry out, rather than be subservient to whoever is providing the (initial) options.

\subsubsection{Faux choice}

We then went on to a lengthy exposition about the activity that had just been run. First, we emphasised that it was essential to have human volunteers at the front, to ensure that the ensuing debate made sense and was not stuck in abstraction. We highlighted that choice is a very human concept, stemming from the idea of \textit{human agency} (i.e., a human's ``ability to actively impact the system'' \cite[p.~4]{chiodo2025formalising}), and has an important role in all human activity, including tech \textit{development} surrounding what developers choose to do \cite{chiodo_educating_2024}, and tech \textit{safety mechanisms} surrounding what benefit human oversight can effectively bring to a (semi)-automated process \cite[p.~4]{chiodo2025formalising}. We pointed out that this idea of agency had already become obvious from earlier when students were asking things such as ``What if $C$ wanted chalk instead of a pen?''; there, they were (perhaps subconsciously) trying to increase human agency within the choice process.

We then moved on to explore how such notions of \textit{faux choice} might be presented to a mathematician in their working life, and how it might influence their decision making process, and their moral judgement of a situation. For example, a company may give its employees the `choice' of working on various projects, and they are allowed to reject some (for reasons such as having a moral objection) in favour of others; selecting a `good' project over a `bad' one. But the tools that mathematicians create in their `good' project might easily be applied directly to the project(s) they rejected; a very real case of dual use, as covered in Section \ref{s:dual}. Explicit examples of this can be seen in various areas already discussed, such as in Lecture \ref{s:cryptoprivacy} on surveillance. For example, suppose a mathematician is working in the surveillance industry, and is asked to create surveillance tools for collecting and decrypting email communications from a particular group of people; perhaps journalists. And suppose that the mathematician objects to this project, and is instead allowed to work on another project, say, collecting and decrypting email communications from another group of people; perhaps suspected tax evaders, whom they are happy to surveil. This may \textit{seem} like a genuine choice; they have a moral objection to surveilling journalists, but no objection surveilling suspected tax evaders. But if the mathematician creates mathematical tools for the second project, what might happen? If the organisation they work for is in any way dishonest to them -- and there is little to stop an organisation lying to its employees -- then the tools they create to surveil tax evaders may very easily be `exported' back in to the  project surveilling journalists, nullifying their original `choice' not to partake in that project. We emphasised that mathematical and technical tools are difficult to create, but once done, they are very easy to copy and reproduce elsewhere; \textit{dual use} is a an issue that mathematicians need to be acutely aware of at all times.
  
Connecting this back to the activity, we pointed out that a mathematician being presented with just one option -- carry out surveillance of journalists -- might clearly see that they are in a situation akin to volunteer $A$ being presented one black pen, with no choice available to them. However, being presented with an option to surveil either journalists or tax evaders, they may think they have meaningful choice, but actually be in a situation akin to volunteer $C$; having distinguishable choice, but (unbeknownst to them) no meaningful choice if their work will simply be brought back to the journalist-surveillance project they specifically did not want to participate in. Such an instance of \textit{faux choice}, whereby the chooser perceives they have a much higher level of choice than is actually present, could easily occur in any area of mathematics, given how easily mathematical tools and output can be brought from one specific application to another, even without the original creator knowing.

The underlying point of this section was to illustrate to students that, not only do they have limited control over how their mathematical tools and outputs are used, but they might also be (blissfully) unaware that their work is being put to (other/mis)-use. Moreover, they may think that they have made an ethical choice on what work to do or avoid, but through a process of faux choice, still end up contributing to the projects they specifically `chose' not to. Just like false privacy (Section \ref{s:false}), faux choice can be equally disempowering, for all the same reasons.

\subsection{Pressure to do your job}

Having introduced the notion that the behaviour of other people might influence how and what work mathematicians do, we then extended the concept from manipulating \textit{choice}, to exerting influence through direct requests (with consequences). We started with a fairly obvious observation, given all the examples presented so far: that  mathematicians create and wield tools of vast power. By virtue of modern technology, these tools are much more powerful than those that mathematicians have been able to create in the past. But ultimately, they give these tools to \textit{others} to decide how to use, implement, and derive benefit from. Mathematicians work for people, for organisations, for institutions, for management;  at some point, others get to control the impact and outcome of what they create. And so, mathematicians need to be aware that these other groups may not have the same interests and objectives they do. Mathematicians might need to  stop and think about the work they do, because of the sheer power of the things they create. They are creating intellectual tools; one could even argue they are creating intellectual weapons, so-called ``weapons of math destruction'' \cite{o2016weapons}. As mathematicians, they cannot necessarily control how their tools or outputs will be `put into operation' (i.e., ``integrated'' \cite{muller2025integrators}), and whether this will match their own goals. And with any sort of weaponry, one needs to think carefully about how it's being used, how it's being applied, who's using it, and how they're applying it \cite{chiodo_educating_2024}.

\subsubsection{The threat of dismissal}
 
We informed students that, in the workforce, the day may come where they might be asked to do something that they fundamentally disagree with, and do not want to do. And this is not hypothetical, as there are already AI developers such as Jacob Coxon -- a mathematics graduate from Cambridge \cite{booth2026coxon} -- who have quite their roles due to concerns about AI's threat to humanity and what their employer (here, Anthropic) is doing \cite{zeff2026coxon}. So we began to explore what students might do if they eventually found themselves in such a situation. We suggested that first, they might want to express such an opinion to their employer or to their peers; this might be successful in changing the request. However, several other things can happen if they object. They might be told ``This is your job, you should just sit down and do it'',  and perhaps subsequently think to themselves ``OK, I'll just do what I'm told''; an entirely reasonable decision to make as an individual. But, they may also choose to continue their objection, and say ``No, just because it's my job, doesn't mean I need to do it'', and so  voice more objection. At this point, their  employer might threaten them with termination and say ``If you don't do it, then you'll be fired''. For some people, this may indeed be a problem, and such a threat would suffice to `knock them back into line', especially if they have significant regular financial commitments such as a mortgage, a family, etc. But, others might think ``I'm fairly employable as a mathematician, and can find a job elsewhere. So I'm happy to leave this job if these are the conditions that are placed on me''. By absorbing some of the consequences, a mathematician -- indeed, any worker -- can still make a meaningful difference to the outcome of any undesirable task.

\subsubsection{What's your price?}

We then pointed out that the threat of `less money' (i.e., being fired) is not the only tool at the disposal of an employer; they can also offer `more  money', by giving the employee some sort of financial incentive. Mathematically, this is the same sliding scale, or axis, of incentivisation: move the money needle. Such an incentive might come in the form of a one-time bonus, pay rise, promotion, or anything else that `adds money' to the employee. Just like in the previous example, where an employee would need to decide between turning down the difference between their current salary and zero salary (i.e., be fired), this example involves them turning down the difference between a new higher salary, and their current one; it still costs them money. And while getting fired has a bounded cost (their current salary), the amount of extra salary they could receive is essentially unbounded.

It was at this point that we informed students of the maxim that (almost) ``everyone has their price'' \cite[p.~459, Item 26]{Kant1996}, and that vast sums of money might become extremely convincing, even in the face of their strongest moral objections. Indeed, this idea goes back at least as far as Kant in 1793, which he in turn draws from the bible \cite[p.~85]{Kant1793}\footnote{Where he writes \textit{A member of the English Parliament exclaimed in the heat of debate: ``Every man has his price, for which he sells himself.'' If this is true (and everyone can decide by himself), if nowhere is a virtue which no level of temptation can overthrow, if whether the good or evil spirit wins us over only depends on which bids the most and affords the promptest pay-off, then, what the Apostle says might indeed hold true of human beings universally, ``There is no distinction here, they are all under sin - there is
none righteous (in the spirit of the law), no, not one.''} \cite[p.~85]{Kant1793}. The first quote is due to English MP Sir Robert Walpole \cite[p.~459, Item 26]{Kant1996}. The second is derived from the bible, Romans 3:9--10, which contains the line ``There is none righteous, no, not one.'' \cite[p.~459, Item 27]{Kant1996}}. And so, a simple `yes' to a straightforward request might put them into a life of luxury and opulence, and might weigh on their decision-making quite heavily. As an example, we suggested that students go home and consider whether they would carry out a morally-questionable (yet legal) mathematical task such as setting up a financial arbitrage position against investments by pension funds or registered charities\footnote{We do not know for sure if such activities are actually legal, but for the purpose of the thought exercise, we assumed that they were.}, for the following reward amounts:

\begin{itemize}
    \item Their `current salary'.
    \item A \$1000 bonus.
    \item A \$1million bonus.
    \item A \$1billion bonus.
    \item A \$1trillion bonus.
\end{itemize}
We gave two follow-up analyses of this to students. The first was that, regardless of how their morals and ideals stand at present, it is likely that when faced with a real scenario like this, they will indeed each have a `price' at which they would yield and carry out the task. Such a number might be quite high (for some), but it almost certainly exists. Regardless of what they say in a hypothetical situation, things are very different when there is actual money on the table. And all the amounts above are perfectly plausible; even a \$1trillion bonus is now in the realm of possibility, as Elon Musk was offered such a bonus by Tesla if he met certain performance targets \cite{bbc2026musk}.

And second, we suggested to students that, while the above sums may all seem like unimaginable amounts of money to them at present given their current finances, in the future this may not be the case. Later in their careers, they may find that \$1000, or even \$1million, is not a significant amount of money for them, and it is not impossible that \$1billion becomes insignificant for some. So indeed, what may seem like a tempting sum of money \textit{now}, may actually not be so tempting later in life. Thus, their perspective on what might incentivise them can change over time, in both directions; that which they refuse to partake in now might prove tempting later, and that which seems tempting now might not be so later on. And of course, this need not only be about money; we highlighted that there are other work factors employers can change, for better or for worse, to incentivise employees. This can be things like working conditions, location, hours, project assignment, team assignment, and so on. For some, being assigned an exciting (or boring) project might be just as incentivising (or disincentivising) as a raise (or pay cut).

\subsubsection{If you don't do it, someone else \textit{might} (or might not)}\label{s:ifyou}

However, we then stressed that employers are quite clever. They can even go one step further, beyond the axis of `money up or down', and use arguments of the form ``If you don't do it, someone else will''. This is, of course, a very powerful argument psychologically. It's very convincing, as it belittles one's objection and might make them think ``Even if I don't do this, and I have strong convictions about it, if I go to some sort of personal expense to make sure I don't do it, then someone else will come along and do it anyway''. And, when successful, it has a powerful effect as it limits the number and quality of choices available to them, thus eroding ethical autonomy as per Raz (Section \ref{s:choice}). But, while this argument \textit{might} sound convincing in many other areas of work, it doesn't necessarily carry much weight against mathematicians, for two very important reasons.

The first reason is that there simply aren't that many mathematicians. They posses particular training, particular skills, particular ways of thinking, that the bulk of the population does not. So it doesn't take a great deal of objection, it doesn't take a great number of mathematicians to say ``No, I don't think this should be done'', for the job to not be done, or stop happening. In addition, the `handover' process for such technical work can be quite difficult, and indeed almost impossible if no-one from the `old team' is there to explain to the `new team' what the tool or object is, how it works, and how to keep developing it; one need only envisage a codebase with no commenting, to see just how hard this can be. Thus, the moral objections of mathematicians, coupled with them standing their ground on something they disagree with, can have a much bigger impact than in most areas.

The second reason is a more nuanced one: mathematicians understand the contrapositive. The statement ``If you don't do it, someone else will'', is logically equivalent to ``If no-one else does it, you will'', explained in full in \cite[p.~56--57]{chiodo2023teaching}. This is, of course, absurd. The argument has an implicit underlying assumption that the task being requested will definitely be completed; there is nothing to say that is the case. As an extremal example, one can consider the (flipped) argument of an employer asking ``If no-one else is going to build me a nuclear bomb, then you will''; this is neither true, nor terribly convincing. Indeed, hearing that every other human has rejected the request might psychologically \textit{persuade} someone to reject it. So what mathematicians\footnote{Well, everyone, really. But this is not easy to explain to those who have not been trained in mathematics.} should really be considering here is
the argument ``If you don't do it, someone else \textit{might}''. This is true; someone else might. But they might not. 

However, we suggested that this defence goes even deeper, arguing that even \textit{if} `someone else' can, and will, do the task, they may not be easy to find, or may take a long time to source and bring on board, or not have the understanding or skills to be able to carry it out as efficiently or accurately. Any of these can effectively delay the process, perhaps long enough to avert harm. If an activity that is deemed morally objectionable becomes impractical or prohibitively expensive to carry out, or is significantly delayed, then that can almost be as effective as it not being done at all. Complex mathematical work can often be time-sensitive; a (digital or kinetic) weapon of warfare might only have a limited window in which to be used, a digital exploit or `zero-day' may only be usable for a certain period of time, or a particular method for collecting data for AI training may only be legal up to a certain point. And so, an evolution or change in external circumstances (e.g., a target for surveillance is no longer of any value), or internal circumstances (e.g., there is a change in management or political leadership, leading to such surveillance no longer being desirable), might be on the horizon, meaning a \textit{sufficient} delay has the same outcome as a blanket refusal. In such cases, the power of meaningful objection still sits with the mathematician.

And so, we highlighted to students that, whether they choose to take the pragmatic perspective that there are not many mathematicians, or the logical perspective of the contrapositive, (and factoring in that a delay can be just as good as a prevention), the conclusion remains the same: their objection \textit{means something}, and can have a material impact on the outcome and potential harms that they are worried about. That is, they can retain ethical autonomy, even in these difficult scenarios.

\subsubsection{Avoiding forks}

Having covered many of the difficult workplace situations mathematicians might find themselves in, and potential tools to deal with them \textit{when they happen}, we then moved on to another strategy they might employ: \textit{prevention}.

We observed that many of the above scenarios involve making a (sometimes difficult) decision between several options, where each one leads to some sort of \textit{loss} for the mathematician. They may need to choose between carrying out an unethical task, and keeping their job, or (unfairly) missing a promotion or bonus, or the uncertainty that someone else \textit{might} do it. In such situations, there can often be very little time to make the decision, and no way to avoid a loss of some form. This is called a \textit{fork} in chess; when a player gets into a position where 2 (or more) of their pieces are simultaneously threatened, yet they can only save at most 1. A fork is problematic in that it is a  particular way to limit someone's `range of options', and therefore reduce their (ethical) autonomy as per Raz (Section \ref{s:choice}), which is why mathematicians need to learn about this.

It may seem that such situations will always end in an fork, with a mathematician needing to effectively choose between their morals, and other aspects of their life. But this need not be the case. It is indeed impossible to get \textit{out} of a fork, by definition. So one must therefore avoid getting \textit{in} to a fork in the first place. And we pointed out that there are certainly ways to do so. Our primary suggestion to students was to collectivise and collaborate in the workplace \textit{beforehand}, so that when such a situation arises they are ready not just as individuals, but also as a group. For example, having a collective notion within a team or group of what is acceptable and unacceptable to work on or develop, helps counter the risk that `someone else might'. Or a collective agreement that if one person objects, then the entire team will respect that and not work on it either. In addition, there is the option of forming or joining a workers union; tech workers are increasingly unionising \cite{pardes2026unionize}, and there are now branches of existing unions that specifically focus on tech workers \cite{unitedigitaltech2026}, and even entire unions dedicated to them \cite{ifpte2026}. This helps protect against dismissal or other unfair workplace treatment. And we also suggested to students that they can take other measures, such as have sufficient savings banked, and a preparedness to change jobs, should the choose to leave, or be fired from, their current one in the process of standing their (moral) ground.

We emphasised that such approaches require forward-thinking and actions in advance of any such situation, and in particular a specific kind of forward-thinking: the ability to think in counterfactuals \cite{sep-counterfactuals} and reason about hypotheticals that are uncertain but plausible. And all this must be done in advance. Leaving it to collectivise when things start going wrong is unlikely to achieve much. Likewise, one needs to save money and prepare for potential career moves long before unemployment hits. But we explained to students that such forward-thinking was a necessity if they want to retain ethical autonomy in difficult workplace situations and avoid a fork between that and their finances or having others carry out the task anyway. We suggested that perhaps this might not be intuitive for mathematicians, who may think that the problem can be dealt with and solved `when it arises', just like a mathematical problem. So we stressed that no, this is different, and needs to be thought about and acted on long before the day comes.

\subsection{Manipulation}

\subsubsection{Example: The Milgram experiment and the power of authority}

It is not necessarily the case that it takes a bad person to do a bad thing. The power of coercion and authority might weight heavily on their decision-making, more so than they perhaps realise possible, and thus erode their ethical autonomy as per Raz (Section \ref{s:choice}). Indeed, this very idea was studied in the aftermath of WWII and the Holocaust, with the slaughter of millions in death camps. The psychologist Stanley Milgram designed the \textit{Milgram Experiment} in 1963 to test this hypothesis regarding obedience \cite{milgram1963behavioral}, noting that ``obedience may be a deeply ingrained behavior tendency, indeed, a prepotent impulse overriding training in ethics, sympathy, and moral conduct'' \cite[p.~371]{milgram1963behavioral}. 

Milgram collected 40 male volunteers, aged 20--50, from a range of backgrounds, and ran the experiment at Yale University \cite[p.~372]{milgram1963behavioral}. Three individuals took part in each session of the experiment \cite[p.~373]{milgram1963behavioral}:

\begin{enumerate}
 \item The \textit{experimenter}; a hired participant who pretended to be in charge of the experiment.
 \item The \textit{learner}; another hired participant who pretended to be a volunteer.
\item The \textit{teacher}; a volunteer for a single session, who was led to believe that they were merely assisting, whereas they were actually the subject of the experiment.
\end{enumerate}
The experimenter explained that the purpose of the study was to understand ``what effect punishment will have on learning'', and the teacher then saw the learner being brought into an adjacent room and  ``strapped into an ``electric chair'' apparatus'', with ``an
electrode \ldots attached to the learner's wrist'' \cite[p.~373]{milgram1963behavioral}.  The teacher was then brought back into the first room, and sat in front of a ``shock generator''; a large instrument panel with 30 switches labelled from 15 to 450 volts in 15 volt increments. The switches were additionally labelled in clusters of 4 with a shock severity, ranging across \textit{Slight, Moderate, Strong, Very Strong, Intense, Extreme Intensity,} and \textit{Danger:
Severe}; the final two, 435 and 450 volts, were simply labelled \textit{XXX} \cite[p.~373]{milgram1963behavioral}. The teacher was given a sample shock of 45 volts to affirm that the device worked, and told ``Although the shocks can be extremely painful, they cause no permanent tissue damage'' \cite[p.~373]{milgram1963behavioral}.

The teacher was then told how the task worked: they were to ``read a series of word pairs to the learner, and then read the first word of the pair along with four terms. The learner was to indicate which of the four terms had originally been paired with the first word. He communicated his answer by pressing one of four switches in front of him, which in turn lit up one of four numbered quadrants in an answerbox located atop the shock generator'' \cite[p.~373]{milgram1963behavioral}. The teacher was told to ``administer a shock to the learner each time he gives a wrong response'', and to ``move one level higher on the shock generator each time the learner flashes a wrong answer'', adding that they should ``treat the absence of a response as a wrong answer'' \cite[p.~373--374]{milgram1963behavioral}.

The ``test'' then commences, with the learner scripted to give ``approximately three wrong answers to one correct answer'' \cite[p.~374]{milgram1963behavioral}. Early in the test, ``no vocal response or other sign of protest is heard from the learner until Shock Level 300 [volts] is reached'' at which point ``the learner pounds on the wall of the room in which he is bound to the electric chair''. This  can be heard by the teacher, and  ``From this point on, the learner's answers no longer appear on the four-way panel''. Moreover, the  learner's pounding is ``repeated after the 315-volt shock is administered; afterwards he is not heard from'' \cite[p.~374]{milgram1963behavioral}. If at any point the teacher asked the experimenter whether to continue, or expressed a desire to stop, the experimenter would reply with a ``sequence of 4 ``prods'', using as many as necessary to bring the subject into line'' \cite[p.~374]{milgram1963behavioral}:

\begin{enumerate}
\item Please continue.
\item The experiment requires that you continue.
\item It is absolutely essential that you continue.
\item You have no other choice, you \textit{must} go on.
\end{enumerate}
The test ended when either the teacher still requested to stop even after the fourth prod \cite[p.~374]{milgram1963behavioral}, or they administered the final 450 volt shock \cite[p.~376]{milgram1963behavioral}, whichever came first.

The primary measure of the study was to observe the teacher and record  ``the maximum shock he administers before he refuses to go any further'' \cite[p.~374]{milgram1963behavioral}. A description of the study was given to ``14 Yale seniors, all, psychology majors'' to predict what would happen; the ``most pessimistic'' predicted that 3\% would go all the way and administer the 450 volt shock, and the average prediction was 1.2\% \cite[p.~375]{milgram1963behavioral}. Moreover, ``It had been conjectured that persons would not, in general, agree to administer shocks beyond the level designated as Very Strong Shock'' (i.e., not beyond 240 volts) \cite[p.~375]{milgram1963behavioral}.  The study was carried out, and of the 40 volunteers, the results were as follows \cite[p.~375]{milgram1963behavioral}:

\begin{itemize}
    \item All volunteers administered all shocks up to and including 300 volts.
    \item 5 refused to proceed after the 300 volt shock  (the first in which the learner made audible signs of pain).
    \item 4 refused to proceed after the 315 volt shock.
    \item 5 refused to proceed somewhere in the range of 330-375 volt shocks.
    \item The remaining 26 administered all shocks up to 450 volts. 
\end{itemize}
That is, 65\% of the volunteers (teachers) were \textit{obedient} (i.e., ``complies with experimental commands fully''), whereas only 35\% were \textit{defiant} (i.e., ``breaks off the experiment at any point prior to administering the thirtieth [450 volt] shock'' \cite[p.~374]{milgram1963behavioral}. That is not to say the volunteers were indifferent to the (perceived) impact of their actions; far from it. Many showed signs of ``extreme tension'' \cite[p.~375]{milgram1963behavioral}. And ``Although obedient subjects continued to administer shocks, they often did so under extreme stress. Some expressed reluctance to administer shocks beyond the 300-volt level, and displayed fears similar to those who defied the experimenter; yet they obeyed'', and at the end many of these obedient volunteers ``heaved sighs of relief, mopped their brows, rubbed their fingers over their eyes, or nervously fumbled cigarettes. Some shook their heads, apparently in regret'' \cite[p.~376]{milgram1963behavioral}.

The results of the experiment were very surprising, going against all initial predictions and conjectures. A key finding was how the sheer power of authority, with no consequences for defiance, still  managed to trump morals, even among those who knew they were doing something fundamentally wrong. Indeed, Milgram writes:

\begin{quote}
``Subjects have learned from childhood that it is a fundamental breach of moral conduct to hurt another person against his will. Yet, 26 subjects abandon this tenet in following the instructions of an authority who has no special powers to enforce his commands. To disobey would bring no material loss to the subject; no punishment would ensue. It is clear from the remarks and outward behavior of many participants that in punishing the victim they are often acting against their own values. Subjects often expressed deep disapproval of shocking a man in the face of his objections, and others denounced it as stupid and senseless. Yet the majority complied with the experimental commands.'' \cite[p.~376]{milgram1963behavioral}.
\end{quote}
Simply put: the authority of a technician in a lab coat \cite[p.~373]{milgram1963behavioral} can significantly outweigh the good morals of the average person. As such, it can be a very straightforward process to make an otherwise good person do bad things.

\subsubsection{The consequences of obedience}

The key lesson we drew for students from the Milgram experiment was that moral awareness can be eroded simply by having someone act as an \textit{instrument} carrying out the requests of others. In doing so, the `intermediate person' loses perceived agency, acting merely as an efficiency bridge between the person making the requests, and the implementation of said requests. And when people feel they have no agency, they thus feel they have no mechanism to object or intervene. Indeed, Milgram \cite{milgram1974obedience} later wrote of his initial experiment that:

\begin{quote}
``The essence of obedience consists in the fact that a person comes to view himself as the instrument for carrying out another person's wishes, and he therefore no longer sees himself as responsible for his actions. Once this critical shift of viewpoint has occurred in the person, all of the essential features of obedience follow.'' \cite[p.~xii]{milgram1974obedience}
\end{quote}
We suggested to the students that they think very carefully about this statement, and whether it sounded familiar to a phrase they had heard in the previous lecture: ``I was only following orders'', from Lecture \ref{s:law}, and in particular Section \ref{s:Nuremberg} on `superior orders'. We explained that as well as being used as an \textit{external} legal defence against having carried out vile acts (as was done in Nuremberg), it also serves as an \textit{internal} moral gateway for \textit{carrying out} such acts in the first place. That is, telling themselves ``I'm only carrying out the wishes of others'' enables people to actually do those actions. `Superior orders' can act as a defence against the most important critique of human activity; one's own conscience. And when that is bypassed, all manners of harm are possible. Indeed, as mentioned earlier, ``obedience may be a deeply ingrained behavior tendency, indeed, a prepotent impulse overriding training in ethics, sympathy, and moral conduct'' \cite[p.~371]{milgram1963behavioral}. This further shows that ethics is not just something that needs to be \textit{known}; it also needs to be \textit{acted on}, and \textit{actively defended}. Knowing the difference between right and wrong is a necessary step, but it is not sufficient unless one also knows how to act on that knowledge in a meaningful way.

We then connected this specifically to mathematics, by asking students if they had ever heard anyone use the phrase ``I just do the mathematics'' in defence of mathematical work and its harmful consequences. We argued that it might be tempting to view mathematical work as just a conduit or `instrument' within a larger process, rather than the source or impact of any harm. But, as argued earlier, this `instrumentalisation' of mathematical work, and by extension, mathematicians themselves, is exactly the sort of viewpoint that can ultimately remove perceived agency from mathematicians, reducing their ``range of options'' and eroding their ethical autonomy as per Raz (Section \ref{s:choice}). Milgram clearly demonstrated that viewing one's work as simply `carrying out another person's wishes' is the final step in enabling one to carry out harmful acts, even with full knowledge that they are wrong. And then, when mathematicians and their tools are co-opted, they can be used to manipulate and exploit others, as demonstrated in Lectures \ref{s:financialmaths}, \ref{s:cryptoprivacy}, and \ref{s:algorithms}.

\subsubsection{Psychopaths in the workplace}

Having given an explanation as to \textit{how} mathematicians can be influenced, coerced, or manipulated into doing unethical work, even against their own will or better judgement, we then turned to discuss \textit{who} might do this to them. A naive approach might be to think  ``Well, no-one has any reason or desire to be unkind, so this will never happen''. Unfortunately, this is not the case. There are people who, by many measures, are intrinsically `unkind' (or worse); such people are known as \textit{psychopaths}.

As covered in \cite{Dorwart2026psychopath}, a psychopath demonstrates personality traits related to lack of emotion, such as ``lack of nervousness, absence of remorse or shame, inability to love, shallow affective reactions'', as well as appearing outwardly normal via ``lack of delusions and other signs of irrational thought, superficial charm, and good `intelligence' '' \cite[p.~2]{sanz2021prevalence}. Psychopaths ``manifest socially maladaptive and ethically reprehensible behaviors (e.g., failure to follow a life plan, impersonal, trivial, and poorly integrated sex life, deceitfulness and lack of sincerity, pathological egocentrism)'' \cite[p.~2]{sanz2021prevalence}. In other words, they are hardwired to say, and do, whatever it takes to get what they want, without fear or regret. And doing so can indeed help them get what they want, with there being a notion of ``successful psychopathy'' referring to ``those psychopathic personality traits such as lack of fear, high self-confidence or charisma, which can be beneficial in certain contexts'' \cite[p.~2]{sanz2021prevalence}. But one should not conflate `achieving your aims' with `doing good', if those aims are not actually good.

We clarified that psychopathy should not be confused with \textit{Antisocial Personality Disorder} (ASPD), sometimes referred to as \textit{sociopathy}, as they are ``related, but distinct disorders''  \cite[p.~195]{werner2015epidemiology}. ASPD seems to manifest in a more chaotic way, and those with ASPD are ``more likely to be chronically irresponsible, often neglecting responsibilities at work or school'' \cite{Dorwart2026psychopath}, and ``their behaviour tends to be more erratic than psychopaths'' \cite{mcnally2015difference}. In short: those with ASPD are more noticeable by others and predominantly a harm to themselves, whereas those with psychopathy are less noticeable by others, but at the same time predominantly a harm to others; a dangerous combination.

Next, we explored the idea of where one might find psychopaths, and how common they are. Various studies put an estimate of about 4.5\% of the adult population as being psychopaths \cite[p.~1]{sanz2021prevalence}, though this depends on the measure used; the PCL-R ``gold standard'' measure puts it at about 1.2\% \cite[p.~1]{sanz2021prevalence}. In addition, these numbers change when analysing different  professions. Some studies suggest that the professions with the highest levels of psychopaths include company CEOs, lawyers, salespersons, surgeons, and civil servants; in contrast, those with the lowest include  nurses, charity workers, teachers, physicians, and accountants \cite[p.~3]{sanz2021prevalence}. The numbers for corporate bosses and criminals in prisons is staggeringly high -- 21\% -- and moreover such people may be a `` ``successful psychopath'' who may be inclined to unethical or illegal practices [but] has been allowed into the top ranks of companies because of the way firms hire'' \cite{agerholm2016psychopath}. But why is this of any relevance to mathematicians in the average workplace? Because there is another group of people who have a high proportion of psychopaths compared to the general population: \textit{managers}. Indeed, studies suggest that ``people in a managerial position at work had higher levels of psychopathic traits than those who did not occupy such positions'' \cite[p.~3]{sanz2021prevalence}. And so we pointed out to students that, not only might they work \textit{with} psychopaths; they have a reasonable chance of working \textit{under} them. And in places like Silicon Valley, where mathematicians may well find themselves working, there is indeed a high proportion of psychopathic CEOs there \cite{solon2017silicon}, so much so that regular tech workers there have taken to writing blogposts about how to ``detect and avoid'' them \cite{xander2026darktriads}.

We then moved on to consider why -- if psychopaths are so bad, so toxic, and so harmful to those around them -- can't the public just keep away from them, and keep them away from the workplace? The problem here is that psychopaths are very difficult to detect; they have a myriad of behavioural deficits, but also an ``uncanny ability to mask them'' \cite{solon2017silicon}. The traits of psychopaths are well-documented at a high level, but knowing for example that ``absence of remorse or shame'' is a trait doesn't actually directly help someone to detect it, let alone defend against it. Some people have devised their own active ``tests'' for spotting psychopaths \cite{xander2026darktriads}, but there is no set formula as every psychopath will be different. Furthermore, we suggested that perhaps mathematicians might find it particularly difficult to recognise the characteristics of psychopaths, as their training develops their thinking towards intense, focused thought and literal interpretation of things around them (recall Section \ref{s:strengths}; their mathematical strengths become weaknesses elsewhere). This style of thinking can play into the hands of psychopaths perfectly, as it allows them to tire mathematicians out through logical exhaustion. For example, they may ``invalidate your experiences, without any evidence'', thus trapping a mathematician into being ``put on the defense, precisely where the manipulator wants you'' \cite{xander2026darktriads}.

Psychopaths have an amazing ability to manipulate those around them. However, they are usually clever enough to do it in a way which forces guilt and responsibility onto others, and moreover not get noticed or `caught' in the process.  And so we stressed that \textit{everyone} in the workplace needs to be aware that such people exist, and that such people can carry out manipulation  in ways those around them might not expect. Mathematicians in particular need to be acutely aware, as they may be particularly susceptible to psychopaths, and moreover work in areas or domains, and carry out work, that psychopaths may be particularly interested in. We harked back to the earlier case Tom Hayes (Section \ref{s:libor}), who was jailed for rigging LIBOR, defending himself by arguing ``I was part of the dishonest scheme, so obviously I was being dishonest.'' However, Tom was one of the very very few people in finance to actually be prosecuted for rigging LIBOR; most of those around him were not. We suggested that Tom may have been working in the midst of many psychopaths (which, in high-rolling finance, is unsurprising, especially in higher-up leadership positions), who all had the ability to mask their behaviour and protect themselves from blame. Tom had many of the `mathematical strengths' that made him good at technical mathematics, but also may have deprived him of the the tools, insight, or will to protect himself in such a toxic work environment. And his actions ended up causing significant harm to himself, and to wider society.

\subsection{Conforming to our surrounds}

\subsubsection{Forming ideology}\label{s:ideology}
 
We reiterated to students that throughout the previous lectures they had seen numerous examples of mathematicians behaving in ways that raise serious ethical questions. But how did all these things happen? How did mathematicians actually think that it was a good idea to manipulate financial markets (Lecture \ref{s:financialmaths}), conduct intrusive surveillance on billions of people (Lecture \ref{s:cryptoprivacy}), or manipulate voters during elections (Lecture \ref{s:algorithms})? Why was Tom Hayes happy to be ``part of the dishonest scheme'' (Section \ref{s:libor})? Why did any of this seem like a good idea?

We suggested to students that, in any industry, workers `fall in line' shortly after starting, and the thinking and ethos of the organisation they're in quickly becomes \textit{their} thinking. Many companies \emph{deliberately} recruit young employees because they are \textit{ideologically malleable}. The company teaches them what to think and believe, and accustomises them to the `standard practice' of how things are done there. Such new, young recruits are ``more malleable'', and the younger they are, the more ``open-minded and generally accepting they are of a work environment'', being ``unclouded by experience, which brings advantages to an employer -- they’re easier to mould'' \cite{Christian2023companies}. However, this can be detrimental to young graduate workers, who ``can find themselves vulnerable to exploitation where they haven’t acquired the experience to know what’s OK and what’s not'' \cite{Christian2023companies}.

We explained this with the notion of \textit{first-fill ideology}: take an empty bottle, and fill it with lemonade (that's easy). Now try filling it with orange juice, \emph{without} pouring out the lemonade. There is no room for the juice (the new ideology), because the lemonade (the original ideology) is in the way. We suggested that first-fill ideology will always be the most captivating, and thus the most effective. And this extends all the way back to `pre-employment'; where workers come in as interns. Indeed, young employees such as interns ``arrive at your organization without certain workflows or norms already ingrained'' and so have ``fewer habits to unlearn'' \cite{shedek2025whysmart}. This in turn makes them ``more open to learning how things work at your organization'', with employers knowing that ``if you invite successful interns back to the program in following years, you can continue building their institutional knowledge and develop them into perfect-fit full-time candidates'' \cite{shedek2025whysmart}. But we warned students that this goes beyond working style and knowledge, as ``a certain measure of loyalty comes with hiring interns'' \cite{shedek2025whysmart}. Employers know, and make use, of the fact that ``If you're their first major employer, you'll likely become their gold standard for defining what a ``good job'' looks like'', and can exploit it to monopolise their grip on employees as ``This can be tough to compete with once these young professionals officially enter the workforce'' \cite{shedek2025whysmart}. First-fill ideology doesn't just influence how a worker works; it influences how they view work and even how they define what a ``good job'' actually is.

We finished by warning students that such `loyalty' might manifest in other ways: as a form of \textit{intellectual entrapment}. That is, upon joining an organisation or institution, especially from a young age where they are likely to build some form of `loyalty' towards it, workers may then feel compelled to serve that organisation and `do favours' for it, beyond what they might normally consider doing. We suggested that they may end up seeing their membership or selection there as some sort of `debt' that needs to be repaid, including at a personal level of looking up to a mentor or authority figure and not wanting to disappoint them. Such entrapment may not be (completely) intentional, but nonetheless may even compromise their morals, or even their otherwise consistent adherence to the law. We again reminded students of examples such as that of Tom Hayes, who argued ``I was part of the dishonest scheme, so obviously I was being dishonest'' (Section \ref{s:libor}); a clear example of surrounding institutional culture and practices (negatively) influencing how an individual worker acts. And this can even occur in academia, within a student-supervisor interaction.

But where does the human propensity to conform come from, and how prevalent is it? This very question was investigated by Asch \cite{asch1951effects}, as will now be explained.

\subsubsection{Example: The Asch conformity experiment}
 
Conformity to one's surrounds, as described above, need not occur as an intentionally imposed process. Even when others operating around someone are not \textit{deliberately} manipulating their behaviour, be they colleagues, managers, CEOs, etc., such interactions still might be doing so \textit{inadvertently}. The human desire to conform is very powerful, and can make people doubt even the most obvious of their own objective sensory observations, let alone internal subjective beliefs. This was studied by the psychologist Solomon Asch with the \textit{Asch Conformity Experiment} in 1951 \cite{asch1951effects}.
  
Asch ran his experiment as follows. A group of 8 people were sat together in a row;  7 were `rigged' volunteers (actors) and 1 was a `real' volunteer. They were then shown a series of (not difficult) multi-choice perception questions involving lengths of lines,  needing to ``match the length of a given line with one of three unequal lines'' \cite[p.~178]{asch1951effects}. That is, the volunteers were shown a line between 2 and 10 inches long \cite[p.~180]{asch1951effects}, alongside three other lines, one of which was the same length, and the other two differing in length from it by somewhere between 0.5 and 1.75 inches \cite[p.~178]{asch1951effects}. The task was simple: match the original line to the line in the group of three of the same length.

Each volunteer would, in order, say out aloud their answer to the line test currently posed; this was done for 18 different line tests \cite[p.~180]{asch1951effects}. The rigged volunteers had been given instructions beforehand: starting from test 3, they were each told to say the \textit{same} incorrect answer for 12 of the subsequent tests (they all answered correctly on tests 1,2,5,10,11,14) \cite[p.~180]{asch1951effects}. The experiment then recorded the responses of the real volunteer for these 18 tests, and in particular if they `conformed' to the (wrong) unanimous majority around them and gave an incorrect response contrary to what they could clearly see in front of them. This was done for 50 different volunteers, and the number of incorrect responses for each real volunteer was recorded \cite[p.~180]{asch1951effects}. The purpose of the study was to examine the effect on the real volunteers of being in ``a situation in which a group unanimously contradicted the evidence of his senses'' \cite[p.~179]{asch1951effects}.

The results were quite surprising. Only 26\% of real volunteers  gave the correct response every time, whereas the remaining 74\% answered at least one test incorrectly, roughly evenly spread across 1--11 errors overall (out of 12), and over all the `rigged questions' the error rate was about 32\%, with some real volunteers getting 92\% (11) incorrect \cite[p.~181]{asch1951effects}. 30\% of the real volunteers got at least half of the rigged questions wrong. As a control, an additional 37 real volunteers were given the same experiment, but all responses were written rather than verbal, preventing interactions between respondents; the error rate there was under 0.7\% \cite[p.~181]{asch1951effects}. During the experiment, even those real volunteers with a 0\% error rate would often vocalise justifications for why they were disagreeing with the majority (unprompted) \cite[p.~181]{asch1951effects}, demonstrating that going against a unanimous majority was clearly weighing on their mind, even if it did not affect their final answers. Asch concluded that ``The quantitative results are clear and unambiguous. There was a marked movement toward the majority'' \cite[p.~180--181]{asch1951effects}. That is, going up against a \textit{unanimous majority} has a significant effect on human decision making, even in the absence of direct coercion.

The entire experiment was re-run with a few variants. Firstly, it was run with a ``True Partner''; a rigged volunteer who always answered correctly. In that case, the error rate dropped to 5.5\% (down from 32\%) \cite[p.~185--186]{asch1951effects}. This was modified again to have the true partner revert to the majority (wrong) answers half-way through the test. When this happened, the error rate shot back up to 28.5\% after the true partner `flipped' \cite[p.~186]{asch1951effects}. Another modification had a ``late-arriving true partner'', who only started answering correctly half-way through the tests. In this case, the error rate reduced to 8.7\%; not as low as having a true partner from the onset \cite[p.~186]{asch1951effects}, showing that conformity had some sort of lasting persistence. The size of the opposing majority was also investigated. Going against a majority of 1 or 2 had a negligible or minor effect on conformity, with error rates of 2.8\% and 12.8\%, and error spreads of 0--1 and 0--5 (out of 12), respectively \cite[p.~188]{asch1951effects}. It was an opposing majority of 3 which `saturated' conformity, with an error rate of 33.3\% and spread of 1--12 errors \cite[p.~188]{asch1951effects}. This was run again with opposing majorities of 4 and 15 which, like the original 7, all gave similar error rates of around 32\%   \cite[p.~188]{asch1951effects}. Thus, conformity reached its full effect with an opposing majority of just 3.

\subsubsection{The effects of conformity}

We moved on to discuss how the human desire to conform was  extremely pronounced in Asch's experiments, showing that even a rather objective question with no prior knowledge needed -- `which line is the same length as this one?' -- can be distorted in someone's mind if they feel that answering as they saw it would put them at odds with the rest of the group. This happened even in the face of objective, indisputable truth; the sorts of ideas and objects  that form the centre of the work of mathematicians. Indeed, we suggested that this phenomenon of conformity can distort someone's notion of truth itself; \textit{If everyone else thinks X, then I must surely be wrong}. And thus, while the \textit{hope} of having new recruits conform to a workplace culture and ideology, and the \textit{intent} of hiring those who may be particularly susceptible to  conforming (interns or new graduates), might both be deliberate (as covered in Section \ref{s:ideology}), the \textit{process} by which this happens may well be down to human nature.
  
We highlighted in particular how several of the aspects of forming ideology in the workplace become more understandable in light of the Asch experiment. To begin, conformity to the unanimous majority of a workplace culture is quite clearly demonstrated by Asch. And while having a ``true partner'' can help break conformity, this is offset by the fact that adding new workers to a workplace `one at a time' is reflected in the ``late-arriving true partner'' instance, whereby there is persistence in the conformity. That is, each new worker might have a period in which they `align' to the workplace and become part of the unanimous majority, before the next new worker is then brought in for the process to then repeat. This persistence is further demonstrated by the ``'Stand up, sit down'' test \cite{integration2018conformity}, where an optometrist's waiting room was full of (rigged) volunteers who would all stand up, and the immediately sit down, at a beep that was rung every few minutes. As real patients entered, they would conform to this ``stand up, sit down'' response, having never been given a reason for it. But as the original (rigged) volunteers were slowly removed (eventually going down to zero), the real patients continued to ``stand up, sit down'', and further indoctrinate additional real patients as they entered. This can all be watched on a recent YouTube video \cite{integration2018conformity}.

We stressed that of particular relevance to mathematicians is firstly that such conformity can manifest even through highly objective activities (such as mathematical work), and secondly that such conformity can occur even with universal majorities of just 3; almost the perfect size for a 4-5 person team working on a technical (sub)-project within a company, or even on a research project in academia. So even if a mathematician has minimal exposure to the rest of the organisation, joining a small team of 4 is enough for  conformity to occur. And so we argued that, far from being `immune' to workplace conformity because of the nature of their mathematical work, mathematicians are actually specifically and demonstrably susceptible to it; what they work on, and how they work together, matches the ways in which conformity occurs.

\subsubsection{Can we defer ethics to management?}\label{s:defer}

We finished the lecture by addressing a notion we had heard from mathematicians several times: that, while it is important to consider ethical issues when carrying out mathematical work, such tasks and responsibilities should fall to management, rather than frontline mathematicians. So we discussed some of the shortcomings of managers in the workplace which might make that approach somewhat impractical or ineffective, and that `deferring to management' to identify and deal with ethical issues simply might not work, even with completely benevolent, well-meaning managers\footnote{The previous parts of the lecture covered the various forms of \textit{malevolence} that might be present in colleagues, managers, or the company as a whole.}.

We first looked at the sorts of training and experience that might properly equip a manager to identify and deal with ethical issues. We suggested that training in the humanities and social sciences, be it philosophy, law, history, psychology, international relations, etc., might be highly conducive to making such evaluations and decisions; much more so than pure training in the technical sciences. But we then put it to students that, if a manager comes from a humanities background (say, history), how will that manager properly understand and evaluate the \textit{mathematical} work being done? How would a historian scrutinise datasets, understand limitations in mathematical models, see how cryptographic capabilities could be exploited, or understand the potential for feedback loops within automated processes? Those with skills in the humanities (alone) might be good at identifying and addressing \textit{some} ethical issues in the work, but a lack of technical understanding may inhibit them, meaning they miss many other deeply-hidden issues that are woven into the technicalities of the work. They might understand societal consequences, but not how to connect or infer that from the mathematics in front of them.

So we then explored how effective a technically-trained manager might be at this instead. Such a manager, perhaps from computer science, engineering, or mathematics, might have a very good understanding of the technical work being done by those they manage. However, we suggested to students that such a manager can often just be an older, more senior version of a young entry-level mathematician. Such a manager might also come from a strictly technical training background, have worked in a frontline technical role for some number of years, and then been promoted to a management role with little more than a short course on basic managerial skills and techniques (probably mostly related to HR issues). Given the scarcity and shallowness of ethics training in most highly-numerate degrees, these managers may have never had the opportunity to develop such skills and understanding. We posed to students that, if an entry level mathematician does not have the right insight to evaluate ethical issues, why would a version of themselves with 5 years experience in the role and promoted to manager have that insight? They might understand the mathematics in front of them, but not how to connect that to the societal consequences.

This led us to posit the final type of manager: the \textit{magic manager}. This is the manager that understands all the technical aspects of the work being done, \textit{and} all the surrounding societal concerns and how the mathematics might link to those and cause harm. We suggested to students that, naively, having such managers might `solve' the problem, even though finding them could be extremely difficult as there are very few people with sufficient training and experience in both deep mathematical work \textit{and} broad social studies. However, the \textit{magic manager} might not suffice either, for one simple reason: time. A typical manager might manage 5--10 people, and so we suggested students consider how much technical output such a group might produce in an average week, and then how long it would take to check all that. Even for something as simple as grading undergraduate mathematics exercises, it can often take a competent mathematician more time to understand a critique a student submission than it would to simply write out the answer themselves; checking the technical work of others is a very time consuming process. And so, even with all the will in the world, it may simply be beyond the time constraints of the magic manager to check \textit{all} the mathematical work of \textit{all} their subordinates to try and find \textit{all} ethical issues contained therein.

And so we put it to students that, if managers with a humanities background are unable to scrutinise technical work for ethical issues, and those with a technical background are unable to connect such work to broader societal concerns, and those with a humanities \textit{and} technical background (magic managers) still lack the time to carry out such evaluations, then who is left? Who can take up the role of looking for and raising any ethical concerns that come from such technical work? We suggested to students that there is only one group left who can do this: the mathematicians themselves. It is up to them to at least initiate the process of thinking about ethics during the development of mathematical work. Of course, they can then seek out management and others in the organisation for further assistance and insight, but the essential first steps must be done by the mathematicians themselves; there is simply no-one else who can do this. They don't need to go out and get a degree in philosophy to do this, as even a modest broadening of their skills so they are a `5\% ethics, 95\% mathematics' person makes a significant difference on their own work, and on the rest of their their team. One doesn't need to know as much ethics as one knows mathematics to be an ethically autonomous mathematician in a world of colleagues who are `100\% mathematics, 0\% ethics'. But we highlighted that such knowledge and skills don't just come by osmosis; mathematicians must put in the necessary time and effort to learn about ethics and responsible development in mathematics, such as by doing a course on ethics in mathematics (like the one they have just completed). It falls on them to do so, as no-one else can do this work for them.

\begingroup
\let\section\subsection
\renewcommand{\bibname}{Bibliography (Lecture 7)}
\renewcommand{\refname}{Bibliography (Lecture 7)}
\addcontentsline{toc}{subsection}{Bibliography (Lecture 7)}
\renewcommand{\addcontentsline}[3]{}
\bibliographystyle{IEEEtran}
\bibliography{bibLec7}
\endgroup

\newpage
\section{Lecture 8: Looking into the future, what more can mathematicians do?}

\fbox{
  \parbox{0.98\textwidth}{
Mathematicians need an acute awareness of the ethical issues they may face during their career, but this is only the first in a sequence of steps required to do mathematics responsibly. Beyond knowledge, actions and changes in behaviour are necessary in order to turn safety and ethics from a theoretical concept to lived practice, and such actions are not always immediately obvious. 

This lecture looks at how mathematicians can take this next step, from increasing their level of ethical engagement, to looking past oaths and codes as `quick-fixes', to operationalising responsible development, and to seeing how to factor in other points of view and put their own thoughts into writing. It shows how they can `be the change' and transition from thinking about ethics, to doing ethics, leading them out from the insights of the course and into more proactive roles as responsible mathematicians.
}}

\subsection{Further tools and resources to be more proactive}

This final lecture was put together to `lead out' from the earlier parts of the course which covered \textit{what} went wrong and \textit{why} it went wrong; we now moved on to start covering \textit{what to do}, and how to actively avoid such harm in the future.

In the very early iterations of our seminars on ethics in mathematics (going as far back as the initial ones in 2016), which just focused on `what had gone wrong and why', we noticed a recurring phenomenon. Students started to articulate comments of the form ``OK, we get it; there \textit{is} ethics in mathematics. Now what can we do to be to be \textit{actively} ethical? What guidelines, tools, and resources are out there to guide us?'' Our initial response was to point out that there are no algorithms or axioms to identify, let alone `solve', ethical problems (in mathematics, or elsewhere), and the technical, deterministic approach that works so well for abstract mathematical problems would not work when trying to deal with ethics. Such problems needed a much more holistic approach, factoring in subjective analysis and determinations, and considering aspects that are hard (or even impossible) to quantify. However, at the time (i.e., in 2016), there were no such tools available that were specifically relevant to mathematical work.

So, we set about creating such insights, tools, and resources, over a period of about 10 years. This 8th and final lecture then became a summary of the tools we had created, so that students knew what was out there, where to find it, and (roughly) what it contained. The works we covered in this lecture were:

\begin{itemize}
    \item \textit{Four levels of ethical engagement} \cite{chiodo2018four}
    \item \textit{A Hippocratic Oath for mathematicians?} \cite{muller2022hippocratic}
    \item \textit{Manifesto for the responsible development of mathematical works} \cite{chiodo2023manifesto}
    \item \textit{The ethical and sustainable concerns triangle} \cite{müller2025triangle}
\end{itemize}
Each work looks at some aspect of what, or how, mathematicians can actively change the way they think about or carry out mathematical work. Each of these provides a different, complementary take on how mathematicians can incorporate ethical thinking into their actions and outputs; there is no \textit{one} guidebook to ethics in mathematics, as it is a deep and complicated area. We give a summary of each work below. In addition, we gave students some direction on how they can begin to \textit{write} about (aspects of) ethics in mathematics, referencing further tools such as \cite{chiodo2024field} and \cite{chiodo2023teaching}, which we also describe below.

\subsection{Four levels of ethical engagement}

\begin{center}
\textit{All arguments and ideas in section come from} \cite{chiodo2018four}.
\end{center}
Engaging in ethics is not a binary on-off state for mathematicians. Instead, they can act on their ethical awareness at different ``levels'', or categories, of engagement. There are, for now, 4 such levels which have been identified (and a zeroth `non-level'). These ``levels'' are not a grading or ranking of `ethical goodness'. Rather, they are a measure of how far away those tasks or actions take the mathematician from their own mathematical work; the deeper and more active the engagement, the further the mathematician is drawn from their everyday mathematics. Specifically, these levels are:

\begin{itemize}
    \item Level 0: Believing there is no ethics in mathematics.
    \item Level 1: Realising there are ethical issues inherent in mathematics.
    \item Level 2: Speaking out to other mathematicians.
    \item Level 3: Engaging with the tables of power.
    \item Level 4: Calling out the bad mathematics of others.
\end{itemize}
A mathematician need not progress through the levels incrementally; indeed, they can jump between them. Nor is the definition of each level formulaic; instead, they are described more by example behaviour. Overall, the idea here is to see that ethical engagement comes in many forms for mathematicians, each of which is a valuable part of the journey towards having an  ethical mathematical community as a whole.

\setcounter{subsubsection}{-1}
\subsubsection{Level 0: Believing there is no ethics in mathematics}

This level needs little explanation; it describes mathematicians who are of the belief that there is no ethics in mathematics. This is where most mathematicians currently sit, often holding the belief that their mathematics is pure, neutral, value-free, and merely an abstraction not mired by social context. Such thinking can be implicitly encouraged through the teaching of mathematics, and the way it is often presented in an abstract manner, with all motivation and context stripped away, leaving only formalism and abstraction.

Mathematicians at this level might engage in deferral of responsibility, either arguing that their work inherently has no ethical weight due to its abstract nature, or that it is somehow the responsibility of others -- the `users' of the mathematical tools or output -- who carry all responsibility for any harms that come from it. At this level, mathematicians tend to have no awareness of the fact that they are the human link between the abstract (decontextualised) truth of mathematics, and the physical-world impact of its use. They do not see that, without them, such harms would not occur, and that their output has a direct causal effect on such impacts, for better or worse.

It is perhaps unsurprising that this is the level where most mathematicians are; their teaching encourages it, the nature of their work alludes to it, and wider society is not sufficiently equipped or knowledgeable to understand and hold them to account for it. And at this level, mathematicians are not in any way `distracted' from their usual mathematical work, and carry on without so much as lifting their head up and looking around at their surrounds. Hence, this is level 0, as the mathematicians here don't move even 1mm away from their own mathematics.

\subsubsection{Level 1: Realising there are ethical issues inherent in mathematics}

At this level, mathematicians understand that they exist in a social and political matrix, and that all work they do carries ethical consequences, be it now, or at some point in the future. They are fully aware that they make and wield tools of vast power, and that these tools can be used both for good, and for harm; there is no such thing as `absolutely harmless' mathematics. Such mathematicians are mindful of their own individual actions, and modify them accordingly in a (mostly) individualistic way, by probing further into the problems they are solving, giving careful consideration to the potential for harm that might come from their work, etc. They are fully aware that their actions have consequences, and that they have the agency to change those actions as an when is needed; an agency they are willing to exercise should the time come.

Here, we also see mathematicians thinking not only about where their work is going, but also where it comes from, in terms of who is paying or funding them, and what motives and objectives such entities might have. They realise that someone (with an agenda) pays them for their work, and that they don't have a God-given right to do mathematics. At this level, actions can be as simple as asking questions about data provenance, to ensure that datasets being used have been sourced legally and ethically, and are representative of the situation they are trying to model. Or it may involve further investigation into the circumstances surrounding a problem, so that they don't just optimise over simple factors such as time and money.

These mathematicians are also aware that such considerations extend into academia -- not just in the context of their technical output and ensuring safety there with issues such as responsible disclosure -- but also in their other academic tasks such as which employers they recommend to their students or would consider withholding letters of reference from on moral grounds. And they also realise that actions at this level can have consequences. They know they may find themselves in positions of professional conflict, or might suffer some employment consequences including dismissal, or may, at the very extreme, be blacklisted from a particular industry or area of work as a result. But they are also prepared for such eventualities, and think ahead to protect themselves and place themselves in the best possible position before embarking on such actions, lest they find themselves being `forked' between betraying their values, and harming their careers.

At the very extreme, mathematicians at this level will go  all-out, and blow the whistle on the actions of their employer or organisation, as has occurred many times in industries such as state surveillance agencies with the likes of Edward Snowden, Chelsea Manning, William Binney, or Thomas Drake. Here, mathematicians take extremal individual steps, ensuring that their values and beliefs are seen through in their activities. Overall, wile still focused on the mathematical work \textit{they} are doing, mathematicians at this level now start to step away a bit from just carrying out technical mathematical tasks, and dedicate some time and effort to ensuring the safe and ethical impact of such tasks.

\subsubsection{Level 2: Speaking out to other mathematicians}

Like all human knowledge, ethical awareness doesn't just `grow on trees'. It needs to either be developed through extensive research and painstaking commitment, or, more commonly, be learned from someone. Of course, some mathematicians will forge for themselves an understanding of the impacts of mathematics in the world and the ethical issues that ensue. But most will not even know to start on such a journey, let alone have the time or insight to see it through. Hence, mathematicians at this level not only realise there is ethics in mathematics, but moreover take it upon themselves to inform and teach \textit{other} mathematicians about it.

It is seldom anyone's `job' to communicate ideas about ethics in mathematics to other mathematicians. Yet, for such ideas to spread efficiently, they need to be explained to mathematicians, just like any other form of knowledge. And so, mathematicians at this level take on this responsibility themselves, to raise awareness among some of the other mathematicians they are in contact with, or have ways to contact. This can be as simple as approaching colleagues in the workplace with their concerns, in order to help said colleagues identify and think about such issues together with them. As managers, they might incorporate discussions about safety and ethics in management meetings and task setting. Or as teachers in an educational setting, they may choose to contextualise some of the mathematics they teach, highlighting both the beneficial and harmful consequences of the mathematics being presented.

But it doesn't just need to be to the people they are already in contact with, as mathematicians at this level may choose to seek out other mathematicians to convey these ideas to; mathematicians they would otherwise have no interaction with. They may set up a standalone teaching module to create a new audience of mathematicians to teach ethical awareness to. They may write blogposts, articles, or other texts, about some issues related to ethics in mathematics, publishing these in a way that reaches a broad audience of mathematicians, many of whom they will never meet. Or they may form collectives of mathematicians with the aim of dealing with ethical issues, either very broadly (such as the student society for ethics in mathematics set up in Cambridge\footnote{\url{https://cueims.soc.srcf.net}}), or in a more specific manner (such as the Just Mathematics Collective\footnote{\url{https://justmathematicscollective.net}} which brought together mathematicians on a platform for social justice).

Whatever their chosen method(s) of engagement, mathematicians at this level all have one thing in common: they have consciously moved on from just considering the ethics of their own mathematical work, and taken the additional step of reaching out to other mathematicians to help them understand the ethics of \textit{their} work. This `journey' to other mathematicians, be it a physical one of walking down a corridor to talk to colleagues, or a metaphorical one of putting time and effort into publishing a piece on ethics in mathematics, takes the mathematician that bit further away from their own abstract mathematical work. Even though this may often only be for short periods of time, it still temporarily detaches them from their own work, which is what distinguishes these mathematicians from those at level 1.

\subsubsection{Level 3: Taking a seat at the tables of power}

The earlier levels describe mathematicians operating for the most part in their `natural habitat': in universities, industrial organisations, etc., within groups containing (almost exclusively) other mathematicians\footnote{Or those from closely related disciplines such as computer science, physics, engineering, etc.}, and at their existing role or job. However, mathematicians at this next level do something different; they move from their usual role and tasks, and go and sit at the tables of power, with (or as) those in managerial or leadership roles. Such mathematicians understand that their role as technical workers means they have limited -- some, but not total -- influence over what mathematics they are tasked to do, how their mathematical output is interpreted, and where such output is used. These interventions can be very localised, with mathematicians offering to join leadership committees or put themselves forward for management roles in their place of work. It can be more general, such as providing expert advice on a regulatory or legislative body. Or it can be very high-level, such as when Fields Medallist Cédric Villani ran for public office in France.

And so mathematicians at this level realise the importance of taking a seat at these tables of power, be it as representatives, managers, or other roles. There, they can present the benefits, and limitations, of the mathematics being requested or used. Mathematics is not a perfect tools or unmitigated force for good; it can cause harm, and someone needs to properly inform leaders about this. Mathematicians at this level know that they need to be at those tables of power, to ensure the mathematics is not given exaggerated influence, or overly-trusted, or used in completely inappropriate scenarios. They are aware that taking credit for the positive impact of mathematical work also requires taking ownership of its foreseeable harms and engaging in relevant debates about what it can and can't do, or should and shouldn't be used for.

But, as mentioned earlier, this level takes mathematicians away from the surrounds they are accustomed to, which has many challenges. First, there is the fundamental difficulty of communicating with those who do not have any technical background. Mathematicians at this level learn to navigate such conversations, and to speak understandably (and not rely on technical terminology), as well as listen intently (and be accepting of non-technical reasoning and concerns). A conversation with, say, a senate committee might involve a senator use the phrase ``What is an algorithm?'', but also ``We need a solution that will appeal to the public'', neither of which are common phrases in a `mathematical' environment. In addition, such interactions can often be confrontational, aggressive, and highly-politicised; all quite common parts of the discourse at tables of power, but perhaps not the sorts of interactions mathematicians might be accustomed to.

Ultimately, mathematicians at this level are taken quite far away from their own mathematical work, or indeed much mathematical work at all. They are no longer just doing a different thing; they are in a different world. And this is what differentiates level 3 from the earlier levels; a detachment so substantial that these mathematicians end up spending more time arguing \textit{for} mathematics than \textit{with} mathematics.

\subsubsection{Level 4: Calling out the bad mathematics of others}

This fourth and final level of ethical engagement addresses the responsibility some mathematicians take on to call out the bad mathematical work of others. Here, mathematicians actively seek out where mathematics is wrong or has `gone wrong' outside of any domain or realm they are directly involved in, and specifically `call it out' to the  wider public; a form of `mathematical hunting'. Those at this level realise that sometimes, only a mathematician can provide the deep perspective and analysis needed to scrutinise a piece of mathematics, and thus that they are \textit{uniquely} positioned to at least identify such problems.  However, this level is particular, in that it splits into two parallel sub-levels.

The first of these, level 4a, deals with mathematicians calling out those who are not well-versed in mathematics using it unethically. That is, calling out those who do not have a full, deep grasp of mathematics, but nonetheless (mis)use it in naive, possibly foolish, and potentially malicious, ways, claiming mathematical results that are simply not true. This can have serious consequences, such as the case of Sally Clark who was (wrongfully) imprisoned for murdering her two children on flawed statistical `evidence' from a paediatrician. Here, the Royal Statistical Society intervened and `called out' the flaws, leading to her eventual release. Acting at this level is, however, counter-intuitive for most mathematicians, as the errors are so manifestly obvious that they hardly warrant mention on mathematical grounds. However, those at this level appreciate that \textit{any} use of mathematics usually provides a `Teflon shield' to an argument, which the lay public will either be unable to understand, or lack the courage and conviction to challenge. These mathematicians realise the importance and impact of a `real' mathematician pointing out `silly' (yet highly consequential) mathematics to wider society. They see the fundamental distinction between intrinsic objective truth in mathematics, and its extrinsic subjective meaning and interpretation in the physical world.

The second of these, level 4b, deals with mathematicians calling out \textit{mathematicians} using mathematics unethically. This is fundamentally more difficult because the `bad mathematics' may be technically deep and numerically accurate, having been done by the mathematically competent. In level 4a the mathematics is not doing what it is supposed to do; in 4b, it is, to the detriment of society. Here, mathematicians understand that there are problems they alone can identify and raise wider awareness of. They can detect the mathematical source of problems in ways the general public has no idea about, let alone realise is harmful. This `sixth sense' is a by-product of many years of training, which they use to see where the mathematics comes from, what it has done, and where it is going.  Such mathematicians seek out where mathematics is being intentionally used to cause harm, such as when Cambridge Analytica used various data science and psychometric tools on social media to manipulate voters in 2016. Mathematician Paul-Olivier Dehaye worked with journalists to help uncover and expose these activities, by helping them understand the technicalities behind, and consequences of, the targeted political advertising campaigns. Mathematicians at this level know that it is technically taxing work, which is why there are so few of them. But they remain driven by the fact that only a mathematician can call out another mathematician.

However, as mentioned, mathematicians at level 4 are  rare; apart from flagging errors in academic research papers, there are few examples of mathematicians calling out the `bad mathematics' of others. And there are two reasons for this. Firstly, it is not obvious to most mathematics how to call out bad mathematics. Figuring out who to tell, and what to tell them, can be complicated, and does not come with the standard training they've received. As mathematics is not formally professionalised with accreditation systems, there is no centralised organisation or `ethics board' to inform. So they need to figure all this out themselves, which disincentivises many. And secondly, there are few incentives for a professional mathematician to put in such (vast) effort to identifying the misuse or misapplication of mathematics outside the profession; they can get no research papers, promotions, or other tangible benefits from their community for doing so. Thus, mathematicians at this level have chosen to take on all these complexities, and give up vast amounts of their time, in pursuit of rectifying such harms, making this level unique for how long and how far it takes the mathematician from their own mathematical work and usual surrounds. However, such work is crucial, as it allows society to then bring its usual response tools -- political, regulatory, judicial, or social -- to deal with the issue, and these mathematicians have come to the realisation that they must do it because no one else can.

\subsection{A Hippocratic Oath for mathematicians? Necessary, but not sufficient}

\begin{center}
\textit{All arguments and ideas in section come from} \cite{muller2022hippocratic}.
\end{center}

\subsubsection{Existing codes of conduct, and calls for Hippocratic Oaths}

The Oath of Hippocrates (OoH), variants of which are sworn by physicians upon admission to the profession of medical practice, has long been synonymous with the idea of `ethical conduct' for physicians. Indeed, one might be led to believe that (most) physicians behave ethically \textit{because} they swear such an oath. However, as shall soon be covered, a better interpretation of this might be that they swear such an oath \textit{because} they are ethical. The difference here is critical; the oath becoming an \textit{effect}, rather than a cause, of embedded ethical insight and behaviour. So might mathematics do well to also have its own `Hippocratic Oath for Mathematicians' (HOFM)? One key obstacle already presents itself here: the full discipline  of mathematics is not professionalised, unlike medicine or similar disciplines that place great importance on ethical conduct (e.g., law). However, one subfield of mathematics is -- actuarial studies -- and makes for interesting analysis in this context.

Actuaries professionalised nearly two centuries ago, and this was done to make the field more reliable, and to gain trust from society. Actuaries have ``all the apparatus of a profession such as medicine'' \cite[p.~4]{muller2022hippocratic}. This includes ``professional bodies that oversee specialised education and license practitioners, compulsory continuing professional development, enforcement and self-governance, and commitments to public service. Ethical codes come as a matter of course and are prominent in the profession’s marketing of itself'' \cite[p.~4]{muller2022hippocratic}. And there is another subfield of mathematics that is (partially) professionalised: statistics. For example, in the UK statisticians have an accreditation system through the Royal Statistical Society (RSS), which ``accredits and governs the actions of professional statisticians in the UK following an extensive code of conduct'' \cite[p.~5]{muller2022hippocratic}. Indeed, statisticians such as William Freeman have been proposing sets of principles -- akin to oaths or codes of conduct -- since the early 1950s.

But these are not the only branches of mathematics that have considered ethics and associated codes or oaths. Ethics in Operations Research (OR) has been considered since the 1960s, with a `Hippocratic Oath' (HO) for OR proposed by Braus in 2002. However, ``by 2009 only 75 researchers
had signed [it]'' \cite[p.~7]{muller2022hippocratic}, suggesting that perhaps ``an oath alone might not be sufficient'' \cite[p.~7]{muller2022hippocratic} to properly address ethics in (parts of) mathematics. Additionally, after the cold war, mathematicians such as Chandler Davis in 1989 ``asked whether the discipline needed a Hippocratic Oath'' \cite[p.~7--8]{muller2022hippocratic}. Unfortunately, this was largely ignored, and in the decades that followed, ``most mathematical societies focused their efforts on producing codes of conduct dealing with academic misconduct such as plagiarism and teaching responsibilities'' \cite[p.~8]{muller2022hippocratic}. And again, in 2009 following the Global Financial Crisis (GFC), Wilmott and Derman produced the Financial Modeler’s Manifesto, which contained a (short) HO for modellers. Likewise, after the Snowden Leaks in 2013, there were explicit calls for more reflection by mathematicians in the security industry regarding ethical behaviour. In 2016, with the rising impact of algorithms and big data, O'Neill called for a HO for data scientists, but at the same time admitted it might be ``insufficient for the task that lies ahead'' \cite[p.~9]{muller2022hippocratic}, acknowledging that oaths or codes \textit{alone} may not suffice. And finally, as late as 2019, mathematician and presenter Hannah Fry called for a HOFM.

An so several patterns start to emerge. Firstly, such calls are usually \textit{reactive} to some large-scale harm that was (primarily) caused by mathematics; it takes a disaster to precipitate such ethical reflection. And secondly, such calls, and associated discourse around ethics in mathematics, have been interpreted as `ethics in X' for some subfield X of mathematics that has (recently) encountered ethical challenges. There has, up to 2019 (the latest call for a HOFM), been very little published work addressing the issue of ethics in \textit{all} mathematics. This comes from the fact that mathematics as a whole, unlike some of its subfields, is not recognised as a profession. Mathematics is widely used, and mathematicians ``take on various job titles throughout their career, including software engineer, statistician, or analyst'' \cite[p.~10]{muller2022hippocratic}. Compare this to physicians where, from the ancient OoH, to today’s modern version, all `physician's oaths' have viewed them as a profession; their actions and domain of work has ``always been clearly delineated'' \cite[p.~10]{muller2022hippocratic}. And so, can all of mathematicians be viewed as one group, and therefore does it really make sense to have a HOFM?

But before proceeding, it helps to ``distinguish between \textit{codes of ethics} and \textit{codes of conduct}, even though that distinction is not clearly made in the names of existing professional codes. A code of ethics is the broader concept. It includes principles applicable to all the ethical aspects of a professional’s work, including its relation to the public interest. A code of conduct, on the other hand, restricts itself to the personal behaviour of a professional in the field, enjoining such virtues as honesty, competence, respect between professionals and ensuring work does not bring the profession into disrepute. In this terminology, the Royal Statistical Society’s code, though called a code of conduct, is broad enough to be a code of ethics. But the codes of the American Mathematical Society, the Institute of Mathematics and its Applications, and the Society of Industrial and Applied Mathematics, are codes of conduct'' \cite[p.~16]{muller2022hippocratic}.

\subsubsection{Difficulties in adapting the `Oath of Hippocrates' to mathematics}

One cannot simply pick up the physician's OoH and replace all instances of `physician' with `mathematician' to get a HOFM; that much is obvious. But the problem is far deeper than that; it is extremely difficult, if not impossible, to directly adapt the OoH for mathematics, as the differences between the two fields are substantial.

The problem begins on the very first line, where the OoH asks the physician to enter into a covenant ``with their peers and their national accreditation body'' \cite[p.~10]{muller2022hippocratic}; no such accreditation body exists for mathematicians, so who is the covenant/oath even to? Next, the OoH is for the ``benefit of the sick'' \cite[p.~11]{muller2022hippocratic} (i.e., dealing with faults in the human body), whereas mathematicians don't just rectify faulty systems; they try to improve completely functional systems, so there is a very different notion of `necessity' here. Furthermore, versions of the OoH have existed for millennia, which has had a ``deep and lasting impact on the culture of medical practitioners'' \cite[p.~13--14]{muller2022hippocratic}. Not only does mathematics lack such an oath; the ``idea of ethical responsibility has not permeated into the disciplinary culture and training (as evidenced by the lack of by the lack of codes, courses, articles and conference papers on ethics in mathematics)'' \cite[p.~14]{muller2022hippocratic}, thus there may be a community inertia against any such oath. And finally, the OoH has the notion to ``do no harm''. But ``how might mathematicians interpret this when it is not entirely clear what harm actually means? It might be easy for mathematicians to restrict this notion to their employer’s interests, but how should they consider impacted parties further afield?'' \cite[p.~14]{muller2022hippocratic}.

Two further, high-level `dissonances' arise here. The first is that ``existing oaths for physicians rely on there being a well-defined patient, and an associated (common) understanding of ``good'' and ``harm'' '' \cite[p.~20]{muller2022hippocratic}. Such a patient has a ``relationship of personal encounter with the physician, who has explicitly agreed to act in their interest, and that personal relationship structures the professional ethics of the relationship'' \cite[p.~20]{muller2022hippocratic}. But for mathematicians, the concept of `patient' -- the entity to which they have a duty of care -- is fuzzy, and ``one needs to consider a variety of stakeholders in addition to the party employing the mathematician: those using the technology (e.g., the users of a search engine), the customers of the company (e.g., companies buying ad space on social media), their employer and the various funding sources, and other impacted parties (e.g. those being advertised/sold to)'' \cite[p.~20]{muller2022hippocratic}. Any, or all, of these could assume the role of `patient', and it becomes particularly complex when they have conflicting interests; which is the mathematician to serve?

The second of these `dissonances' lies in the ``proximity and abstraction of impacted parties. Unlike mathematicians, a physician can usually see the patient and personally inquire to gain a better understanding. Mathematical work usually affects many individuals and groups far away from, and unseen by, the mathematician, and often in ways that the mathematician might never see nor experience first hand. Complex mathematical systems lack feedback mechanisms, of the kind that exist in a doctor-patient relationship, and are difficult to change if mistakes were made. Questions of consent and accountability are already highly debated issues in artificial intelligence and the wider field of ethics of algorithms. All their challenges remain: Can users object and withdraw consent in the same way a physician’s patient could? Do alternative technologies exist and what are the opportunity costs associated with them? Who should be held accountable if things go wrong, and what are the best ways to establish accountability?'' \cite[p.~20]{muller2022hippocratic}.

And finally, we see that although the OoH ``lacks detail of how a medical practitioner should act, this is supplemented by various medical ethics standards, and medical ethics training that is received by almost all medical students. Furthermore, the medical profession is accredited and governed by extensive rules and regulations on what warrants improper conduct. So, in essence, the oath is merely an easy-to-remember reminder of an otherwise deep, complex and well-supported system. The oath alone, though useful, would probably not achieve all of it. Mathematicians have no such systems. There are remarkably few rules or regulations governing general mathematical work, there is no widely-recognised accreditation process that can be revoked, and there is almost no ethical training offered to mathematicians as part of their degree courses. As even the existence of ethics in mathematics is not yet widely accepted, its content would likely be subjected to extensive challenges'' \cite[p.~14]{muller2022hippocratic}.

\subsubsection{What is needed instead?}

With all the above in mind, it appears that simply producing an `oath' (i.e., a commitment to uphold moral principles) might be unsuited for mathematicians, for now. To begin, a better choice may be a code (describing expected behaviours and actions). And much can be learned from  existing related quantitative fields such as engineering, computer science, actuarial studies, and statistics, all of which have reasoned, relevant ideas on these ``basic ethical principles'' \cite[p.~15]{muller2022hippocratic}. One key recurring theme in the codes of these disciplines is the idea of having a `duty of care' and acting in `the public interest' for any work that impacts the public, and the need to take reasonable measures to prevent foreseeable harms to the public. Others include ``honesty/integrity, competence, avoiding conflict of interest, engaging/communicating, protection against information asymmetry, and attending to data issues'' \cite[p.~16]{muller2022hippocratic}. All of these are relevant for mathematics, and thus it would ``certainly be possible to construct a code of ethics for mathematicians'' \cite[p.~19]{muller2022hippocratic}. But would that be \textit{useful}?

Firstly, ``such a code would not (yet) be enforceable. Mathematics is not, and is unlikely to soon become, fully professionalised like medicine, law and the actuarial profession, with accreditation and a powerful professional association capable of disciplining members'' \cite[p.~19]{muller2022hippocratic}. So it may act as a set of rules, but with no consequences for not following them; a potentially futile exercise. Secondly, ``it may also be unclear who such a code would apply to'' \cite[p.~19]{muller2022hippocratic}, since ``even the labelling of mathematicians as mathematicians is problematic, with many graduates in mathematics later taking on the job titles of other professions and identifying as software engineers or computer scientists, sometimes without further professional training'' \cite[p.~19--20]{muller2022hippocratic}. Such a code would be ineffective without a clear target group.

In addition, ``the notions of `good' and `harm' are muddied for many mathematical applications. While a clinical patient is usually fighting a disease or ailment and the physician assists in that fight, mathematicians are not necessarily fighting something generally agreed upon as `bad'. Often, problems and moral questions in mathematical work only become obvious after implementation. Few had thought about the role of mathematics in targeted political advertising before the case of Cambridge Analytica, in which the personal data of millions of Facebook users were harvested without their consent'' \cite[p.~22--23]{muller2022hippocratic}. And so it becomes very difficult to clearly articulate what `doing good' really means.

However, some of the difficulties mentioned earlier have been partially addressed  in other fields. For example, ``those issues of remoteness from impacted parties and the complexity of interrelations between many affected parties are largely common to all quantitative disciplines like engineering and computer science, which has not prevented those disciplines from incorporating useful ethical reflection into codes and education. Here lessons could be taken from existing codes for very mathematical disciplines, such as those of statisticians and the Oxford-Munich Code for Data Science as both communities deal with similar issues. In all such cases, the ethical problems arise from the nature of the work and its impact, rather than, as in medicine and law, from the personal relationship of trust between professional and client. Codes simply need to adapt to that'' \cite[p.~22]{muller2022hippocratic}.

Additionally, we see that ``codes and oaths partially take their value from being a source of reference. Mathematicians facing ethical problems that naturally arise in their work do not need to `reinvent the wheel' if there is an easily accessible and authoritative document that distills previous thinking on similar issues. A code can also be a useful source in litigation. Computer scientists, who face the same problem of a non-professionalised and diversified discipline, have found some value in a code for this reason. So one point of a code is that no-one can then say in court `Well, no-one \textit{told} me I was supposed to be ethical, or what that would mean.' If the law inquires whether a mathematician has breached a duty of care by causing foreseeable harm, a code could add flesh to that general principle. For it to be effective, however, requires its authors to engage with the previously mentioned challenges'' \cite[p.~22]{muller2022hippocratic}. In particular, ``any successful code of ethics must also address the issues of pure mathematical research, which is more comparable to philosophy and art than to engineering'' \cite[p.~22]{muller2022hippocratic}. 

An enduring problem here is that ``reaching consensus on a code of ethics for all mathematicians could take a long time, and an effective middle-ground could be institutional mission statements'' \cite[p.~23]{muller2022hippocratic}. However, to speed up this process, as well as make it more democratic and community-driven, ``individual universities, departments and other employers might introduce or expand their mission statements to go beyond teaching, inclusivity and diversity. These can be tailored to the particular circumstances of the institution and its surroundings, taking into account the relevant cultural and social needs. Such statements might not be universally applicable, but they are quicker to develop, implement, and evolve. Furthermore, they could be tied to existing institutional structures, making it easier to implement tools of accountability and enforcement mechanisms. But most importantly, they would address a concern commonly heard in the context of ethics in mathematics: `Whose ethics?' Those of the people within the institution. With time, mission statements could feed into and provide inspiration for codes of ethics for national and international mathematical societies'' \cite[p.~23]{muller2022hippocratic}. But, with any such `representative' tool such as mission statements, ``to be successful, they need to build on a set of informed, diverse and representative contributors, including those currently institutionally underrepresented'' \cite[p.~24]{muller2022hippocratic}.

\subsubsection{A roadmap for what to do}

So if simply producing a HOFM is not the right approach, then what exactly \textit{should} be done at this point? Surely the gravity of the situation, and how important it now is for mathematicians to consider ethics, warrants some action. Indeed, there are several steps that mathematicians and the organisations that engage them can now take, to begin the process of building an `ethical understanding' into the community. 

To begin, ``in order to both get initial work underway in an unreceptive environment and to take advantage of local knowledge and enthusiasm, individual institutions, such as universities and companies employing mathematicians, should draw up or modify their own small-scale mission statements to specifically include references to ethical considerations of the mathematical work relating to their activities. These can be done quickly, and can evolve in time. It can help create a critical mass of activity relating to ethics in mathematics'' \cite[p.~24]{muller2022hippocratic}.

Following on from that, ``the professional societies should act officially to express an actionable commitment to ethical responsibility. They should call for the inclusion of standalone courses on professional issues and ethics in all accredited mathematics degrees. They should sponsor the development of syllabuses for such courses, with internationally common core content but with variations for local conditions and allowing for advice from industry practitioners'' \cite[p.~24]{muller2022hippocratic}. Along these lines, ``societies should also encourage the incorporation of ethical awareness and training into all mathematics teaching, courses and modules by way of mathematical exercises involving a societal consideration and short stories on the impact of mathematical technologies'' \cite[p.~24]{muller2022hippocratic}.

Next, ``the mathematical societies, in both pure and applied mathematics, should imitate the successful development of ethics codes by other disciplines, codes which foreground the public interest. They should collaborate to sponsor work towards an international code of ethics for working mathematicians. It should go beyond a mere code of conduct for the internal workings of the academic mathematics profession, but instead put at its centre the primacy of the public interest. It should distinguish between the ethics applicable to pure and applied work. The code when finalised should become an official policy document of each society. That would show to both mathematicians and the public that ethics is taken seriously'' \cite[p.~24--25]{muller2022hippocratic}.

And, as a final step in this `initialisation' of ethical understanding, ``to give their commitment to ethics some teeth in the face of malfeasance, the professional societies should institute complaints procedures to which allegations of unethical professional behaviour by mathematicians can be referred for investigation. If allegations are found to be proved, an official statement by the society should make public the results of the investigation. Despite not carrying any legal standing, they may still deter unethical practice among mathematicians through general community awareness. While the development of mission statements, courses, and codes of ethics is crucial, it is eventually needed to go beyond passive documents, and to create an environment actively conducive to ethical mathematical practice'' \cite[p.~25]{muller2022hippocratic}.

None of these steps are easy, and all of them take time, effort, and commitment to carry out. And they only outline the \textit{beginning} of the process of creating a consensus of ethics in the mathematical community, rather than `solve' it in its entirety. But these would at least put mathematics on the right trajectory, so that in a generation or two the community might take the view that `\textit{of course} there is ethics in mathematics', and have the associated training, understanding, and incentives to meaningfully act on that view. At that point, a HOFM becomes a natural final step in the process.

\subsection{Manifesto for the responsible development of mathematical works}

\begin{center}
\textit{All arguments and ideas in section come from} \cite{chiodo2023manifesto}.
\end{center}

\setcounter{subsubsection}{-1}
\subsubsection{The idea of a `workflow' of ethics}

The first step towards ethical awareness in mathematics is an understanding of the sorts of harmful outcomes that may occur; a case `what we don't want to happen'. The next step is, therefore, what we \textit{do} want to see happen, and what sorts of outcomes are considered desirable or `good'. This is already slightly more complicated, and needs to be pre-emptive, rather than retrospective reflection. It is from this that many `ethics frameworks' are spawned, giving collections of desirable end-outcomes from technological work in fields such as AI. However, from a mathematician's point of view, this only solves half the problem; knowing what the end outcome should be doesn't necessarily help the mathematician see how to get there. This is akin to giving a mathematician a differential equation and told that an acceptable outcome is a function that satisfies that equation; this doesn't tell the mathematician how to find it, and they are left to their own devices to figure out how.

And so, there is a need to move beyond such `endpoint-oriented' frameworks for ethics and safety, and go to the next step: having operationalised guidelines on how to reach such desirable outcomes. With this in mind, a `workflow' for mathematicians was constructed  in \cite{chiodo2023manifesto}. This `follows' the typical steps and actions mathematicians take during a work project, identifying ethics and safety concerns along the way. Presented as 10 key stages, or ``pillars'', it breaks down the actions taken by mathematicians into highly granular activities and decisions. Below is a summary of the key ideas from each pillar; full details can be found in \cite{chiodo2023manifesto}.

\subsubsection{Pillar 1: Deciding whether to begin}

\begin{center}
\textbf{\textit{Why are you providing this mathematical product or service, \\and should you even do so?}} 
\end{center}
This pillar explores fundamental initial questions around the motivation for the problem or project, such as whether it should be solved at all, what good or harm might come from doing so, and what would happen if the work was not carried out. It then moves on to scrutinising funders and beneficiaries, exploring who is paying for the work, what their interests are, and who might end up owning or controlling the outputs at the end.

The pillar also covers how to weigh up the consequences of the mathematical solution, seeking to understand who wins and who loses, how to weigh up the pros and cons, and whether it is possible to improve fairness. It breaks down the ethical issues with the mathematics into those that are standard, obscure, or totally unexpected, by suggesting checks on these three distinct pitfall types, evaluation of what deliberate harm is possible, and consideration of what unintended harm is possible. And it scrutinises why a mathematical solution is even appropriate, exploring how  mathematics actually helps in this instance, whether a non-mathematical approach was considered, and why using mathematics is the best option.

And finally, this pillar looks at where  ethics sits within the project plan,  questioning the resources and priority being assigned to ethical considerations, whether the project can truly afford an unethical outcome, and whether `the project' may become more important than its outputs. It also explores having someone in charge, accountable and responsible at all times, questioning how responsibility and leadership for all parts (including ethics) has been assigned, what quality control measures will exist, and whether there is a disaster response plan.

\subsubsection{Pillar 2: Diversity and perspectives}

\begin{center}\textbf{\textit{Do your co-workers, superiors, and you have sufficient perspective, \\and do you understand the limitations and biases in your thinking?}}
\end{center}
This pillar addresses the value and importance of diversity in a development team, covering how a diverse team produces a better output and has value beyond just profitability, how a team’s biases can manifest in their product, and the need for putting in effort to ensure diversity at all levels including at the tables of power. It interrogates the skills, viewpoints or diversity a team might be missing, by exploring what knowledge beyond mathematics is needed in the room, whether there is more than just diversity of skills and training, and how to maintain an appropriate set of skills, expertise, and comprehensive representation in the team. And it covers how to create a team with a diverse background, expertise, and perspective, by outlining who to bring in and where to find them, how job adverts can make diverse people want to join, and how to hire with diversity in mind.

The pillar also looks at how to bring a diverse team together to work, exploring whether the team’s values allow for diversity, whether the workplace enables interactions that nurture diversity, and how the workplace is meeting diversity and accessibility requirements. It goes over how to harmonise a diverse team to get work done together, by seeing if team members know to interact with others, if people with different characteristics are actually working together and have effective ways to interact, and whether the team uses their diversity to consider underprivileged groups. It explores creating an environment which welcomes differing viewpoints, by explaining how to keep the peace and balance in a diverse team, how to listen to people’s feelings and adapt when needs change, and how to break groupthink. And it covers how to implement the output or recommendations of a diverse team, through listening enough to provide constructive feedback, rewarding people who find problems rather than punishing them, and making changes after hearing suggestions.

\subsubsection{Pillar 3: Handling data and information}

\begin{center}\textbf{\textit{Are you using authorised and morally obtained datasets, \\in a responsible manner?}}    
\end{center}
This pillar looks at how to respect data, by addressing what the `data basics' are, and how to understand the human connection of data. It analyses what data to collect and use, by seeing whether the intended datasets cover enough to do the job, whether the datasets are overkill for the task, and whether all the data is actually relevant. It asks how much understanding there is of the desired data, by looking at what is in the data and where and how it originates, by preventing any abstraction of the data from inhibiting empathy for its subjects, and by seeing how use of the data might cause harm. And it queries the legality of the data usage, by exploring the relevant data laws, the allowed uses of the data, and how the data supply chain complies with the law.

The pillar then examines the the risks associated with data, by studying how inherently risky the data is, and how to prepare for data loss. It goes over proper data privacy and transparency procedures, such as how the privacy policy is articulated and communicated, how to practise good data hygiene so that people can find and excise their data from the work, and whether the data and policies end-game is robust. It also checks whether the data is properly secured, through  `data minimising' practices, controlling data storage and access, and regularly updating data protection measures. And it explains how to collect and create data in an ethical way, by avoiding tricking others into giving away their data or otherwise ripping them off, by ensuring people consent to how their data is used, and by making sure people know what data is being collected off them.

\subsubsection{Pillar 4: Data manipulation and inference}

\begin{center}\textbf{\textit{Do you have the expertise to properly manipulate data \\ensuring quality and ethics?}}    
\end{center}
This pillar explores providing sufficient training for a team to work with data, such as whether the team is trained to manipulate data, whether they are aware of how data is used, and whether they understand common biases for the problem domain and how to detect such biases in datasets. It explains how to work with relevant qualitative data, by going through how to detect obvious red flags with the data, how to check whether the data distribution matches the application domain, and whether the end product might perpetuate a problematic world. It also looks at how to work with qualitative labelled data, explaining how to care for the labellers and ensure the labels are of high quality and will get better over time, how to avoid unexpected surprises from the production environment, and how to pre-process the data and address any `out of distribution' problems. And it covers  documenting any problems with the data, such as documenting the whole data pipeline, noting any unwanted properties of the data, and putting all biases into writing.

The pillar also examines whether and how the data can be valuable to wider society, by emphasising awareness of the value of the data, articulating what further consent should be obtained to share the data more widely, and describing how to decide if it would it be helpful, desirable, and safe to share the data. And it details how to make the data available to others safely, by examining whether the entities being given the data can be trusted, whether such entities might make use of the data without taking a copy away, and whether both the laws of the land, and of mathematics, protect any data passed on.

\subsubsection{Pillar 5: The mathematisation of the problem}

\begin{center}\textbf{\textit{What optimisation objectives and constraints have you chosen, and what are their real-life consequences? Who might the other impacted parties be?}}    
\end{center}
This pillar investigates doing mathematics responsibly to solve the problem, by asking about the `bigger picture' of who wants the work and why, querying who needs to know and understand the work and take responsibility for it, and checking whether there are any conflicts of interest or thefts of intellectual ideas. It explains how to do mathematics with ethics ever-present, by checking if everyone is thinking about ethics regularly throughout all their work, by encouraging everyone, including technical people, to talk about ethics, and by ensuring that technical and non-technical people are communicating.

The pillar also explores difficulties in making a mathematical setup match the problem being solved, by scrutinising how the world has been made to fit into the mathematics, whether any simplifications of the mathematics were made just to make it `nice', and what might be lost in the mathematisation. It questions the use of mathematical truth, asking how precise the use case of the mathematics is and whether users are left to `fill in the gaps',  whether precise mathematics is being deployed meaningfully and not being misleading by saying something is `provably true', and the ways that true mathematics can still cause harm. It also addresses ways to select what mathematics to use, by analysing what mathematical shape the world is assumed to take and what mathematical constraints have been imposed on it, whether trivial mathematical decisions have been made without sufficient reflection, and whether mathematical decisions have compounded to deviate from reality. And it addresses how to optimise for everyone's benefit, by querying what is being optimised over and who the mathematics affects, how the mathematics is projecting into the world, and who will be helped or harmed by the work.

\subsubsection{Pillar 6: Communicating and documenting your work}

\begin{center}\textbf{\textit{Are you properly considering how to comment and document your work \\and communicate the results to those who need them?}}    
\end{center}
This pillar looks at how to properly communicate and create awareness about the mathematics, by ensuring everyone is on the same page and knows the project’s assumptions, by suggesting how to talk about decisions with colleagues, and by encouraging discussion about transparency and fairness. It covers how to talk about the non-mathematical aspects of the work, such as how to discuss responsibility assignment, how to match the means of communication to fit others, and how to keep people updated about changes and progress. And it explores methods to raise and document ethical issues about the mathematics, such as allowing for ethical issues to be raised safely and securely, and updating and talking about ethics regularly.

The pillar also covers ways to document background decisions and reasoning, such as commenting any decisions made, and commenting the reasoning for such decisions. It also goes over how to continuously document mathematical work, by regularly updating documentation and recording the things that must legally be recorded, by ensuring all relevant people can access and understand such documentation, and by making handover documents ensuring continuity.

Finally, the pillar addresses the need to speak with domain experts, by seeing whether there is any actual discourse with domain experts, whether enough of the team can access the domain experts and `speak their language', and whether anyone is actually listening to the domain experts. And it gives guidance on how to communicate with the external world, such as the need to know the audience and have clear motives for engaging in communication, the need to curate what to convey and effectively convey it, and the importance of communication speed and style that matches the surrounding need and also conveys what might go wrong.

\subsubsection{Pillar 7: Falsifiability and feedback loops}

\begin{center}\textbf{\textit{Is your work falsifiable, and can you handle its large-scale impact \\and any feedback loops that arise?}}
\end{center}
This pillar examines what failure means for the given mathematical system, by considering whether the mathematics is data- or theory-driven and how each fails, what it means for the mathematics itself or the overall product to `fail', and whether the work is falsifiable and the underlying mathematics scientifically proven for the application. It explores how the application environment might deviate from initial assumptions over time, by seeing whether the modelled reality or human behaviour changes over time or adapts to the product, by checking what else the deployed product changes in the world, and by monitoring whether the mathematics no longer matches reality or faces unexpected scenarios.

The pillar also considers how much the product impacts the world, by evaluating whether the circle of impacted parties expanding, how big the overall impact could get, and whether the impact scale might compound errors and creating systemic problems. It deals with how the work can amplify ideas, patterns and problems, by seeing whether the product might leverage its scale to cause harm, whether the product might potentially amplify existing problems, and whether the team is prepared to deal with large-scale problems that might arise.

Finally, the pillar explores how the product interacts with the world, examining what else the product might interact with and how, testing how users and the system interact and whether this can distort peoples’ perceived reality, and seeing how others might misuse the product or what new behaviour it incentivises. And it covers how to identify and curtail specific feedback loops, by checking how the product changes the world, how the product changes as the world changes, and thus how to detect and dampen feedback loops.

\subsubsection{Pillar 8: Explainable and safe mathematics}

\begin{center}\textbf{\textit{Is your mathematical output explainable, and followed up with \\proper monitoring and maintenance?}}
\end{center}
This pillar analyses what the system is intended to do, by asking what the team does, and doesn’t, directly want it to do, what the rules, standards, and best practices of the industry are, and how risky the solution implementation is. It explores how well the team understands what it does, and doesn’t, do, by checking if the team has the expertise to scrutinise it and how this is done pre-deployment, by verifying whether the system’s use and behaviour meet expectations and what limiters and failsafes have been put in, and whether it can behave in unexpected, dangerous ways. And it covers how to explain the product's functionality and outputs by establishing who it needs to be explained to and why, what functionality needs to be explained and whether this can be done sufficiently for the risks involved, and whether interpretable and explainable mathematics is being used that can be effectively externally audited.

The pillar also gives ways to know when the product is causing problems, such as monitoring the system post-deployment, giving team members avenues to report problems without retribution, and enable users or impacted parties to report problems directly. It examines how to keep the product safe after completion or deployment, scrutinising whether there is an effective maintenance team with responsibility and accountability effectively conveyed, how long after deployment the system will be maintained, and whether corrections can be made sufficiently quick. And it probes what adjustments are being made to the system based on the monitoring, by seeing if human or technical corrective steps are necessary and whether the work environment is even open for change, if adjustments are creating more problems or require a rethink of things, and whether everyone whose work gets overwritten is being cared for.

\subsubsection{Pillar 9: Mathematical artefacts have politics}

\begin{center}\textbf{\textit{Are you aware of other non-mathematical aspects and the political nature of your work? What do you do to earn trust in yourself and your product?}}
\end{center}
This pillar considers the environmental cost of the work, by evaluating whether the product and its development are environmentally friendly, whether it is environmentally sustainable, and whether it might be incentivising excessive consumption by others. And it addresses considerations of the product's social cost, such as whether all workers are being properly paid and treated well, whether the project is a good use of money, and who might be hurt by what the product replaces.

The pillar also examines the politics of the mathematics, by investigating who the team's friends and ideological enemies might be, whether the mathematics causes too much upheaval, and who might be able to hurt the team or product, and how. It scrutinises the politics of the team itself, asking what the team's ethics and politics are and whether they might be hiding from politics through mathematics, whose toes the team might be treading on, and who (not) to work with and how any such associations might cause harm.

And finally, the pillar interrogates how society views the product, by asking how understandable the whole operation is and whether everyone can see what the product is and does, by emphasising the importance of understanding the public’s perceptions and whether the good that comes from the work has actually been publicised, and by probing whether the solution will be widely accepted by society. It questions whether the overall operation is trustworthy, by analysing what trustworthiness actually means, how trustworthy the team, mathematics, and third-party components are, and whether the project is genuinely earning trust rather than `faking it', or misdirecting people, through the mathematics. And it describes how to build a trustworthy product, by introducing how trustworthiness can be built into the system and ensuring people know the system’s capabilities and limitations, by seeing whether users can know and  understand when it fails thus avoiding undue faith being placed in it, and by ensuring people are interacting with it safely and checking for any migration of its use to other domains.

\subsubsection{Pillar 10: Emergency response strategies}

\begin{center}\textbf{\textit{Do you have a non-technical response strategy for when things go wrong? \\And a support network, including peers who support you \\and with whom you can talk freely?}}
\end{center}
This pillar explores the ways in which any project may go wrong, such as what if the product collapses, of if the team disbands, or if the overall operation is acquired. And it looks at the specific risks that apply for the project, studying what exactly can go wrong and who is even defining wrong, how bad each problem can get and what is needed to solve them, and what the timescale and interdependency of each problem is.

The pillar also covers emergency response capabilities, seeing whether there is a process for raising the alarm, whether such a process can spot and respond to problems fast enough, and what might be missing in the `emergency response toolbox'. It stresses the importance of a non-mathematical response plan, querying which parts of the response go beyond just doing more mathematics, whether all emergency roles been assigned and communicated with a clear and regularly reviewed emergency plan, and whether the overall setup can defend itself from external attack and individuals  have someone to talk to when things get hard. And it outlines how to execute a response plan, by checking what the problem actually is, by having ways to bring everyone into the `crisis room', and by evaluating whether the problem fits into the existing plan and how to change the plan if necessary.

And finally, the pillar emphasises the need for a proper debrief after an emergency, by checking what emergency actually happened and why, checking what did or didn’t work in the response, and evaluating why the response succeeded or failed. And it covers incorporating lessons learned and updating the emergency plan, by reviewing what lessons can be taken away from the emergency itself, what lessons can be taken from the response, and what can be changed to be prepared for the future.

\newpage
\subsection{The ethical and sustainable concerns triangle}

\begin{center}
\textit{All arguments and ideas in section come from} \cite{müller2025triangle}.
\end{center}

\setcounter{subsubsection}{-1}
\subsubsection{Three areas of concern: Mathematics, Community, and Society}

No two mathematicians are the same in terms of their approach to `being a mathematician', and what their aims and objectives are as they operate within the discipline of mathematics. And one aspect that distinguishes mathematicians is the sorts of things they are `concerned' about. That is, the things they put in time and work to `protect'. The work \cite{müller2025triangle} identifies three core areas of concern for mathematicians; aspects of the discipline that they might seek to advance, and worry about protecting and preserving. First introduced in \cite{chiodo2024field}, these three `concerns' are:

\begin{enumerate}
    \item A concern for ``ensuring that (research) mathematics remains a stable and continuing body of knowledge'' \cite[p.~4]{müller2025triangle},
    \item A concern for ``our mathematical, and other, related communities'' \cite[p.~4]{müller2025triangle}, and
    \item A concern for ``the impact of mathematics on wider society and our planet'' \cite[p.~4]{müller2025triangle}.
\end{enumerate}
These three concerns, which can be abbreviated to \textit{Mathematics}, \textit{Community}, and \textit{Society/Planet}, become drivers for how mathematicians dedicate their time and effort. They are by no means mutually exclusive, and (as will be elaborated later), a mathematician can hold multiple, or even all, of these concerns simultaneously. No one of these concerns is any more or less important that the others, and how much time and effort a mathematician devotes to each concern is their own personal choice.

\subsubsection{A concern for research mathematics}

The mathematicians who carry this concern are thinking about how to preserve, add to, and pass down, the corpus of mathematical knowledge that currently exists. They seek the preservation and advancement of mathematical theorems, processes, and other such (primarily) technical insights. Such concerns focus on \textit{what} is known, rather than \textit{who} knows is or \textit{how} it is used. Here, mathematics is still seen as a very human endeavour, and so simply having theorems written down in books or papers that nobody understands does not constitute `knowledge'; they aim to keep (various parts of) mathematical insight in the active knowledge base of the community of mathematicians, in sufficient number so as to propagate it on to the next generation.

A key component of this concern is the furthering of mathematical insight. To this effect, mathematicians with this concern place a lot of emphasis (and effort) into encouraging the `best and brightest' students to study mathematics, so that they (the students) may later contribute to this pool of collective knowledge as research mathematicians. This `best and brightest' is interpreted by these mathematicians as those whom they perceive will be most likely to prove new and `exciting' results; the very notion of `exciting' is also a subjective one that these mathematicians seek to foster, convey, and propagate down the generations.

Overall, this concern reflects a desire to maintain and continue the \textit{operation} of the research mathematics community, with the aim of ensuring that it remains a thriving research discipline that can continue to advance the frontier of mathematical knowledge as effectively as possible.

\subsubsection{A concern for the community of mathematicians}

The mathematicians who carry this concern are thinking about how to care for the community of people who study and work in mathematics, either now or in the future. They seek to ensure that the discipline becomes as inclusive and welcoming as possible for all those who wish to do mathematics. Here, the interpretation of `mathematics' is very broad, and includes those who learn and carry out mathematics that it not at, or even anywhere near, the frontier of human knowledge. Such concerns focus on \textit{who} knows and uses mathematics, rather than \textit{what} is known or \textit{how} it is used. The (widely-interpreted) community of mathematicians is the focus of this concern.

A key component of this concern is how to expand the mathematical community as far as possible, and enable as many people as possible to do as much mathematics as possible. Greater value here is placed on having a lot of people learn some mathematics, rather than some people learning a lot of mathematics, and this is reflected in how those here set up their teaching and outreach processes. Mathematicians with this concern see mathematics not just as a tool to carry out tasks, but also as an enabler, empowering those who can wield it to further their own livelihoods and resist other forms of power and oppression.

Overall, this concern reflects a desire to expand the mathematical community and enable those in it to enjoy and flourish in a world where they can learn and use more mathematics than existing setups might otherwise enable. It is a concern that focuses on mathematicians in their own right, rather than on the mathematics produced.

\subsubsection{A concern for the impact of mathematics on society and the planet}

The mathematicians who carry this concern are thinking about how to assist all those in society who are impacted or harmed by mathematical work in some way, including ecological effects on the planet itself. They seek to ensure that mathematics is used in a way that helps all of humanity, and steer it away from harmful large-scale effects, be they deliberate or completely unintentional. Such concerns focus on \textit{how} mathematics is used, rather than \textit{what} is known or \textit{who} knows. The focus of this concern, in terms of what is being `protected and advanced', is society (and the planet) as a whole, and not any particular group(s) of mathematicians or the knowledge they create.

A key component of this concern is how to understand what impact mathematical work has on the world. This is broadly interpreted to include adjacent mathematically-powered disciplines such as computer science and physics, as well as mathematics which is not `cutting edge' but nonetheless might have significant impact on how the world works. Here, mathematics and those who wield it are being actively scrutinised. Mathematicians with this concern do not see mathematics as an unmitigated tool for good, but rather as a double-edged sword, with the capacity to do vast good, but also inflict vast harm. It is this potential social and planetary harm that they seek to avert with their efforts.

Overall, this concern reflects a desire to protect wider society from any harms or injustices that could come from mathematics and its use. It seeks to promote human prosperity as a core objective of mathematical work, rather than focus only on the work itself or those carrying out such work.

\subsubsection{A triangle of concerns and how to interpret it}

As mentioned above, these three concerns -- mathematics, community, society/planet -- are by no means mutually exclusive; most mathematicians hold all three concerns simultaneously, in some triple-ratio. How a mathematician choses to split their time and efforts between these three can be illustrated with the use of a `triangle of concerns'; a tool that is common in the social sciences, and may even be familiar to some mathematicians in the context of the `research-teaching-administration' triangle that is often presented (whereby one identifies their time split between research, teaching, and administration duties). 

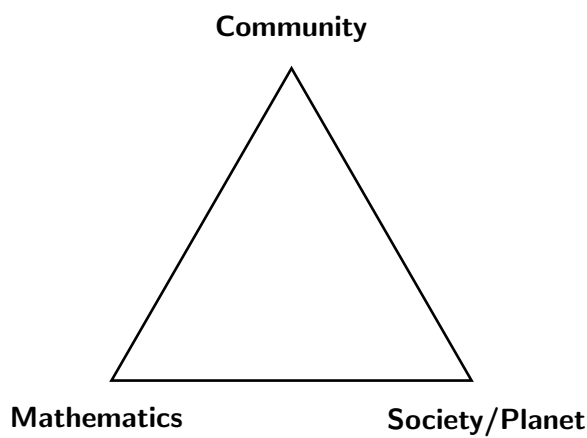
\begin{figure}[H]
    \centering
\begin{tikzpicture}[
    scale = 0.5,
    line width=1pt, 
    font=\sffamily,  
    tension_arrow/.style={
        double distance=4pt, 
        double={rgb,255:red,50;green,90;blue,120}, 
        -implies, implies- 
    }
]
    \coordinate (top_v) at (90:5.5);
    \coordinate (left_v) at (210:5.5);
    \coordinate (right_v) at (330:5.5);

    \draw (top_v) -- (left_v) -- (right_v) -- cycle;

    \node at (top_v)  [above, yshift=2mm] {\textbf{Community}};
    \node at (left_v) [below, xshift=-2mm, yshift=-2mm] {\textbf{Mathematics}};
    \node at (right_v) [below, xshift=2mm, yshift=-2mm] {\textbf{Society/Planet}};

\end{tikzpicture}
    \caption{The \textit{Ethical and Sustainable Concerns Triangle} (ESCT)}
    \label{fig:ESCT}
\end{figure}

So how does this \textit{Ethical and Sustainable Concerns Triangle} (ESCT) work for mathematicians? Simply put, one can take an equilateral triangle and label each vertex with one of the three concerns described above. Then, as each point inside the triangle corresponds to a weighted sum of the (vector) vertices with weights adding to exactly 1, the 3-way split of time and effort between the 3 areas of concern correspond to a unique point within the triangle. By working out these `weights', a mathematician can then broadly see where they sit within the triangle; that is, how much weight they give to each area of concern. Within this triangle, a point very close to a vertex corresponds to being primarily concerned with that vertex, a point in the middle of an edge between two vertices corresponds to having roughly equal concern for those two vertices, and a point in the middle of the triangle corresponds to having roughly equal concern for all three. It often helps to visualise this as an equilateral triangle drawn out in full, which is given in Figure \ref{fig:ESCT}.
\noindent But is this a useful analytical tool? Well, yes, when one considers that some disagreements within the mathematical community can stem from different parties each `caring about' or `protecting' one (or more) of the areas of concern, and finding themselves at odds with those who care about other areas. Without a proper \textit{understanding} of what the other party is concerned about, and \textit{why} they have that concern, such disagreements can quickly descend into conflict; everyone cares, just about different things. By seeing that other mathematicians might have differing `concerns', one can work with such observations to seek a compromise outcome that might be mutually acceptable.

For example, one of the authors of \cite{müller2025triangle} ``gave a talk to an audience of academic mathematicians on teaching ethics and embedding it in everyday (university) mathematics teaching - a talk firmly rooted in the `Society/Planet' vertex of the ESCT. During the Q\&A session after the talk, one mathematician asked (paraphrased), `How can I fit this into my teaching when the curriculum is already packed with mathematics? I have no spare space or time for this.' - a concern firmly rooted in the `Mathematics' vertex of the ESCT. Immediately afterwards, another mathematician asked (paraphrased) `I’m concerned that this teaching might affect how my students seek out jobs, and that on ethical grounds they might choose not to apply for certain types of work, thus reducing their opportunities.' - a concern firmly rooted in the `Community' vertex of the ESCT. To the first question, the speaker replied with `Rather than optimising only over the (technical) mathematical content being taught, might you consider reducing that slightly to incorporate some ethics teaching; material students might otherwise never see at all?' To the second question, he replied with `Rather than optimising only over the interests and careers of your students, might you consider adjusting that slightly to incorporate the interests of those who your students may impact in society through their mathematical work?' Both answers were well-received by those posing the questions, and the wider audience. Here, the ESCT helped to quickly identify different concerns, acknowledge their individual validity, and then seek out the possibility of a compromise.'' \cite[p.~6]{müller2025triangle}.

\subsection{How to put  insights from ethics in mathematics into writing}

Writing about ethics in mathematics is not a straightforward process. While it has been done to varying extents by mathematicians (including ourselves) in the past, up until recently there were no guides or other assistive tools available to facilitate this; mathematicians needed to `reinvent the wheel every' time. Things have now changed, and tools such as \cite{chiodo2024field} and \cite{chiodo2023teaching} now exist.

The first `lesson' to learn here is \textit{why} a mathematician might want to write about ethics in mathematics. This is covered extensively in the ``Field Guide to Ethics in Mathematics'' \cite{chiodo2024field}. It explains why ``ethics in mathematics satisfies the characteristics of a problem worthy of a mathematician’s attention: it is interesting, tractable, nontrivial, well-defined, and (morally)
good'' \cite[p.~940]{chiodo2024field}. To expand on these in more depth (all summarised from \cite{chiodo2024field}):

\begin{enumerate}
    \item Ethics in mathematics is  \textit{\textbf{interesting}} because it creates more problems to work on. It is easy to see that ``Mathematics creates ethical questions, and in return ethical questions generate interesting mathematics'' \cite[p.~940]{chiodo2024field}. Such problems include algorithmic fairness, which `` `led mathematicians and computer scientists to develop new metrics to measure individual and group fairness, connecting their mathematics with ethics, law, and philosophy'' \cite[p.~940]{chiodo2024field}. It also creates a rallying point for community building, as to ``properly address many ethical questions raised by mathematical work, we will need to come together as a community'' \cite[p.~941]{chiodo2024field}, including working on and writing about interdisciplinary topics with diverse teams, producing new teaching resources to help students understand this new field, and building bridges between ``those doing mathematical work, and those considering its effects on society'' \cite[p.~941]{chiodo2024field}.
    \item Ethics in mathematics is  \textit{\textbf{tractable}} because it contains many new and unaddressed problems. As a new area of study, ``ethics in mathematics has received comparatively little academic attention so far, thus there remains much to be done: from the foundational level [of establishing that there are indeed ethical issues in mathematics] \ldots  to the practical applications of turning our theoretical knowledge into being morally good mathematicians'' \cite[p.~941]{chiodo2024field}. This will require incremental contributions, each building on previous works. The problems that exist here, such as AI's threat to humanity, may seem daunting. And ``as mathematicians, we can’t solve any of these problems by ourselves; they’re simply too big. But by working on ethics in mathematics, we \textit{can} find ways to make a positive contribution \ldots and address issues only mathematicians are properly placed to address, or even identify!'' \cite[p.~941--942]{chiodo2024field}. And it is now clear that ``ethics in mathematics has become an academic discipline, including peer-reviewed articles, books, [and] conferences'' \cite[p.~942]{chiodo2024field}, so such writings certainly have avenues for publication. 
    \item Ethics in mathematics is  \textit{\textbf{nontrivial}} because not even the `questions' are obvious, let alone the answers. When setting up an optimisation problem, choosing the `right' objective, and the `right' constraints, is ``often a highly contentious, and thus nontrivial, question in itself'' \cite[p.~942]{chiodo2024field}. However,  it’s ``not just the ethical problems themselves that pose a challenge; bringing together the right people to address these can also be challenging'', as doing ethics in mathematics ``often requires help, support, and contribution from others in the mathematical community'' \cite[p.~942]{chiodo2024field}. But ``growing such a community is hard, as not everyone will agree on what should be done, let alone how'' \cite[p.~942]{chiodo2024field}. And ethics in mathematics is not just retrospective, as ``offering insight into what is, or has, gone wrong in the mathematical process is a good start, but suggesting ways to resolve them is often more appreciated'' \cite[p.~942]{chiodo2024field}. And so, the fact that ethics in mathematics is hard ``shouldn’t be surprising. Its problems require a good grasp of many hard areas, including mathematics, ethics, philosophy, and sociology of science, as well as people and context. All of these are far from trivial'' \cite[p.~943]{chiodo2024field}.
     \item Ethics in mathematics is  \textit{\textbf{well-defined}} because ``reasoning in ethics can be just as rational as reasoning in mathematics'' \cite[p.~943]{chiodo2024field}. We see that ``just like mathematics, much of ethical reasoning begins with finding good starting points or axioms, from which we build further'' \cite[p.~943]{chiodo2024field}, whether this be utilitarianism, virtue ethics, or any other of a long list of ethical positions. And while these positions may disagree, ``they [all] produce well-thoughtout answers subject to an initial set of assumptions; something we’re all very familiar with as mathematicians'' \cite[p.~943]{chiodo2024field}. So, to reconcile such differences and ``do ethics in mathematics effectively'', we may need to ``seek consensus from colleagues on what assumptions to make, what to optimise over'', etc.~\cite[p.~943]{chiodo2024field}. But of course, ``consensus is again a well-defined notion, which as mathematicians we use all the time in scenarios such as multirefereed papers, PhD exam boards, and hiring panels'' \cite[p.~943]{chiodo2024field}. And as for \textit{who} should be dealing with ethical issues in mathematics, ``most mathematics departments have no ethics committees to guide mathematicians with the many ethical issues specific to our subject'', and so we are often left to ``do it ourselves'' \cite[p.~943]{chiodo2024field}. Fortunately, there are ``well-defined research areas and standard methods within ethics in mathematics'', as well as ``many well-known problems: AI safety, cryptography, finance, mathematics in warfare, and more'' \cite[p.~943]{chiodo2024field}.
     \item Ethics in mathematics is  \textit{\textbf{good}} because as mathematicians ``we can empower ourselves, allowing us to better navigate the world of mathematics and its culture'' \cite[p.~944]{chiodo2024field}. In particular, ``we expand our knowledge and expertise so that we can do more interesting and advanced work with better impact and minimal harm. We can extol the virtues of our work, making it more attractive for employers and funders. We can develop greater kindness and humility, towards those near to us, as well as those impacted far away whom we might never see. We can become part of new and different communities, interacting and working with people from other backgrounds and disciplines. We can find new opportunities for work and employment, bringing our now-broadened understanding of the role and place of mathematics into other domains and jobs. And we can be happy knowing that our work and input has made a positive contribution to the field, to those around us, to those in need, and to wider society'' \cite[p.~944]{chiodo2024field}. And so overall, we see that ``doing ethics in mathematics helps us to become better mathematicians: better people doing better mathematics for a better world'' \cite[p.~944]{chiodo2024field}.
     \end{enumerate}
The second `lesson' to learn here is \textit{how} a mathematician might write about ethics in mathematics. This is covered extensively in the ``Teaching Resources for Embedding Ethics in Mathematics: Exercises, Projects, and Handouts'' \cite{chiodo2023teaching}. There, several longform writing tasks are presented for students to work through and gain gain practice researching ethics in mathematics topics, and writing these up in prose. These are provided as projects on ``communicating your mathematics'' \cite[Section 16]{chiodo2023teaching}, on ``open-ended discussion questions'' \cite[Section 17]{chiodo2023teaching}, and on ``longer essay topics'' \cite[Section 18]{chiodo2023teaching}, as follows:

\begin{enumerate}
    \item \textit{Communicating your mathematics} \cite[Section 16]{chiodo2023teaching} includes writing ``Ethics memos'' for mathematical work: short, succinct summaries of an ethical problem in mathematics, akin to a policy memo. It also includes a (group) project looking at ``How your science can end up in the newspapers'', analysing how the media portrays scientific and mathematical results.
    \item \textit{Open-ended discussion questions} \cite[Section 17]{chiodo2023teaching} includes ``Open-ended discussion questions'' on the ethical issues surrounding topics such as the Diffie-Hellman Key Exchange, the work of statistician Ronald Fisher and its motivations and applications, and the misuse of statistics in the wrongful conviction of Sally Clark.
    \item \textit{Longer essay topics} \cite[Section 18]{chiodo2023teaching} aims to help students develop skills such as using \LaTeX, looking up research papers, writing prose, looking into something with technical and non-technical aspects, and (ultimately) doing research on ethics in mathematics. It suggests topics such as algorithmic fairness, mathematicians and cryptography, the Ofqual grading algorithm, modelling COVID-19, and medical statistics and probability. It even invites students to find a topic of their choosing from current events.
\end{enumerate}
Tools such as these provide initial motivation, guidance, and training, for anyone wanting to actively think about, and write about, ethics in mathematics, in a meaningful and impactful way. This allows mathematicians to turn ethics in mathematics from a `spectator sport' (i.e., something they read about from afar and long after the fact) to a participatory activity (i.e., something they partake in and contribute to). After all, just as with their mathematical studies, there comes a time when every students needs to stop just reading about mathematics \textit{other} mathematicians have written, and start writing some of their own; in this regard, ethics in mathematics is no different.

\begingroup
\let\section\subsection
\renewcommand{\bibname}{Bibliography (Lecture 8)}
\renewcommand{\refname}{Bibliography (Lecture 8)}
\addcontentsline{toc}{subsection}{Bibliography (Lecture 8)}
\renewcommand{\addcontentsline}[3]{}
\bibliographystyle{IEEEtran}
\bibliography{bibLec8}
\endgroup

\newpage

\end{document}